\documentclass[10pt]{amsart}
\usepackage[T1]{fontenc}
\usepackage{amssymb}

\newcommand{\vsp}{\vspace*{0.5cm}}

\begin{document}
\begin{flushleft}

\title{On partially hypoelliptic operators. \\ Part II: Pseudo differential operators}
\author{T. Dahn \\ \emph{Lund University}}

\begin{abstract}
This study is an attempt at generalizing the class of partially hypoelliptic differential operators to a class of pseudodifferential operators, Symbol ideals are formed on the set of lineality and we discuss suitable topologies that allow the generalizations.
\end{abstract}
\maketitle

\section{Introduction}
The results in this study are an attempt to generalize results on formally partially hypoelliptic differential operators (cf. \cite{jag_0}) to more general operators. For a discussion on representations of pseudo differential operators we refer to (\cite{jag_0}). For a discussion on fundamental results on the class of differential operators, see (\cite{jag_I}).

\vsp

The generalization is based upon symbols in a (geometric) ideal $f \in(I)=(\mbox{ ker h })$, and $f(\zeta)=F(\gamma)(\zeta)$, where $\gamma$ is considered in a dynamical system (or pseudo base), $F$ is a lifting operator and $\zeta$ in a domain of holomorphy. The micro-local analysis is based upon the set of lineality and we will discuss generalizations of this set of invariance. We consider for instance the weighted lineality,
of points $\eta$, such that $\frac{Q(x+\eta)}{Q(x)} \tau_{\eta} h \tau_{\eta}=h$, for a self-adjoint and hypoelliptic polynomial $Q$ and where $\tau_{\eta}$ denotes translation. A longer discussion on this set appears in the end of this study.

\vsp

Allthrough the study we will use the concept of monotropy (cf. \cite{Cous}). A monotropic function is a continuous function $f$ such that there exists a holomorphic function with zero's within $\epsilon$ distance from the zero's to $f$. We consider monotropic ideals such that $h=\tau_{\epsilon}w$, where $w$ is algebraic in the sense that $w(x^{2})=w^{2}(x)$ and we write $h \sim_{m} w$. A longer discussion on the monotropic functions and functionals is given in (\cite{jag_0}).

\section{Symbols in a geometric ideal}
\subsection*{ Reduced operators }

Let's define a real set, called the set of lineality, as the zero's to a holomorphic function $\varphi_{\lambda}(i t \eta)=P_{\lambda}(\xi+i t \eta)-P_{\lambda}(\xi)$,
$$\Delta_{ \mathbf{C}}=\{ \eta \in
 \mathbf{R}^{\textit{m}}; \quad \textsl{P}(\xi + \textsl{it} \eta) - \textsl{P}(\xi)=\textsl{0} \quad  \forall \xi \in
 \mathbf{R}^{\nu}, \quad  \forall \textsl{t} \in  \mathbf{R} \}=\textsl{Z}_{\varphi_{\lambda}}$$
We can define a pseudoconvex neighborhood of $Z_{\varphi_{\lambda}}$ using
$ ( 1 + c\mid \zeta \mid )^k \mid \varphi_{\lambda} ( \zeta ) \mid^2 = C h(\zeta) $
for a holomorphic function $h \in  H(U')$, $U \subset U'$, for an open set $U'$, for positive constants $c,C$ and $k$.
Let $g( \zeta ) = C h( \zeta ) - ( 1 + \mid \zeta \mid )^k \mid \varphi_{\lambda} ( \zeta ) \mid^2$, be the holomorphic
function defining $U$. $P_{\lambda}$ denotes for $\lambda \in  \mathbf{C}$ the operator $P-\lambda$.
\newtheorem{b_Treves}{Lemma }[section]
\begin{b_Treves} \label{b_Treves}
If $P_{\lambda}$ is reduced for complex lineality (with respect to one dimensional translations), thus  $\Delta_{ \mathbf{C}}=\{ 0 \}$ ,
 then there are positive constants $c,C$, such that $(1+ \mid \xi \mid)^c \leq C \mid P_{\lambda}(\xi) \mid$ for all $\xi \in  \mathbf{R}^{\nu}$.
\end{b_Treves}

Assume $ \mathcal{H}_{\textit{c}}$ is the class of polynomials satisfying the inequality in Lemma \ref{b_Treves}, for some constant $c$.
For $h(\zeta)=0$, since $g^t$ is reduced, for some positive
integer $t$, there is a positive $\sigma$, such that $\mid \xi
\mid^{\sigma} \leq \mid g(\xi) \mid$, for large $ \xi \in $
$ \mathbf{R}^{\textit{m}}$. Thus, $\mid \xi \mid^{t \sigma -k''} \leq C
\mid \varphi_{\lambda}(\xi) \mid^{2t} $, for some positive $k''$.
For $\mid h(\zeta) \mid > 0$, the inequality for $h$, is
immediate. Assuming the defining polynomial is self-adjoint, we
conclude that $\varphi^t_{\lambda} \in  \mathcal{H}$$_{\frac{1}{2}(t
\sigma - k'')}( \mathbf{R}^{\textit{m}})$. This means, for a
sufficient number of iterations, $\sigma > k''/t$ and the polynomial is reduced. The constant $k'' < 1$, so there is a
positive integer $N$ (possibly smaller than $t$), such that
$\sigma > 1/N$ and $\varphi^N_{\lambda}$ is hypoelliptic in $L^2$.
Note that this hypoellipticity is not dependent on $\lambda$.

\newtheorem{b_Ruc}[b_Treves]{ Proposition }
\begin{b_Ruc} \label{b_Ruc}
If $P(D)$ is a constant coefficients, partially hypoelliptic operator, there is an iteration index
$N_0$, such that $P(D)^N$ is hypoelliptic in $ \mathcal{D}'$, for all $N \geq N_0$.
\end{b_Ruc}

We can prove a similar result for self-adjointness. The notation $P \prec Q$ indicates that
the quotient $ \mid (P/Q)(\xi) \mid$ is bounded as $\mid \xi \mid \rightarrow \infty$ $\xi$ real and $P \prec
\prec Q$ that the same quotient goes to zero as $\mid \xi \mid \rightarrow \infty$.
\newtheorem{b_ortho-strength}[b_Treves]{ Lemma }
\begin{b_ortho-strength} \label{b_ortho-strength}
Assume $P$ and $Q$ constant coefficients differential operators such that $P \prec Q$, $Q$
hypoelliptic and $(P\varphi,Q\varphi)=0$ for all $\varphi \in N(Q)^{\bot}$. Then $P \prec \prec Q$.
\end{b_ortho-strength}

\section{ Geometric symbol ideals }
We have the following problem: For which class of operators, do we have that the iterated symbols
define hypoelliptic operators?

\vsp
We will use notions and results from  \cite{Oka}. An ideal of holomorphic functions, in the sense of
 \cite{Oka}, is an ordered set of $(f,\delta)$, where $f$ is holomorphic in $\delta \subset
 \mathbf{C}^{\nu}$, that is closed for addition and multiplication with a holomorphic function, on the
intersection of the domains corresponding to the functions. Two ideals are said to be equivalent,
if for every point in the domain, any function in one ideal for this point, is in the other for the
point. A local pseudo-base for the ideal $(I)$, is a finite system $ \mathcal{K}$ of holomorphic functions $F_i \in (I)$
such that for every $f \in (I)$, $f \equiv 0 ( \mbox{ mod } \mathcal{K})$ in a neighborhood of a given
point. If an ideal has a pseudo-base, then any equivalent ideal has a pseudo-base. A geometric ideal,
is an ideal, defined as holomorphic functions and domains $(f,\delta)$, such that $f \equiv 0$ on
$\Sigma \cap \delta$, where $\Sigma$ is a given characteristic set (analytic set).

\newtheorem{b_CondHE}{Condition}[section]
\begin{b_CondHE}[HE] \label{b_(HE)}
A class of holomorphic functions $ \mathcal{K}$, is said to have the property $(HE)$, if it satisfies the criterion on
complex zero's for hypoellipticity, that is if $f \in (I)$ and $\zeta$ is such that
$f(\zeta)=0$ and $\mbox{ Im }\zeta$ is bounded, then $\mbox{ Re }\zeta$ is bounded.
\end{b_CondHE}
(cf. (\cite{Ho_LPDOII}))
\newtheorem{b_Cartan}[b_CondHE]{ Lemma }
\begin{b_Cartan}
Given $(I)$, an ideal of holomorphic functions with the property \ref{b_(HE)}, we have a local
pseudo-base.
\end{b_Cartan}
Proof:
The condition \ref{b_(HE)} defines a geometric ideal, why the Lemma follows from a theorem by H.
Cartan  \cite{Oka} Ch VIII, sec.3 $\Box$.

\vsp

We will in this section study only entire functions of finite type. The ideals can thus be studied
in the topology $\mbox{ Exp }_{\rho,A}( \mathbf{C}^{\nu})$ and we assume that $\rho$ is
chosen so that we have a product topology.
Any $C^{\infty}( \mathbf{R})$-function, real-valued and of exponentially finite (real) type, is the real part of
an analytic function, in fact such a function is a real-analytic function itself. The same result
holds in several variables for $\mbox{ Exp }A$ with product-topology, by iteration of the
one-dimensional result. Thus, given a real-valued function in $C^{\infty}( \mathbf{R}^{\nu})$ of (real) exponential
type $A$, we have $f=\mbox{ Re }F$ in $\mbox{ Exp }A$. 

\vsp

Assume $E(x,y)$ a kernel in $C^{\infty}(\Omega \times \Omega), \ \Omega \subset  \mathbf{R}^{\nu}$, that can be
represented, by the real part of a holomorphic function $f$ and thus has a representation in $C^{\infty}$ by a real
analytic function. Assume $F$ has complex zero's (in both variables separately ), that include the set of common complex zero's,
for an ideal $(I_{HE})$, defined by the property \ref{b_(HE)}. Assume $ \mathcal{K}$ a local pseudo-base of holomorphic functions
$F_1, \ldots, F_p$ for $(I_{HE})$. According to R\"uckert Nullstellensatz ( \cite{Oka})
$$ F^{\lambda} \equiv 0 (\mbox{ mod }  \mathcal{K})$$
for some real number $\lambda$. This means, that we
can determine holomorphic functions $A_1,\ldots, A_p$ in $H(\Omega \times \Omega)$, such that $F^{\lambda}=\sum^r_{j=1} A_jF_j$,
in $H(\Omega \times \Omega )$.

\vsp

As the regularizing operators are compact, we note that for reduced and regularizing operators the
square root of the operator is defined. More precisely,
since the operator ideal of compact operators is idempotent, if we consider a reduced $F \in H \cap L^2$ and
 $F=\big[F^{1/2},F^{1/2}\big]$, we have that $F^{1/2}$ corresponds to a compact operator and if we
 use that $H \cap L^2$ is nuclear, we see that $F^{1/2} \in H \cap L^2$.

\vsp

Assume $F_1,\ldots,F_p$ a pseudobase of hypoelliptic polynomials, $V=\cap_jZ_{F_j}$, then locally
$Z_{F_1} \cap\Omega_1=Z_{F_2} \cap \Omega_1$, for a domain of holomorphy $\Omega_1$. Since the pseudo-base
elements are assumed hypoelliptic, it is sufficient to consider the real part $P_j=\mbox{ Re }F_j$ and we have
$P_1=HP_2$ for instance with $H$ real analytic in $\Omega_1$. Using parametrices to these operators, adjusted
to $ \mathcal{E}'$, we see that $E_2=HE_1$, where $E_j$ are entire and of type $0$ for $j=1,2$. This means that,
$H$ is entire and of type $< 0$ (study $ I=HP_2E_1$) and is thus bounded and by Liouvilles' theorem, constant
( if the argument is divided into two inclusions, we get $H'=1/H''$ constant )
We conclude that $P_1 \sim P_2$ and thus $F_1 \sim F_2$. The same argument can be applied for any pseudo-base
elements. We can also use this argument to prove that it is sufficient to consider $(I_{HE})$ as generated by
hypoelliptic and equivalent ( in strength ) operators and over constant coefficients.

\vsp

A lacunary point for an ideal $(I)$, is a point such that for all functions $f \in (I)$, $f$ is not
in $(I)$ for this particular point. We have seen (cf. \cite{jag_I}) that for partially hypoelliptic, constant
coefficients polynomials, the set of lacunary points relative the corresponding ideal $(I_{HE})$, is a
(locally) finitely generated set.

\vsp

The symbol $E$, can now be represented as $$E=\frac{1}{2}\big[
b_{\Gamma_{+}}f + b_{\Gamma_{-}}\overline{f} \big]$$
Over the lineality, we have that $E=\mbox{ Re }f$ and outside this set we have that the difference is in $\mbox{
Exp }_{0,\parallel \cdot \parallel}$ and in $C^{\infty}$.
Particularly, if the symbol allows real support, the limit is well defined and uncomplicated.
The set of complex lineality, can be defined as before, but we will have to prove a proposition
analogous to $b_{\Gamma}$ being portable by $\Delta_{ \mathbf{C}}$.

\section{ Lineality for an ideal }
Assume the ideal $(I)$ defined by a homomorphism $h$, that is $(I)=(\mbox{ ker h})$.
Consider now the lineality , $\Delta(h)$, corresponding to $h$ with $h \sim_{\infty} \mbox{ Re }h$,
where $\sim_{\infty}$ refers to equivalence in strength (cf. \cite{jag_I}). If $h$ is injective, naturally
$\Delta(h(f))=\Delta(f)$ for all $f \in (I)$. Otherwise $\tau_{\eta}h(f)=h(f)$ is equivalent with
${}^t\tau_{\eta}f-f \in J_h$ and we define $\Delta_I(h)=\{ \eta \quad \tau_{\eta}f-f \in J_h \quad \text{ for
all  } f \in (I) \}$. We assume in this context, that $\eta$ defines a full line.
Assume further that $h$ maps real-analytical functions on real-analytical
functions. Using the theory on analytic functionals (cf. \cite{Mart} ), if $b_{\Gamma}h$ is
quasi-portable by $\Delta_I(h)$, then it is portable by $\Delta_I(h)$. Quasi-portability means particularly that
the limit $\lim_{\Delta_I(h) \ni \eta \rightarrow 0} \tau_{\eta}h(f)$ exists, independently of the
choice of $\eta \in \Delta_I(h)$ and assuming $\Delta_{I}(h)$ is analytic, we can conclude that $b_{\Gamma}h$ is portable by $\Delta_I(h)$.
Strictly speaking, if the ideal is over $H(U)$, for an open set $U$, the porteur is $U \times
\Delta_{I}$, but we will write the porteur as $\Delta_{I}$, since $b_{\Gamma}$ does not depend on
$U$.

\vspace*{.5cm}

For a reduced $F$, algebraic over $\tau_{\eta}$, the equality $\tau_{\eta}(BF)=\tau_{\eta}(B)\tau_{\eta}(F)=BF$ should be treated as
$\tau_{\eta}B\tau_{\eta}=B$. Assume $B=I_{B'}$, that is given by an integral operator with kernel
$B'$, then $\tau_{\eta}B=B\tau_{-\eta}$. We assume also, $B'(x+\eta,y-\eta)=B'(x,y)$, that is
$\eta$ is in the lineality for the kernel, separately in each variable $\eta \in \Delta_{x,y}(B')$.
Particularly, if $\eta \in \Delta(h)$ and if $\mbox{ ker }\sigma_h$ gives the coefficients to
$h(f)=0$ on the form of integral operators $B$ with kernel $B'$. Then, for all $B \in \mbox{ ker
}\sigma_h$, we have $B'(x+\eta,y-\eta)=B'(x,y)$.

\vspace*{.5cm}

Assume $f$ is not reduced, but partially hypoelliptic and $h(f)$ with the same properties.
$h(f)=\sum_jA_jF_j$ with $A_j \in H$ and $F_j$ reduced and $f=\sum_j B_jF_j$ with $B_j \in H$.
Finally, $h(F_j)=\alpha_jF_j$ for constants $\alpha_j$ and all $j$. Then, $\eta \in \Delta(h(f))$
$\Leftarrow$ $\eta \in \cap_j \Delta(A_j)$. If further $h(B_j)=\beta_jB_j$ for constants $\beta_j$ and
all $j$, then $\Delta(h(f))=\Delta(f)$. Otherwise, $\frac{h(B_j)}{B_j}=\beta_j \in H$ not all
constant, but with kernel $\beta_j'$. Thus, $\eta \in \Delta(h(B_j))$ implies that $\beta_j'$ are invariant
in both variables separately for translation with $\eta$. Finally,
$$ \eta \in  \Delta_{x,y} \big( \frac{h(B_j)}{B_j} \big) \Leftrightarrow \tau_{\eta}B_j - B_j \in
J_h \quad \forall j$$

\section{ Schwartz-type topology }
We make the following proposition, assume $(I)$ a symmetric ideal of holomorphy with topology of Schwartz-type (cf. \cite{Mart})
and a compact translation.

\newtheorem{e_minimal}{ Proposition }[section]
\begin{e_minimal} \label{e_minimal}
If $\psi \in (I)$ and $\psi \sim_{m} 0$ in $\mid \zeta \mid$-infinity. Then we have that $\{ d_{\zeta}
\psi=\psi=0 \}$ is nowhere dense in $N(I)$.
\end{e_minimal}

Proof:\\
Assume $\psi$ is approximated by reduced $\psi_{j} \in (I)$ such that
$d_{\zeta}\psi_{j}
\neq 0$. Let $X=N(I)$ and assume $\psi$ defined on an open, small, complex set $\Omega \supset X$.
Assume further that $V_{j}=\{ \zeta \in \Omega \quad \psi_{j}(\zeta)=\psi(\zeta) \}$ and that as
$\psi_{j} \rightarrow \psi$, we have $V_{j} \uparrow \Omega$. In case $V_{j}$ is connected, we have
an infinite zero on $V_{j}$ for $\psi_{j}-\psi$. Through the condition $\psi \sim_{m} 0$, if the
$\mid \zeta \mid$-infinity lies in a connected component in $\mbox{ lim }V_{j}$, then
$d_{\zeta}\psi_{j} \sim_{m} d_{\zeta} \psi$, why the zero for $\psi$ is simple in the $\mid \zeta
\mid$- infinity.$\Box$

\vspace*{.5cm}

Assume $\psi(x)=\tau_{x}\phi/\phi$, where $\phi$ is considered on the real space.
Let $\Gamma=\{ x \quad \psi(x)=1 \}$ and let $\Gamma_{j}$ be the corresponding sets for the
approximating sequence $\psi_{j}$. As $\psi_{j}$ are reduced, the sets $\Gamma_{j}$ must be
isolated points and as $\psi_{j} \rightarrow \psi$, we must have $\Gamma_{j} \rightarrow \Gamma$,
why if $\mbox{ lim }V_{j}$ connected, the set $\Gamma$ must be a point.

\vspace*{.5cm}

Note that a condition $\psi \rightarrow 1$ in the $\mid \zeta \mid$-infinity should be compared
with the condition of slow oscillation, when $\mbox{ Im }\psi \rightarrow 0$, why we always assume
$(I)$ symmetric. Note that $(I)$ with compact translation can be seen as radical and thus as
finitely generated, why we expect the minimally defined situation.

\newtheorem{e_red}[e_minimal]{ Lemma }
\begin{e_red}
Assume existence of $\{ \psi_{j} \}$ such that $\psi_{j} \rightarrow \psi$ and $\psi_{j}$ reduced
in $(I)$ and such that $d_{\zeta}\psi_{j} \neq 0$ $\forall j$. Let $\Gamma=\{ x \quad \psi=1 \}$.
Then we have that $\Gamma$ is an isolated point.
\end{e_red}

\newtheorem{e_appr}[e_minimal]{ Lemma }
\begin{e_appr}
There is an approximating sequence $\{ \psi_{j} \}$ according to the previous Lemma.
\end{e_appr}

The second Lemma is immediate, since $(I_{RED}) \subset (I)$ and if the condition $\psi \sim_{m} 0$
is satisfied for $\psi$ such that $\mbox{ Re }\psi \rightarrow 1$ and $\mbox{ Im }\psi \rightarrow
0$ as $\mid \zeta \mid \rightarrow \infty$, that is we are in the proposition considering the set
$\{ d(\psi-1)=(\psi-1)=0 \}$, then the approximating sequence exists as the
translation is assumed compact. Considering the first Lemma, if we do not have a connected
component in the limit-point, there is at least one point in $\Omega$ where all the jets for
$\mbox{ lim }\psi_{j}$ and $\psi$ coincide. Assume now $d_{\zeta} \psi_{j}- \psi_{1} \rightarrow 0$,
with $d_{\zeta} \psi_{j} \neq 0$ and $\psi_{1} \sim_{m} d_{\zeta} \psi$. We then have modulo
monotropy that $\psi_{1} \neq 0$ in the limit point, why $\psi$ has a simple zero there. In the same
manner we see that $\Gamma$ is a point ( modulo monotropy ). Finally, we obviously have $\{ \psi_{1}=\psi=0 \}$
is nowhere dense in $X$, why the proposition follows modulo monotropy.

\vspace*{.5cm}

Consider the case with topology of Schwartz-type and weakly compact translation. More precisely,
consider $\mbox{ rad }(I)$ and $(I)$ is assumed as above with compact translation. If $(I)=\mbox{ ker }w$
for a homomorphism $w$, that we assume algebraic modulo monotropy, that is $w$ is such that $w^{N}$
is locally injective. According to proposition \ref{e_minimal}, if $\psi \sim_{m} 0$ and $w(\psi)=0$, we
have that $\{ d_{\zeta}w(\psi)=w(\psi)=0 \}$ is nowhere dense in $N(I)$. Further $\psi^{N} \in (I)$
means that $\psi \sim_{m} 0$ and $\{ d_{\zeta} w(\psi^{N})=w(\psi^{N})=0 \}$ is nowhere dense in
$N(I)$ and for a suitable $N$ we see that $\{ d_{\zeta} \psi=\psi=0 \}$ is nowhere dense in $N(I)$,
where the argument is modulo monotropy.

\section{ Lacunary points }
There will be significant lacunary cases in this approach. For instance any case, where
$b_{\Gamma}$ is not portable by the lineality, is lacunary for the symbolclass.
More precisely, a lacunary case, is when
the symbols $g_{(\alpha)}(i t \eta)$, in the development of $g(\zeta)-g(x)$, are
only locally constant in a neighborhood of $0$. We study neighborhoods of $\Gamma_{g}$, such that there
exists constants $C_{t,\epsilon}$ $$
\Lambda_{t,\epsilon}=\{ \eta \ e^{-\epsilon \parallel x \parallel} \mid g(x+i t \eta)-g(x) \mid \leq
C_{t,\epsilon} \}$$ where $g(\zeta)-g(x)$ has real type $0$. Obviously, for any finite $t$, these sets
consist of all $\eta$. We can find a-priori-estimates in
 \cite{Mart},\\
 $ \log \mid g_{(\alpha)}(i t \eta)- g_{(\alpha)}(0) \mid \leq \gamma_{\alpha}(i t \eta)$
where
$$
\gamma_{\alpha}(i t \eta) = \frac{1}{( \pi \rho^2)^{n-k}} \int^{\rho}_0 \int^{2 \pi}_0\!\!\! \ldots \int^{2 \pi}_0 \log \mid g_{(\alpha)}(i t \eta + r e^{i \theta}) - g_{(\alpha)}(re^{i \theta}) \mid
\times$$
$$\times d \theta_{k+1} \ldots \theta_{\nu} d r_{k+1} \ldots d r_{\nu}$$
Immediately, $\gamma_{\alpha}(0)=0$ and $ \mid \gamma_{\alpha}(i t \eta) \mid \leq \epsilon t$,
for some small positive $\epsilon$. Further, $\mid \gamma_{\alpha}( i t \eta) \mid \leq
\frac{t}{\rho} \log M_{\alpha}(t + \rho)$, where $M_{\alpha}(R)=\sup_{\parallel t \eta
\parallel_{\infty} \leq R} \mid g_{(\alpha)}(i t \eta)-g_{(\alpha)}(0) \mid$ and as usual
$M_{\alpha}(R)^{\frac{1}{\mid \alpha \mid}} \rightarrow 0$ as $\mid \alpha \mid \rightarrow
\infty$, for all $R$. Finally, we have the following result (cf.  \cite{Horm} Cor. 4.4.14). The Lelong number
$v_{\varphi}(\zeta)$ is defined for a pluri-subharmonic function $\varphi$ on an open set $\Omega
\subset  \mathbf{C}^{\nu}$ as $$v_{\varphi}(\zeta)=\lim_{r \rightarrow 0} \int \frac{d \mu(\zeta +
w)}{c_{2n-2}r^{2n-2}}, \text{ where } d \mu= \Delta \varphi / 2 \pi.$$ Let $V= i \Gamma_{g} \cap
\Omega$ and $\varphi=\frac{1}{2} \log \mid g(\zeta)-g(x) \mid^2$, then $v_{\varphi}(\zeta)=0$ for $
\zeta
\in \Omega \backslash V$ and $v_{\varphi}(\zeta)$ on $V$ gives the order of zero in $\zeta$.

\vsp

Note that in the case where $g_{(\alpha)}$ is only locally constant, we consider support-sets for
$\widehat{g}_{(\alpha)}$ for $ x \in  \mathbf{R}^{\textit{m}}$ on the form $x \in C^*_{[t_0,t_1]}$ if
and only if $t x \in  C^*_{[t_0,t_1]}$ for all $t$ such that $t_0 \leq t \leq t_1$. Then obviously,
$\widehat{g}_{(\alpha)} \in L^1( C^*_{[t_0,t_1]})$. Further, for $\eta \in \Delta_{t,\epsilon}$
$$ \mid g_{(\alpha)}(-i t \eta) \mid \leq \mid \int_{C^*_{[t_0,t_1]}} e^{- i < x, i t \eta >}
\widehat{g}_{(\alpha)}(x) d x \mid \leq e^{H_{t_0,t_1}(\eta)} \parallel \widehat{g}_{(\alpha)}
\parallel_{L^1}$$
and using the norm for locally summable Fourier-transforms, $M_{\beta}$
\begin{equation} \label{b_reg-cone}
\mid D^{\beta} g_{(\alpha)}(-i t \eta) \mid \leq e^{H_{t_0,t_1}(\eta)} M_{\beta}(g_{(\alpha)})
\end{equation}
where the last expression is well-defined, assuming that $g_{(\alpha)}$ is not constant on the support-set
that we are considering. Particularly, using Cauchy's inequalities
$$ \mid g_{(\alpha)}(i t \eta) \mid \leq M t^{-\mid \alpha \mid}$$
why $\mid g_{(\alpha)}(i t \eta) \mid \rightarrow 0$ as $\mid \alpha \mid
\rightarrow \infty$ and as $t \rightarrow \infty$. As $t_1 \rightarrow \infty$ in (\ref{b_reg-cone}), we see
that $g_{(\alpha)} \rightarrow 0$, so we can consider support-sets that are actually regular cones.
Conversely, we have to consider $\alpha$ larger than the Lelong-number on the support-cones, to avoid degeneracy
in (\ref{b_reg-cone}). We get expressions like
$$ \mid D^{\beta} ( g(\overline{\zeta})-g(x) ) \mid \leq C \sum_{d}^{\infty} e^{H_{t_0}(\eta)}
M_{\beta}(F_{(\alpha)}) \mid x^{\alpha} \mid$$ for $F_{(\alpha)}=(g_{(\alpha)}(i t
\eta)-g_{(\alpha)}(0))$.

\vsp

We say that $x'$ preserves a constant value, if $F(x',y)$ is analytic for $x'$ inner to a closed contour
(assume product of contours) and $y$ is outer to another closed contour, and if $F(x',y)=f(x',y)-f(x',x_n) \rightarrow 0$ as $y
\rightarrow \infty$. Assume for instance $C,C',C_1,C_1'$ simple, closed contours in the respective
planes $x',x_n$ such that $C'$ inner to $C_1'$ and $C$ outer to $C_1$. Assume $F$ analytic for $x'$
inner to $C_1'$ and $y$ outer to $C_1$, $a$ a point outer to $C_1'$ and $b$ a point inner to $C_1$.
Then according to Cousin ( \cite{Cous} \S24, Prop. 4), for any small positive $\epsilon$,
there is a polynomial in $(x'-a)^{-1},(y-b)^{-1}$, such that for $x'$ inner to $C$, $y'$ outer to
$C'$,
$$ \mid F(x',y)-Q(\frac{1}{x'-a},\frac{1}{y-b}) \mid < \epsilon $$
and we can produce estimates, for $\eta \in \Delta_{t,\epsilon}$ as $\mid f(x',x_n+i t \eta)-f(x',x_n) \mid < \epsilon + \sum_{\alpha}
A_{\alpha} t^{-\alpha}$, for constants $A_{\alpha}$ and $x'$ inner to $C$.

\subsection{ Lacunary points for the induced topology }

Inducing a topology on the geometric ideal of operator symbols will affect the lacunary points for
this ideal. In an attempt to understand this problem, we will modify the concept of a quasi-portable
functional limit.
$T$ is said to be quasi-portable ( with respect to $X$) by an open set $\Gamma \subset \Omega$, if there exists a $\mu_T
\in H'_{\Gamma}(\Omega)$ such that $T=i_{\Gamma,\Omega}(\mu_T)$ where $i_{\Gamma,\Omega}: (X)'(\Gamma) \rightarrow (X)'(\Omega)$. Let $i^*:(X)(\Omega) \rightarrow (X)(\Gamma)$. If a measure $\mu$ with respect to $X$,
is portable with respect to $X$ by $\Gamma$, then $$\mid <\mu_{x},\phi> \mid \leq C \parallel \phi \parallel_{X(\Gamma)}, \text{ for
all } \phi \in (X)(\Omega),$$ why in this case the measure is also strictly portable
(with respect to $X$) by $\Gamma$. Assume more generally $\Gamma_{ \mathbf{R}}$ a real, convex
set. According to the theory of analytic functionals, if
$\Delta_{\Gamma}$ is the indicator to $\Gamma_{ \mathbf{C}}$ and $h_F(\eta)=\limsup_{t
\rightarrow \infty} \frac{1}{t} \log \mid F(t \eta) \mid$, for $\eta$ real close to $\Gamma$, $\mid \eta
\mid=1$ and $F= \mathcal{F}$ $T$ (the Fourier-Borel transform) to $T \in H'$ and for
$\Gamma_{ \mathbf{C}}=\Gamma_{ \mathbf{R}}+i{\Gamma}_{ \mathbf{R}}$ which we assume closed,
$(i \Gamma_{ \mathbf{R}})^{\circ}=(i \Gamma_{ \mathbf{R}})^{\bot}$.
Assume $\mu$ a measure with respect to $X=H \cap L^2$, that is
\begin{equation} \label{b_meas-def}
 \mid \mu(\phi)(x) \mid \leq C \parallel \phi \parallel, \quad \forall \phi \in X(\Omega)
\end{equation}
$\Omega$ an open set in $ \mathbf{R}^{\nu}$ containing $x$. Let $\Gamma_y=\{ ty, \ t \in  \mathbf{R} \}$ and $$\Gamma_C=\{ y, \
\mid \mu(\phi)(x+i\Gamma_y)) \mid \leq C \parallel \phi \parallel \quad \forall \phi \in X(\Omega+i\Gamma_y) \}.$$ We can assume $\Omega= \mathbf{R}^{\nu}$ and
$T^{\Gamma_y}=\Omega+i\Gamma_y$. Over $\Gamma_C$ we obviously have
\begin{equation} \label{b_meas-limes}
\lim_{t \rightarrow 0} \mu(\phi)(x+ity)=\mu(\lim_{t \rightarrow 0} \phi)(x+ity)
\end{equation}
and we claim that $b_{\Gamma}\mu$ can be considered as quasi-portable by $\Gamma$, the union over all
finite constants $C$ in (\ref{b_meas-def}). The limits (\ref{b_meas-limes}) are trivial over the
lineality-set, that is $\Delta_{\phi}=\{ y; \mu(\phi)(x+ity)=\mu(\phi)(x) \forall t \ \forall x
\}$. Outside $\Gamma_C$ in a conical neighborhood, $\Upsilon \subset \Gamma$, since the limit is trivial in $X(T^{\Upsilon})$, we have $$\mid \lim_{t \rightarrow 0} \mu(\phi)(x+ity') \mid \leq C' \parallel
\phi \parallel \quad \forall \phi \in X(T^{\Upsilon}),$$ why the limit on the left side of
(\ref{b_meas-limes}) exists, for a new constant $C'$. Assume finally, that $\mu$ is reduced in $L^2$,
then for all $\phi \in X(T^{\Upsilon})$ and $(\tau_{y'}\mu)(x)=\mu(x+ity')$ for an $y' \in \Upsilon
\backslash \Gamma_C$, $\mid (\tau_{y'}\mu)(\phi)(x) \mid \leq C'' \parallel \phi \parallel$ and
$\phi \in X(T^{\Upsilon})$, where we have used that the translation will be bounded in
$X$. We see that to any $\Upsilon \subset \Gamma$, there is a translated measure, defining the same
functional limit $b_{\Gamma}\mu$.

\newtheorem{b_mes}{ Proposition }[section]
\begin{b_mes}
Assume $\mu$ a measure with respect to $L^2 \cap H$, then $b_{\Gamma}\mu$ is quasi-portable by the
lineality for $\mu$. If $\mu$ is a measure as above, such that also $\widehat{\mu}$ is a measure
with respect to $L^2 \cap H$, we have for  a $\phi \in L^2 \cap H$ and $ \mathcal{F}\textsl{T}=\widehat{\mu}(\phi)$ $(x)$, $T$ is semiportable by
$\Gamma_{ \mathbf{C}}=\Gamma_{ \mathbf{R}} + i \Gamma_{ \mathbf{R}}$, the lineality to $\widehat{\mu}(\phi)$ when $\parallel \phi \parallel \leq 1$.
\end{b_mes}
If for a closed, convex set $\Gamma$, the functional $T$ is portable by every convex neighborhood
of $\Gamma$, then $T$ is said to be semi-portable by $\Gamma$.

\newtheorem{b_semi-p}[b_mes]{ Lemma }
\begin{b_semi-p}
 $T$ is semi-portable by $\Gamma_{ \mathbf{C}}$ if and only if $h_F(x+ity)=0$ for $iy \in
(i\Gamma_{ \mathbf{R}})^{\circ}$ or if and only if $h_F(x+ty)=0$ for
$y \in (\Gamma_{ \mathbf{R}})^{\circ}=(\Gamma_{ \mathbf{R}})^{\bot}$.
\end{b_semi-p}
 It is thus, for our application ${X'}^{(0)}$
of no significance if we consider real or complex translations.

\vsp

Assume now $F=\lim_{t \rightarrow 0} \big[ F(x+ty)-F(x) \big]$, where
$F$ is for instance $\widehat{\mu}(\phi)(x)$, for a measure $\mu$, such that $\widehat{\mu} \in {X'}^{(0)}$,
Assume $\Gamma$ the lineality to $\widehat{\mu}(\phi)$, if $\parallel \phi \parallel
\leq 1$, we have $h_F(x+ty)=0$, for $y \in \Gamma$, for all $x$, why $T$ is semi-portable
(with respect to $X$) by $\Gamma_{ \mathbf{C}}$. If $h_F(\eta) \leq 1$ for all $\eta \in
\delta \Gamma$, where $\delta \Gamma=\{ y^* \ \mbox{ Re }<y^*,y> = 1 \ y \in \Gamma \}$, then $T$ is portable
by $\Gamma$ and $h_F(\eta) \leq \Delta_{\Gamma}$ for all real $\eta$.
Finally, if $i^*$ is defined as before with
$\Omega= \mathbf{R}^{\nu}$ and $\mu \in  \mathcal{C}'( \mathbf{C}^{\nu})$ is such that
$\widehat{\mu} \in X$ and $\mu=\mu_1+\mu_2$ with $\parallel i^* \widehat{\mu}
\parallel_{X( \mathbf{R^{\nu}})}=0$. We then have that $\parallel \widehat{\mu}
\parallel_{X( \mathbf{R}^{\nu})}=\text{ constant} (=1)$ implies existence of a $\widehat{T} \in
(X)'(\Gamma_{ \mathbf{C}})$, such that $i^*_{ \mathbf{C}}(\widehat{\mu})=\widehat{T}$. Further,
$$\mid < \widehat{\mu},\phi > \mid \leq C \parallel \widehat{\mu}
\parallel_{X( \mathbf{R}^{\nu})} \parallel \phi \parallel_{X( \mathbf{R}^{\nu})} \leq C'
\parallel \phi \parallel_{X( \mathbf{R}^{\nu})}.$$ Thus, $\widehat{\mu}$ is a measure with
respect to $X$, strictly portable by $\Gamma_{ \mathbf{C}}$.

\section{ A global base for partially hypoelliptic operators }

Define for $F_j \in  \mathcal{K}$ $j=1, \ldots,p$, a set $$G_{\Omega}=\{ (x, F_j(x)) \ \mid d_{\Gamma}f \mid \leq 1 \}$$
where, $d_{\Gamma}f=b_{\Gamma}f-f$, can be defined on real or complex arguments. Define an ideal
$(I_{const})$ over $G_{\Omega}$, as $g=b_{\Gamma}f-\sum_jA_jF_j$, that is $b_{\Gamma}f \equiv g ( \mbox{
mod } \mathcal{K})$ (assume for now that we are dealing with entire functions) and where $const$ is referring
to constant coefficients $A_{j}$ in the pseudo-base development. The zero's for $(I_{const})$
are given as points such that $b_{\Gamma}f \equiv 0 (\mbox{ mod }  \mathcal{K} )$. We then have that
$G_{\Omega}$ can be approximated from the outside, by polynomially defined sets. Further, $f \equiv d_{\Gamma}f (\mbox{ mod } \mathcal{K})$
where according to Lemma \ref{b_e-type}, $d_{\Gamma}f$ is of type $0$. Over the set of
lineality, we have $f \equiv g ( \mbox{ mod } \mathcal{K} )$. Further, it is possible to construct a
local pseudo-base for the projection of $(I_{const})$ on the first variable in $G_{\Omega}$. If $f$ is of type $A$ and $g$ is
developed in a local pseudo-base for $(I_{const})$, according to the remainder theorem, it can be given as a
polynomial in $F_j$'s with holomorphic coefficients of type $A$.

\vsp

Consider, for this situation, lacunary points as points such that $d_{\Gamma}f \notin (I_{HE})^Â$
Note that if $f \in (I_{HE})^Â$ is of type $A > 0$, then $f$ is reduced. That is if
$f^{(\alpha)}(it\eta)$ is constant $\neq 0$, for all $\alpha$, then $\eta =0$. Since assuming $\eta
\neq 0$, means \\ $\sum_{\beta \neq 0}$ $c_{\beta}(f^{(\alpha \beta)}(0))t^{\beta}$  $\neq 0$, for some
$\beta$, thus $f^{(\alpha)}$ is not a constant $\neq 0$. If on the other hand $f^{(\alpha
\beta)}(0)=0$, for all $\beta$, the type is $0$!.

\vsp

The Cousin integral ( \cite{Oka}) is defined as
$$ \Phi(x,y)=\frac{1}{2 \pi i } \int_L \frac{f(x,z)}{z-y} d z $$
where $x$ is assumed as inner to $\gamma=\gamma_1 \times \ldots \times \gamma_n$, where $\gamma_j$ are
closed contours in respective plane, where $L$ is a line containing the point $y_0$ in the $y$-plane
(one-dimensional). $L$ is assumed inner for a closed contour $\Gamma$ and $f$ is assumed analytic
(regular) for $(x,y)$ inner to $(\gamma,\Gamma)$. Assume $l$ a line through $y_0$, which does not
intersect $L$ in any other point and which is inner for $\Gamma$. Then $\Phi$ is analytic (regular) for $x$
inner to $\gamma$ and $y$ not on $L$, but otherwise arbitrary. The analytic continuation
$\widetilde{\Phi}$ can be shown to be analytic in $y_0$ and on $l$, on one side of $L$.
$\widetilde{\Phi}=\Phi \pm f$ (cf.  \cite{Cous}).

\newtheorem{b_lin}{ Proposition }[section]
\begin{b_lin}
Assume $\Delta_{ \mathbf{C}}$ the lineality set and $r(I_{HE})^{A}$ a radical ideal of holomorphy taken over
entire functions of type $A > 0$. Then there exists an ideal of holomorphy $J_A$, with the zero's given
by $\Delta_{ \mathbf{C}}$, such that $J_{A} \sim r(I_{HE})^{A}$.
\end{b_lin}

Proof: Assume $g \in J_A$, since this is a geometric ideal, we have $g=\sum_j c_j\phi_j$ and we
assume that in a sufficiently small neighborhood of $\Delta_{ \mathbf{C}}$, $\phi_j \in r(I_{HE})^Â$, for all $j$.
Thus $g^t \in (I_{HE})^Â$, for large $t$. Conversely, if $g \in r(I_{HE})^Â$ and $g=0$ on
$\Delta_{ \mathbf{C}}$, then $g^t=0$ on $\Delta_{ \mathbf{C}}$ and $g^t=\sum_i a_iF_i$ in a neighborhood of
$\Delta_{ \mathbf{C}}$. We may assume $F_i=0$ on $\Delta_{ \mathbf{C}}$, for all $i$, which gives a development of $g^t$ in
the pseudobase for $J_A$, which is a radical ideal.$\Box$

\vsp

Assume $g$ entire and in $J_A$ and that $g \rightarrow 0$ on a line $L_1$. Assume $L_2$ another
line in the same coordinate-plane, it is then a corollary to Lindel\"of's theorem, that $g
\rightarrow 0$ on the sector between and including $L_1$ and $L_2$. That is, since $g$ is of type $0$,
there is a $F$ analytic for $1/ \mid z \mid > \epsilon$ such that in this region, $zg(z)=F(\frac{1}{z})$
and the corollary can be applied on $zg(z)$. The corresponding Cousin-integrals
$$ \varphi_j(x,y)=\frac{1}{2 \pi i}\int_{L_j} \frac{g(x,z)}{z-y} d z \qquad j=1,2$$
have the property that $\varphi_1+\varphi_2=0$. The argument can be repeated, so that $\Phi=\sum_p
\varphi_p=0$ describes the entire $\Delta_{ \mathbf{C}}$. According to Cousin (cf.  \cite{Cous}) there is a "global"
function which coincides with $\varphi_p$ on the corresponding inner sector for each $p$ and such
that $\widetilde{\Phi}-\Phi$ is analytic. We claim that the second
Cousin-problem is solvable in this situation and that the ideal $J_A$ has a global pseudobase.
Particularly, if $\Delta_{ \mathbf{C}}=\{ 0 \}$, we have a global pseudo-base for $(I_{HE})^Â$.

\section{ The R\"uckert Nullstellensatz }

We have established that $\mbox{ Re }P \in  \mathcal{H}_{+}$, where $ \mathcal{H}_{+}$ denotes $ \mathcal{H}_{\textit{c}}$ for positive $c$
 and $\mbox{ Im }P \prec \prec P$ implies
 $P \in  \mathcal{H}_{+}$, but the converse does not hold, that is $P \in  \mathcal{H}_{+}$ does not imply
  that $\mbox{ Im }P \prec \prec P$, more precisely $I_{+}^* \subset I_{+} \subset r(I_{+})$.
    According to the real version of Nullstellensatz, the given ideal $I$ is radical if and only if
     $P_1^2 + \ldots P_k^2 \in I$ implies $P_i \in I$, for all $i$. This is consistent with our results on
      hypoellipticity for self-adjoint polynomial-operators in $L^2( \mathbf{R}^{\nu})$. Particularly, for a radical
       ideal of real symbols, locally $ \mid f \mid^2 \in I \Rightarrow \mbox{ Re }f,\mbox{ Im }f \in I$ and in our case,
        the ideal is locally given by real zero's. If for an ideal $I$, we let $I^*=\{ \mbox{ Re }f \ f \in I \}$
        and $I_{ \mathbf{C}}=I^*+iI^*$ the "complexified" ideal, it not difficult to see that
        $I_{ \mathbf{C}}=I$ if and only if $I$ is radical. For symbols as entire functions, we give the following result

\newtheorem{b_Ruckert}{ Proposition }[section]
\begin{b_Ruckert}
  Assume $F=(F_1,\ldots,F_p)$ a global pseudo-base for $I_{HE}^{*,A}$, the ideal of hypoelliptic symbols, self-adjoint
   and of type $A$ and $V=\cap_j N(F_j)$. If $f$ is an entire function of type $A$, such that $f=0$ on $V$,then according to R\"uckert's Nullstellensatz, we can find constants
    $A_1,\ldots,A_p$, such that
  $$ (\mbox{ Im }f )^{\lambda} = \sum_j A_j F_j \quad \text{ in } \Omega$$
  for $\Omega$ a holomorphically convex neighborhood of $V$ and for $\lambda > 0$ real.
\end{b_Ruckert}

In the polynomial case, assuming $\mbox{ Re }f \in I_{HE}^{*,A}$ and $\lambda > 1$, for $f$ as in the proposition, we have
$f \in I_{HE}^{A}$. If $I_{PHE}^{*,A}=r(I_{HE}^{*,A})$ then this must be a radical ideal, and it
is locally defined by its real zero's.

\section{ The Whitney closure }

Assume, in analogy with Kohn's approach (\cite{Kohn}), $\alpha(T)=\frac{\psi(\zeta +
T)}{\psi(\zeta)}$. For an analytic homomorphism $b$, $b(\psi)(\zeta +T)=\frac{1}{\alpha(T)}b(\psi)(\zeta)+C_{T}$.
If $C_{T}=0$ we know that $\alpha$ is holomorphic and non-constant on a bounded set, symmetric
around the origin. Let $\Sigma(b(\psi))=\{ T \quad b(\psi)(\zeta + T)=b(\psi)(\zeta)+C_{T} \}$.
Define the Whitney-closure, $\overline{(I_{C_{T}=0})}=\{\psi \quad b(\psi)(\zeta +
T)=b(\psi)(\zeta)+C_{T} \}$.

\newtheorem{f_J}{ Proposition }[section]
\begin{f_J} \label{f_J}
There is an ideal $(J)$ such that $rad(J) \sim \overline{(I_{C_{T}=0})}$
\end{f_J}
Let $\Omega_{0}=\{T \quad b(\psi)(\zeta + T)=b(\psi)(\zeta) \}$ with $\Omega^{(j)}_{0}=\{ T
\quad b(\psi^{j})(\zeta + T)=b(\psi^{j})(\zeta) \}$. We can prove that $\Omega^{(j)}_{0} \downarrow
\{ 0 \}$ as $j \uparrow \infty$. Consider $\psi=0$ on $\Omega^{(j)}_{0}$, then $\psi \in
I(\Omega^{(j)}_{0}))$, where $\Omega^{(j)}_{0}=N(J)$ for some $J$ and $\psi \in rad(J)$. Add a constant, such
that $\psi-c=0$ on $\Omega^{(j)}_{0}$, this gives the result ( modulo monotropy ).
We will now prove a correspondent to regular approximations of periodic points.

\newtheorem{f_per}[f_J]{ Proposition }
\begin{f_per}
Consider the
sets $V=\{T \quad B(\zeta + T)=B(\zeta) \}$ and $V'=\{T \quad DB(\zeta + T)=DB(\zeta) \}$,
where $D=D_{\zeta_{1}} \ldots D_{\zeta_{n}}$. In the finitely generated case, we have that $V \cap V'$
nowhere dense in $V'$. This means that there are $T_{j} \in \Sigma(B)$ such that $B(\zeta +
T_{j})=B(\zeta)+c_{j}$ and $0 \neq c_{j} \rightarrow 0$ as $T_{j} \rightarrow T_{0}$ and $T_{0} \in V$.
\end{f_per}

Proof:\\
Assume $(J)$ with Schwartz-type topology and with a weakly compact translation, in the sense that
we assume $\tau_{T}$ is compact over $(J)$ and weakly compact over $rad(J)$. Let $B=h(f)$, $f \in
(J)$ and reduced, then modulo monotropy, $B^{2} \sim_{m} h^{2}(f)$ and if $h^{N}=1$ for $N \geq
N_{0}$ we see that $\tau_{T}$ is compact over an iterate of $B$. Further $\tau_{T} B=h(\tau_{T}'
f)$, where $\tau_{T}'=\tau_{T'}$ is also a translation. For $T \in V'$ we have that
$B(\zeta + T) - B(\zeta)=c_{T}$ for all $\zeta \in \mathcal{B}$, a bounded set in $\Omega$ pseudoconvex and where the constant $c_{T}$ is dependent on
$T$. Assume further $T \rightarrow T_{0} \in V$. Our proposition is that there is a sequence of
$T_{j}$ in $V'$ such that $T_{j} \rightarrow T_{0}$. For $h$ algebraic, the result follows from
the conditions on $f$ and proposition \ref{f_J}. For $h=\tau_{\epsilon}w$ and $w$ algebraic, there
is a $T_{j}$, $\tau_{T_{j}} B^{2}-B^{2}=C_{T_{j}} \rightarrow 0$ as $T_{j} \rightarrow T_{0}$.
Modulo monotropy, this means that $h(\tau_{T_{j}} B-B) \sim_{m} C_{T_{j}} \rightarrow 0$.
The condition $T_{j} \in V'$ means that $\tau_{T_{j}} B-B=C_{T_{j}}'$ and we have $h(C_{T_{j}}') \sim_{m}
C_{T_{j}} \rightarrow 0$ as above.$\Box$

\section{ A monotropy-ideal }
\label{sec:c_mono}

Cousin \cite{Cous} uses the concept of monotropy as a weak sufficiency condition for solving the
Cousin-problems. We will develop some results for pseudo-differential operators using this concept.
Assume $A$ a set of isolated points, for a continuous function $g$, then there is a holomorphic function, with zero's on
$A \cap U$, for an open set $U$. Let $A_{\epsilon}$ be a set such that $P \in A_{\epsilon} \cap U$ $\Rightarrow \exists P' \in A \cap
U$, such that $\mid P - P' \mid < \epsilon$. Let $f_{\epsilon}$ be the holomorphic function that
has zero's on $A_{\epsilon} \cap U$, we use the notation $f_{\epsilon} \sim_{m} g$.

\vspace*{.5cm}

Rouch\'e's theorem can be used to define a "monopropy" ideal $J^{\epsilon}_{h}$ such that, $f \in
J^{\epsilon}_{h}$, if there exists a homomorphism $g$, such that
$$\mid h(f) - g(f) \mid < \epsilon \mid h(f) \mid \text{ on } \delta \Omega$$
for a small $\epsilon > 0$, where $\Omega$ is a small disc-neighborhood of a point in
$N(J_{h})$. This condition means that $h(f)$ and $g(f)$ have the same number of zero's (counted with
multiplicity) in $\Omega$. The condition particularly means that  $h(f) \sim_{m} g(f)$.
Assume further that $g$ is an algebraic homomorphism, $g^2(f)=g(f^2)$ that commutes with $h$, such that $g(f)-f^2 \in J_{h}$,
that is $h(f^2)=hg(f)$. If we have $h(f)-g(f) \in J_{h}$, then $h(f^2)=h^2(f)$. Otherwise, if $f
\in J^{\epsilon}_{h}$ with $Z_{f} \subset Z_{h(f)}$, $h^2(f)-h(f^2)=h(h(f)-g(f)) \sim_{m} 0$, that is $h^2(f) \sim_{m} h(f^2)$
over $J^{\epsilon}_{h}$. In the simplest case $g$ can be chosen as the identity operator.

\vspace*{.5cm}

Given our analytic homomorphism $h$ over an ideal $(I)$, there is an algebraic homomorphism
$g$, such that locally $h \sim_{m} g$ over $(I)$.
Assume $h$ is a homomorphism with $Z_{f} \subset Z_{h(f)}$ and $h^N$ injective. Further,
$\Omega=\text{nbhd}Z_{f}$ and $h(f) \neq 0$ on $\delta \Omega$, where $\Omega=\Omega(R)$ and
$1/R < \mid 1 + \sum_j \alpha_j \mid$, where $h^j(f)=\alpha_jf$. Assume further that $R$ is so
large that $RC > 1/ \epsilon$, for a given $\epsilon > 0$, where $C$ is such that $\mid h^N(f) -f
\mid = C \mid f \mid$. Thus, there exists an open set $\Omega=\Omega(R)$ such that $\mid h(f)-f \mid <
\epsilon \mid f \mid$ on $\delta \Omega$. If $\alpha=h(f)/f$ is holomorphic on $\delta \Omega$ and
$\mid \alpha \mid > \epsilon$ on $\delta \Omega$, the existence of an algebraic homomorphism
(dependent on $\Omega$) follows.

\vspace*{.5cm}

For $h$ such that $h^2$ is injective, we have existence of a square $g$ with respect to $h$, $h^{-1}gh=\alpha$, for
$\alpha(f)=f^2$. Thus, $h^{-1}g^2h=\alpha^{2}$ and so on. If $h \alpha=\alpha h$, we have
$\alpha^{-1} g= \sqrt{g}=1$. Again $h(f)-g(f)=h(f)-f^2$ and $h(f)-g(f) \in J_{h}$ $\Leftrightarrow$
$h^2(f)=h(f^2)$. Let
$$I^{\epsilon}_{f}=\{ g \quad \tau_{\epsilon}g=0 \quad Z_{f} \cap \Omega \quad \exists \epsilon
\quad \text{small} \},$$ where $\tau_{\epsilon}$ denotes translation and $\Omega$ is a small domain of
holomorphy around $Z_{f}$. If $h$ has integral representation, we have $h \text{:}
I^{\epsilon}_{f} \rightarrow I^{-\epsilon}_{f}$ continuously. This means that $h$ does not spread
zero's much. For $h$ such that $h(f) \sim_{\infty} \mbox{ Re }h(f)$ with $h^2$ injective and $f$ reduced, we have $h(f)^2(\zeta)=0$
implies $h(f)(\zeta)=0$. Assume $g$ the square with respect to $h$ as above and
$h^2(f)=g(\tau_{\epsilon}f)^2$, then $h^2(f)=0$ implies
$f=0$. In the same way $\tau_{\epsilon}h^2(f)=0$ implies $\tau_{\epsilon}f=0$, why $\epsilon=0$.
Note that $h^2$ is reduced if and only if $g^2$ is reduced. Thus, with these conditions, $h^2(f)=h(f)^2$.

\newtheorem{c_monotropy}{ Proposition }[section]
\begin{c_monotropy} \label{c_monotropy}
Assume $h$ a homomorphism with $\mbox{ Re }h \sim_{\infty} h$ and $h^N$ injective, for $N > N_0$,
$N_0$ some positive number. Then we have over $J^{\epsilon}_h$, that $h$ is an algebraic homomorphism.
\end{c_monotropy}

Proof: \\
The condition $f \in J^{\epsilon}_h$ means that there exists an algebraic homomorphism $g$, such
that $h^{N-1}(f)=g^{N-1}(\tau_{\epsilon}f)$. Note that $\epsilon$ depends on both $f$ and
$\zeta$. Then $h^{N}(f)=g^{N-2}(\tau_{\epsilon}^2f^2)=g^{N}(\tau_{\epsilon}f)$ means that
$\epsilon=0$, thus $h^{N}(f)=h(f)^{N}$. $\Box$

\vspace*{.5cm}

Define $(I^{m}_{h})=\{f \quad h(f)=\tau_{\epsilon}g(f) \}$ for an algebraic homomorphism $g$ and
for some $\epsilon$. We are assuming $h \tau_{\epsilon}=\tau_{\epsilon}' h$, where also
$\tau_{\epsilon}'$ is assumed to be translation. Let $B_{\epsilon}=\{f \quad \tau_{\epsilon}f/f
\mbox{ not constant } \}$. We form $(I^{m}_{h})$ over $(I)=\{f \quad B_{\epsilon} \mbox{ bounded and symmetric with respect to the origin
}\}$. We are assuming the translation (weakly) compact over $(I)$, why $g(B_{\epsilon})$ has the
same properties over $(I)$. If $f \in (I^{m}_{h})$, we have $\frac{\tau_{\epsilon}' f}{f}=\frac{g^{-1}
h(f)}{f}$. Further, if $h(f) \sim_{m} g(f)$, we have (with a different $\epsilon$) $h(f)^{2}
\sim_{m} g(f^{2})$, that is the representation is invariant for iteration. Assume $h(f^{2}) \sim_{m}
h(f)^{2}$ (algebraic translation), then $f^{N} \in (I^{m}_{h})$ implies $f \in (I^{m}_{h})$, that
is we have a radical ideal.

\section{ Measure zero domains}
We know that polynomial operators are never regularizing as pseudo-differential operators and this
means that necessarily polynomials do not have an infinite zero in the infinity (it has not complex
dimension). Assume our system has polynomial right hand sides and consider an analytic surface
$S=\{ y=G(x) \}$, for a polynomial $G$. If $\int_{(I)} (1 + \mid G'(x) \mid^{2}) d x=0$ and
if $(I)$ is a neighborhood of the origin, we have that the measure for $(I)$ is zero. If
$\int_{(I)} \mid G'(x) \mid^{2} d x=0$, then as we shall see, there is a segment in $(I)$ such that
$G'(x)=0$ on this segment. Then according to Hurwitz theorem, $G'(x) \equiv 0$ on $(I )$. The conclusion is thus $\int_{(I)} d x >0$ $\implies \int_{(I)} \mid G'(x)
\mid^{2} d x > 0$. Thus for all polynomials we have that the zero-set does not have positive area.

\vspace*{.5cm}

Particularly, assume $w$ is a polynomial in $x,y$ given our conditions and $\int_{(\tilde{I})} w d
\sigma=\int_{(\tilde{I})} Q_{x} - P_{y} d \sigma$ (we assume $\int_{(I)} w dx dy=\int_{(\widetilde{I})} w d \sigma(t)$).
Given a finite integral on both sides, we see that
there exist $\omega_{j} \neq 0$ such that $\omega_{j} \rightarrow \omega$ uniformly on compact sets. Further,
that if $(I)$ is a neighborhood of the origin
$$ \int_{(\tilde{I})} \omega d \sigma=0 \Rightarrow \sigma(\tilde{I})=0$$
Note also that if $\eta(x)=h(x)/x$ and $\int_{(I)}(1 + \mid \eta'(x) \mid^{2} ) d x < \infty$, we have that the
surface $\{ y = \eta(x) \}$ is normal. However, if $\int_{(I)} \mid x \eta'(x) \mid^{2} dx=0$
and $\int_{(I)} d x >0$, according to Hurwitz theorem, we have that $x \eta'(x) \equiv 0$ in $(I)$ and $\frac{d h(x)}{d
x}\equiv \frac{h(x)}{x}$ in $(I)$.

\vspace*{.5cm}

 Assume $\Omega$ an open subset of a
Stein-variety $V$. We can define the quasi-porteur corresponding to an analytic functional in
$H'(V)$ as an open set where the functional has a local representation. A porteur is a set where we have representation of $T \in H'(X)$ by a measure with compact support. If
$\Omega \subset X$ is a quasi-porteur for $T \in H'(X)$, we have for all neighborhoods $W$ of $\Omega$,
existence of a measure $\mu_{W}$ with compact support in $W$, such that $T(\phi)=\int_{W} \phi d
\mu_{W}$ for all $\phi \in H(X)$.  If $\Omega$ is a complex-analytic quasi-porteur for $\sigma$, we can determine
a measure $\mu$ with compact support on the boundary of $\Omega$, $\Gamma'$, such that
\begin{equation} \label{f_rep}
\int_{\Omega} \phi d \sigma = \int_{\Gamma'} \phi d \mu \qquad \phi \in H(X)
\end{equation}
such that $\mid \int_{\Gamma'} \phi d \mu \mid \leq C \sup_{\Gamma'} \mid \phi(x) \mid$ for
$\Gamma' \supset \mbox{supp } \mu$. If $T$ is represented by $\mu$ in $H'$, we can say for the restriction
to a line $L$, that $\widehat{T}(\zeta)=0$ means that $\mu \mid_{L}=const. \delta_{L}$.
In order to, starting from a dynamical system determine the porteur, we note
that every subset $W'$ such that $\int_{W'} d \sigma=0$, is a subset of the quasi-porteur. If we limit
ourselves to the case with analytic $\Omega$'s, this is also a subset of the porteur.

\section{ Symmetry}
\label{section:f_symmetry}
Consider $I_{\Gamma}(F)=\int_{\Gamma} F(x,y) d x$ and $I_{\Gamma}^{*}(F)=\int_{\Gamma}
\overline{F}(x,y)d y$. Assume  for a domain $\Omega \subset (I)\times (I)$ and $\Gamma=bd \Omega$
$$\int_{(\Omega)} \frac{d^{2} F}{d y d x}- \frac{d^{2} F}{d x d
y}d x d y = \int_{\Gamma} \frac{d F}{d x} d x + \frac{d F}{d y }d y=0$$ This
means that $\frac{d^{2} F}{d y d x}=\frac{d^{2} F}{d x d y}$ over
$(\Omega)$. Further, $I_{\Gamma}(F_{x})+ I_{\Gamma}^{*}(\overline{F}_{\overline{x}})=0$ under the
condition above, why if $\overline{F}$ is holomorphic, we have $I_{\Gamma}(F_{x})=0$ and $F$ is
constant over $(\Omega)$. Note that if $I_{\Gamma}(F)=0$ for a polynomial in $x$, $F$, then outside $(\Omega)$,
we have the polynomial case in $y$. Further, we can define the operators
$T_{\Gamma}(u)=\int_{\Gamma}F(x,y)u(y)d y$ and
$T_{\Gamma}^{*}(v)=\int_{\Gamma}\overline{F(x,y)}v(x)d x$. We can define a normal operator as an
operator such that the image of a bounded sequence of $x$ has a subsequence that converges
uniformly on compact sets.

\newtheorem{f_zero_integral}[f_J]{ Lemma }
\begin{f_zero_integral}
If, for a holomorphic function $F$, we have that $\int_{(\tilde{I})} F d \sigma=0$ we have that one of
the following propositions is true, either $\int_{(\tilde{I})} d \sigma=0$ or $F \equiv 0$ on $(\tilde{I})$.
\end{f_zero_integral}

In the polynomial case, we have seen that if $(\Omega)$ is a neighborhood of the origin, we have that
the first condition holds. Note that if $\int_{(\Omega)} \omega d x d y < \infty$ and assuming
$\omega$ polynomial in $x,y$, then also the measure for singularities must be zero. If $F$ is holomorphic and
$\int_{(\tilde{I})} d \sigma=0$, if also $\int_{(\tilde{I})} F d \sigma=0$, we must have that $F$ does not have an
infinite zero in $(\Omega) \ni 0$ according to Hurwitz. Does $F$ have to be a polynomial?. If $F \equiv 0$
on $(\Omega)$, we must have $F \in \mbox{ rad }I(\Omega)$, so $F$ is a polynomial in $x,y$. More precisely,
for $t=t_{1}+it_{2}$, if $t_{1},t_{2}$ are reduced, we have that $\int_{(\tilde{I})} d
\sigma(t)=0$ implies that $(\tilde{I})$ is algebraic. According to Nishino (cf. \cite{Ni}),
these sets are removable, for analytic continuation.

\vspace*{.5cm}

If $(I) \ni f_{j} \rightarrow f \in (\overline{I})$ with uniform convergence on compact sets, then
$N(f_{j}) \rightarrow N(f)$ (cf. \cite{Ni}). We will consider the sets
$V_{j}=\{ \zeta \quad \int_{(\Omega)} f_{j} dx dy=0 \}$ and $\widetilde{V}_{j}=\{ \zeta \quad \int_{(\tilde{I})} f_{j} d \sigma(t)=0 \}$.
For $f_{j}$ we thus have that either $\int_{(\tilde{\overline{I}})}d \sigma(t)=0$ or $ f \equiv 0$ on $(\tilde{\overline{I}})$.
If $d L_{T} \sim_{m} 0$, we know that there are $T_{j}$ such that $\frac{d L_{T_{j}}}{dx} d x + \frac{d L_{T_{j}}}{d y}d y \rightarrow 0$
as $j \rightarrow \infty$ and given that $\frac{d^{2}L_{T_{j}}}{d x d y}$,$\frac{d^{2}L_{T_{j}}}{d y d x}$ continuous,
$\int \big( \frac{d^{2} L_{T_{j}}}{d x d y}
- \frac{d^{2}L_{T_{j}}}{d y d x} \big) d x d y \rightarrow 0$ as $j \rightarrow \infty$. Thus, if $T_{j}$
is chosen such that the integrand holomorphic, we have that $(\overline{\Omega})$ is a domain for symmetry.
Assume $\int_{\tilde{I}} d \sigma(t) > 0$, then we have through the lemma that $ \widetilde{V}_{j}
\subset N(f_{j})$. Further, if the integral is taken over $N(f_{j})$, we must have that
$\int_{\tilde{I}}f_{j} d \sigma(t)=0$. This corresponds to convergence $\widetilde{V}_{j} \rightarrow
\widetilde{V}$ as $ j \rightarrow \infty$.

\vspace*{.5cm}

We can form the ideals $ (J)=\{ f \quad \int_{(\tilde{I})} f d \sigma(t)=0 \text{ on } \widetilde{V}
\}$ and the ideal under topology for normal convergence, where the measure that involves an algebraic homomorphism,
is considered as compact. Existence of a
normal and regular approximation, gives according to Hurwitz theorem the same conclusions and we
denote the corresponding ideal $(J^{*})$, which is considered as weakly compact. More precisely, assume
$I_{(I)}(w_{j}) \rightarrow I_{(I)}(w)$, uniformly on compact sets in $\zeta$. Further, if $w_{j}'
\rightarrow w$ as a normal and regular approximation, then
$I_{(I)}(w)=\lim_{j}I_{(I)}(w_{j}')=\lim_{j}I_{(I)}(w_{j})$.  We must have that
$\lim_{j}I_{(X)}(w-w_{j}') \neq 0$ implies $\sigma(X)=0$. Assuming $w$ algebraic in $x,y$, we must have
that $W=\{ (x,y) \quad \lim_{j}w_{j}(x,y)=w(x,y) \}$ has measure zero. For instance, let
$w_{j}(x,y)=w(x_{j},y_{j})$. In the same manner, $W'=\{\zeta \quad \lim_{j} w_{j}(x,y)(\zeta)=w(\zeta) \}$ has measure zero with respect to a relevant measure
$\exists F$ $\int_{W}d \sigma(t)=\int_{W'} F(d \zeta)$. Note that since the integral is invariant
for monotropy, we do not necessarily have that the $I_{\Omega}(w)$ are holomorphic.

\vspace*{.5cm}

Assume $F$ holomorphic and $F_{1}$ continuous and $F_{1} \sim_{m} F$ on $(I)$
and that $\int_{(\tilde{I})} (F_{1} - F_{j}) d \sigma(t) \rightarrow 0$ as $j \rightarrow \infty$,
for $F_{j}$ a holomorphic regular (in the sense of \cite{Ni} ) approximation. Since the integral is
invariant for monotropy, we have $\int_{(\tilde{I})}(F - F_{j}) d \sigma \rightarrow 0$
as $j \rightarrow \infty$, why if $F_{1} \neq 0$ on $(I)$ it follows that
$\int_{(\tilde{I})} d \sigma(t)=0$

\newtheorem{f_normal}[f_J]{ Proposition }
\begin{f_normal}
Assume $T_{\Gamma}$ defined as above with $T_{\Gamma}^{*}T_{\Gamma}=T_{\Gamma}T_{\Gamma}^{*}$ and $\int \mid F \mid^{2} d \sigma < \infty$, then
$T_{\Gamma}$ is a normal operator.
\end{f_normal}

Now define the operators
$$T_{\Gamma}(u)=\int_{\Gamma} G(x,y)u(y) d y \text{   and   } T_{\Gamma}^{*}(u)=-\int_{\Gamma}\overline{F}(x,y)u(x)d x$$ so that
$T_{\Gamma}(1)-T^{*}_{\Gamma}(1)=\int_{\Gamma}G dy+\overline{F}d x$. If $F,G$ are real and $\int_{\Gamma}F d x+G d
y=0$, then $- dG/d x=d F / d y$ and through mass-conservation $d G/d y=-d F/d x$. Thus, $d G/ d T=d
F/d T$. If $G_{x} \sim_{m} F_{y}$, then $\int_{\Gamma} F dx + G dy=0$ and if $G={}^{t}F$ with the notation
$(F \frac{d}{d y})^{t}=\frac{d}{d x} {}^{t}F$, so that $G_{x}-F_{y}=(F \frac{d}{d y})^{t}-(\frac{d}{d
y}F) \sim_{m} 0$, we have that $< T_{\Gamma}(u),v>=<u,T_{\Gamma}(v)>$. Define
$S_{\Gamma}^{*}(u,T)=-\int_{\Gamma}F(x+T,y)u(x)d x$ and $T_{\Gamma}(u)=\lim_{T \rightarrow 0}
S_{\Gamma}$.

\vspace*{.5cm}

Assume symmetry in the tangent space $TS=\{ \frac{d}{d x}G \sim_{m} \frac{d}{d x}M \}$, for $M$
real, that is $G \sim_{m} M + iC(y)$, where $C$ is independent on $x$. For instance, if $F-F^{*} \in (I)(\Omega)$
where $\Omega$ is a domain for symmetry, we can form the closure $(\overline{I})(\Omega)$, then we have symmetry
in the tangent space over $\Omega$. If $G={}^{t}F$, we thus have
$\int_{\Omega} \big( \frac{d {}^{t}F}{d x}- \frac{d F}{d y}\big)u \otimes v d \sigma=0$, for all
$u,v \in H(X)$, for ${}^{t}F,F$ in $TS$. Thus, $S_{\Gamma}=S^{*}_{\Gamma}$ and $T_{\Gamma}$ can
be considered as the limit of self-adjoint operators.

\vspace*{.5cm}

Assume $0 < \int_{\Omega} d \sigma < \infty$ and $\int F d \sigma < \infty$. Assume
$\int_{\Omega}(F-A)d \sigma=0$ for a constant $A$. Further, $F(\phi)=\int_{\Omega}F(x,y)\phi d
\sigma$ and $F F^{*}=F^{*} F$. We can prove that $F^{*}(\phi)=\overline{A}\int_{\Omega} \phi d \sigma$.
If $\int_{\Omega} \mid F \mid^{2} d \sigma < \infty$ then $ S=\{ F \overline{F}=\mid C \mid^{2} \}$
is normal (cf. \cite{Ni_92}). If $\{ F=C \}$ is normal, then also $S$ is normal. Determine a homomorphism $h$, which we assume
self-adjoint and such that $h^{2}$ is locally injective and finally
\begin{equation}
      \left \{
\begin{array}{lr}
h(F)=0 \Rightarrow \mid F \mid^{2}=\mid C \mid^{2} \\
dh(F)=0 \Rightarrow F=C
\end{array} \right.
\end{equation}
We can assume that $\{ F=C \} \cap \{ \mid F \mid= \mid C \mid\}$ is minimally defined. More
precisely, if $P$ is a singular point, we have that $P$ can be reached through a regular and normal
approximation from all directions. Particularly, $F(\phi)=A\int \phi d
\sigma$ can be approximated by $\mid F \mid (\phi)$. For $\phi$ in a bounded set corresponding to
the normal tube (cf. \cite{Ni}), the corresponding operator $T_{\Gamma}(\phi)$ is normal.

\vspace*{.5cm}

We now claim that if the vorticity $w \in (\mathcal{B}^{\textit{m}})'(\Omega)$ and if $f,g \in
\mathcal{B}^{\textit{m}}(\overline{\Omega})$, where the monotropy-condition is on the boundary $\Gamma$
to $\Omega$, then if $\int_{\Omega}w_{j}(f,g)d \sigma \rightarrow \int_{\Omega} w(f,g) d \sigma$
through a normal and regular approximation and if $ww^{*}=w^{*}w$ and $\int \mid w \mid^{2} d \sigma < \infty$, then the corresponding
operator $T_{\Gamma}$ is normal.
If $\gamma_{T}=(f_{T},h(f_{T}))$ is the regular
approximation and if the dependence of $T$ is polynomial in $\frac{d \gamma_{T}}{d T}$, according
to Hurwitz theorem, since polynomials do not have zero's of infinite order, the zeroset for the
polynomial must have measure zero. We can assume the parameter space in a domain of holomorphy or
simpler, given existence of regular and normal approximations, we assume algebraic dependence of
the parameter $T$ in the tangent space, why all normal approximations algebraically dependent of
the parameter $T$ in the tangent space, can be assumed regular (by adding a regular approximation).

\section{ The nil-radical with respect to $d \sigma$ }
As long as we are considering characteristic sets of measure zero, analytic continuation can be
assumed. If $f^{N} \equiv 0$ on $V$ such that $\sigma(V)=0$, then we have seen that $f^{N}$ has a polynomial
representation. If $\sigma(V) > 0$, we have that $ f \in \mbox{ rad }I(V)$, but as $I(V)$ is a
geometric ideal, we have $f \in I(V)$. Let $$\mathcal{Z}=\{\textit{g} \quad \int \textit{g d} \sigma=\textsl{0} \}$$ and $X=\{
(x,y) \quad \int g d \sigma=0 \quad g \in \mathcal{Z} \}$ and $\Omega=\{ \zeta \quad \int g d \sigma=0
\quad
g \in \mathcal{Z} \}$. If $g=const.$ on a set $V$ such that $d g=0$ on $V$ and if $g_{j}
\rightarrow g$ in a normal and regular approximation, such that $d g_{j} \neq 0$ and
$\int(g_{j}-g)d \sigma \rightarrow 0$ as $j \rightarrow \infty$ then we see that $X=\{ (x,y) \quad
\int (g_{j} - g) d \sigma=0 \quad j \rightarrow \infty \}$ is of measure zero. Assume that $\Omega$
is a domain of holomorphy and $f \neq 0$ on $\Omega \backslash V$, $\int_{\Omega \backslash V} f d
\sigma=0$, since $\int_{\Omega \backslash V} d \sigma =0$ if $f^{N} \equiv 0$ in $\Omega \backslash V$,
$f^{N}$ has polynomial representation there. If $\mathcal{Z}$ $\ni f_{j} \rightarrow f$ in a normal and
regular approximation, then $f \in \mathcal{Z}$, so $\mathcal{Z}$ is closed under normal convergence. 

\vspace*{.5cm}

We can easily consider $({L}_{f})=\{ g \quad \int_{(I)}(f-g)d \sigma=0 \}$, for $f$ such that
$\int_{(I)} f d \sigma < \infty$, so the case with finite area can be compared to the ideal
$\mathcal{Z}$. We have earlier seen that $f^{N} \equiv 0$ on $(I)$ implies $f=c$ on a set $X$.
Let $X_{1}=N(\{ f \quad \int_{(I)} f d \sigma=0 \})$ and $X_{2}=N(\{f \quad \int_{(I)}f^{2} d
\sigma=0\})$ and so on. Then $\ldots \subset X_{N} \subset \ldots \subset X_{2} \subset X_{1}$
and $X_{1}=\cup_{j} X_{j}$
\newtheorem{f_red_sigma}[f_J]{ Lemma }
\begin{f_red_sigma}
With the conditions above and if $X_{1}$ is a domain of holomorphy, we have $\sigma(X_{j})
\downarrow 0$ as $j \uparrow \infty$.
\end{f_red_sigma}
It is easily seen that existence of a $f$ such that $\int_{(I)} f^{N} d \sigma=0$ and $f=const \neq 0$
on $X_{N}$ implies $\int_{X_{N}} d \sigma=0$. This gives a modified polynomial class
\begin{equation}
(R_{\Omega}) \qquad \int_{\Omega} f d \sigma=0 \Rightarrow \sigma(Z_{f})=0
\end{equation}
If $\Omega$ is of measure zero, then holomorphic $f$ are in $(R_{\Omega})$. If $\sigma(\Omega) >
0$, then polynomials are in $(R_{\Omega})$. For pseudo-differential operators, it may be more
interesting to consider
\begin{equation}
(R_{\Omega}') \qquad \int_{\Omega} f(x + T)-f(x) d \sigma=0 \Rightarrow \sigma(\Delta_{\mathbf{C}})=0
\end{equation}
where $\Delta_{\mathbf{C}}=$ $\{ T \quad f(x+T)-f(x)=0 \quad \forall x \}$.

\vspace*{.5cm}

We can also consider the class of finitely generated symbols
\begin{equation}
(FR_{\Omega}') \qquad \int_{\Omega} f^{N}(x + T)-f^{N}(x) d \sigma=0  \Rightarrow
\sigma(\Delta_{\mathbf{C},(\textsl{N})})=0
\end{equation}
where $\Delta_{\mathbf{C},(\textsl{N})}$ $=\{ T \quad f^{N}(x + T)-f^{N}(x)=0 \quad \forall x \}$.
This defines, modulo monotropy, the class of symbols that become reduced with respect to the
measure after finitely many iterations.

\section{ Ideals of holomorphy with induced topology }
It can be proved, that solvability for the second Cousin problem for an ideal of holomorphy $(I)$,
depends on choice of topology for $(I)$. Assume $(I)$ a geometric ideal and $\Omega$ a neighborhood of its
zero's. For the first Cousin problem to be solvable we require that $\Omega$ is a holomorphically
convex domain. For this domain, we have locally a pseudobase of holomorphic functions $F_1, \ldots,
F_{\nu}$. We are interested in topologies for which this base can be selected as global. Such a
topology is said to have Oka's property.

\vsp

Given a covering of $\Omega$, a $\nu-$ dimensional domain of holomorphy ($\nu >1$) by open sets
$\{ \Omega_i \cap \Omega_j \}$ and associated continuous functions $g_{i,j}(\zeta)=e^{-A\rho_{i,j}(\zeta)}$
for a constant $A$, where $\rho_{i,j}$ a real or complex norm on $\Omega$. Assume $\mid \rho_{i,j}
\mid=1$ and let $\rho_i(\zeta)=\zeta \chi_{\Omega_i}$, where $\chi_{\Omega_i} \in C^{\infty}_0(\Omega)$
are of modulus $\leq 1$ and $=1$ on $\Omega_i$. Then $g_{i,j}(\zeta)=e^{-A \rho_i
\overline{\rho_j}}$. The topological space $( \Omega_i, g_i )$ has Oka's property and $\parallel G
\parallel_{\rho,A}=\sup_{i,j} \mid g_{i,j}(\zeta)G(\zeta) \mid$ is a Banachspace over $(I)$.

\vsp

The choice of $\rho_{i,j}$ can be varied, for instance if
\begin{equation}
c_1(1 + \mid \zeta \mid)^A \leq \mid G(\zeta) \mid \leq  c_2(1 + \mid \zeta \mid)^B \mid \eta \mid^{-k}
\end{equation}
where $\zeta=\xi + i \eta$ and $0 < \mid \eta \mid < r$ a positive constant $r$ and where $A,B$ are arbitrary constants.
If $G$ is analytic in the tube $T^{C}$, $C$ an open, connected cone in $ \mathbf{R}^{\nu}$, then for a compact
subcone $C'$, $b_{\Gamma}G$ exists uniquely in $ \mathcal{S}'$. If
\begin{equation} \label{b_inv}
c_1'(1 + \mid \xi \mid)^{A'} \leq \mid G(\zeta) \mid \leq c_2'(1 + \mid \xi \mid)^{B'}\mid \eta
\mid^{-k}
\end{equation}
for constants $A',B'$ and $\zeta=\xi+i\eta$,
we have a topology which can be compared with (included in) the topology for invertible measures. In both
these cases, $\Omega$ has Oka's property.

\vsp

 Assume the topology on $\Omega$ is defined through
$g_i(\varphi)(\zeta)=\alpha I_{\widehat{F}}({}^tP\varphi)(\zeta)$ for a fixed $\varphi \in
 \mathcal{D}_{\textit{L}^{\textsl{2}}}$. For $\zeta \in \Omega_i \cap \Omega_j$, the corresponding $$g_{i,j}(\varphi)(\zeta)= \alpha
\beta^{-1} I_{\widehat{F}}( {}^t(P/Q) \varphi)(\zeta)$$ can be split into $g_i,g_j$ as above. The corresponding
topological space has Oka's property.
Assume $(I)_{ \mathcal{\widehat{D}'}_{\textit{L}^{\textsl{2}}}}$ an ideal of holomorphy in the topology of\\
$ \mathcal{\widehat{D}'}$ ${}_{L^2}$, that is $f \in (I)_{ \mathcal{\widehat{D}'}_{\textit{L}^{\textsl{2}}}}$ means that
$f=\sum_j P_j \widehat{F}_j$ for $\widehat{F}_j$ reduced in both variables separately and in $H \cap L^2$ and $P_j$ polynomials. We can form the
quotient ideal of holomorphy $$(Q)=\{ g \ \exists \phi \in H \ \phi g \in (I)_{ \mathcal{\widehat{D}'}_{\textit{L}^{\textsl{2}}}}
\}.$$ If $(Q)$ is equipped with ${ \mathcal{\widehat{D}'}_{\textit{L}^{\textsl{2}}}}$-topology, $\phi$ can be chosen
as polynomials.
 According to the theory of integral equations, given $g,F \in L^2$,
there is a $\varphi \in {L^2}$, such that $g=I_{\widehat{F}}({}^tP \varphi)$, for a polynomial $P$. Assume $Q'$
a polynomial such that $Q'I_{\widehat{F}}=I_{\widehat{F}}{}^tQ$, if ${}^tQ\psi \in L^2$ we have a $\varphi \in L^2$ such
that ${}^tQ\psi={}^tP\varphi$. This gives a representation of $g \in
(Q)_{ \mathcal{\widehat{D}'}_{\textit{L}^{\textsl{2}}}}$
 Note that a symbol in the ideal $(I)$ equipped with this
topology, will get a representation $f^{\lambda}(\zeta)=\sum_j I_{T_j}(\varphi)(\zeta)$ for
$T_j \in ( Q )_{ \mathcal{\widehat{D}'}_{\textit{L}^{\textsl{2}}}}$, for a $\varphi \in
L^2$ and we can assume the support for $\varphi$ adjusted to the corresponding subset in $\Omega$.

\section{ The hypoelliptic radical for induced topology }
\label{sec:b_topology_ideals}

Assume now $P$ a partially hypoelliptic, self-adjoint operator and $T=P(D)F$, where $F \in
L^2( \mathbf{R}^{\nu} \times  \mathbf{R}^{\nu})$. Further,
$$ \tau_h({}^tP\varphi)(\xi)={}^tP(\xi+h)\varphi(\xi+h) \qquad {}^tP \varphi \in L^2$$
$$ g_{h,P}(\varphi)=\tau_h({}^tP \varphi)-{}^tP\tau_h\varphi \qquad {}^tP \varphi \in L^2$$
For a kernel $F$, such that $I_{\widehat{F}}$ reduced (left semi-bounded), since $I_{\widehat{F}}$ is a compact
operator, if $\varphi \rightarrow 0$ in $L^2$, then $I_{\widehat{F}}(\tau_h \varphi) \rightarrow 0$
and the reducedness implies $\tau_h \varphi \rightarrow 0$ in $L^2$. Thus $\tau_h$ is a compact
operator for this topology.

\vsp

For test-functions ${}^tP \varphi \in L^2$, we have $$I_{P\widehat{F}}(\tau_h \varphi)=
I_{\widehat{F}}({}^tP \tau_h \varphi)=I_{\widehat{F}}(\tau_h({}^tP \varphi )) -
I_{\widehat{F}}(g_{h,P}(\tau_h \varphi))=I_1 + I_2$$
where we can assume $g_{h,P}$ is reduced. If $\varphi \rightarrow 0$ then ${}^tP \varphi
\rightarrow 0$ in $L^2$ and we have $I_1 \rightarrow 0$. As $I_{P\widehat{F}}(\tau_h \varphi)
\rightarrow 0$, we have $I_2 \rightarrow 0$ and $\tau_h \varphi \rightarrow 0$. Thus $\tau_h$ is a
compact operator also in this topology. Without the condition that $I_{\widehat{F}}$ is reduced, we
only have that translation is a weakly compact operator.

\vsp

According to Malgrange (cf.  \cite{Mlg} Theorem 2, Prel. ) if ${}^t\tau: F' \rightarrow E'$ a linear continuous
operator between $( \mathcal{F})$-spaces such that $\tau(E)$ weakly closed in $F$ then ${}^t\tau$ is a homomorphism and conversely.
For the spaces where we have proved that the translation operator is weakly closed, we can apply the
argument in section (cf. \cite{jag_0}) on $a_{\Gamma}={}^t\tau$ and we see that the symbol $p^N$ is reduced for
the corresponding topology. We now use that $A_{\lambda}^N-P_{\lambda}^N=H_{\lambda,N}$ with
$H_{\lambda,N}$ regularizing and assuming $A_{\lambda}^N \sim \mbox{ Re }A_{\lambda}^N$ (equivalence of ideals), we see that
$A_{\lambda}^N$ is hypoelliptic with the initial topology. Note that for self-adjoint polynomial
operators, $N=2$ will do.

\vsp

Let's return to our initial approach of $f=\mbox{ Re }F$ in $C^{\infty}$ and $\mbox{ Exp
}$.  Except for lacunary cases, we assume that $F^N$ is self-adjoint and hypoelliptic in an ideal $(I_{HE}^*)$, with induced topology
of Schwartz-type. If the ideal of holomorphy is completely symmetric, that
is $ F \in (I_{HE}) \Leftrightarrow F^* \in (I_{HE})$, we can assume $ F \sim F^*$, which implies
$\mbox{ Re }F^N \sim (\mbox{ Re }F)^N \sim f^N$. More precisely, let $$r(I_{HE}^*)=\{ F \in
H( \mathbf{C}^{\nu}) \quad \textit{F}^{\textit{N}} \in (\textit{I}_{\textit{HE}}^*) \quad \exists \textit{N} \}.$$ That is entire functions for which the $N$'th
iterate corresponds to a self-adjoint and hypoelliptic operator. Let $$S_r= \mbox{ Re }r(I_{HE}^*).$$
The symbol $f$ can naturally be compared with the ( H\"ormander- ) class of symbols $S^m$, where we
prefer to vary the complexification. For instance let $G=g+ih$ be a symbol and $S^m_1$ the symbol-class with
standard complexification. If $g \in S_r$ and $G^N \in S^{Nm}_1$, then $g \in S^m$ and $G \in
S^m_1$. We have noted that the complexification that corresponds to operators hypoelliptic in
$ \mathcal{D'}$, $S^m_2$, must be different. Assume $G$ corresponds to an operator hypoelliptic in $ \mathcal{D'}$, $G^N=g_N+ih_N \in S^{Nm}_1$ for some $N$. If also $g \in S_r$,
we must have $g^N \sim g_N \in S^{Nm}$. Assume existence of a $h' \in S^m$ such that $(h_N)^{1/N}
\sim h'$ (a sufficient condition is that $h_N$ is reduced), this defines $G \sim g+ih'$ (equivalence of ideals)in
$S^m_2$.

\newtheorem{b_S-r}{ Proposition }[section]
\begin{b_S-r}
The non-standard complexification of the symbol-ideals $S_r$, $S_r+iS_r^{1/N}$, produces (modulo $\sim$)
the class of hypoelliptic operators.
\end{b_S-r}

\section{ Contractions of formulas }
\label{sec:contraction}
Assume $(I)$ an ideal with Schwartz-type topology and a weakly compact translation. Consider for
$\zeta \in \mathcal{B}$, a bounded domain in $\Omega$ pseudoconvex
\begin{displaymath}
   \left. \begin{array}{ll}
   B(\zeta) \quad \longrightarrow   &  \alpha(T)B(\zeta) \\
   \downarrow   &  \downarrow \\
   B(\zeta + T) \longrightarrow & B(\zeta) + c_{T}
   \end{array} \right.
\end{displaymath}
If $B$ is reduced, we have that $B(\zeta + T)=\alpha(T)B(\zeta)$, where $\alpha(T)$ is constant for
$\mid T \mid > \mid T_{0} \mid$ and we see that $B$ does not change sign in the $T$-infinity for
$\zeta \in \mathcal{B}$. If $B=h(\psi) \in (I)$ with $\psi$ reduced, for $h$ algebraic, the earlier
results can be repeated. For $h=\tau_{\epsilon}w$ and $w$ algebraic, the results follow modulo
monotropy. Particularly, for $\psi$ reduced, $\eta=h(\psi)/\psi$ does not change sign in the infinity.

\vspace*{.5cm}

Assume $B$ reduced. Then the condition $B(\zeta+T)/B(\zeta)=1 + c_{T}/B(\zeta)$, for $\zeta \in \mathcal{B}$
gives $c_{T}/B \sim_{m}0$. Assume now that $dB(\zeta + T)=dB(\zeta)+dL_{T}(\zeta)$. We will put the condition
that $dL_{T}/dB \sim_{m}0$. It is sufficient through R\"uckert's theorem, if $dB$ is reduced, that
$dL_{T} \sim_{m} 0$. Consider $$ \Delta_{P/Q}=\{ T \quad \frac{P}{Q}(\zeta + T)-\frac{P}{Q}(\zeta)=0 \quad \zeta \in
\mathcal{B} \}$$ (the set can without problem be extended to $\cup^{\infty}_{j=1} B_{j}$) with $B_{j}$
as above)
\begin{equation} \label{f_contr}
\frac{-dL_{T}/dy}{dL_{T}/dx}=\frac{P}{Q}
\end{equation}
 Let's denote
the set of $T$ such that the last equality holds $\Omega_{L}$. Thus $\Omega_{L}=\Delta_{P/Q}$. We
say that $L_{T}$ is reduced for contraction if $$T \in \Omega_{L} \qquad
 \frac{-dL_{T}/dy}{dL_{T}/dx}=\frac{dx}{dy}=\frac{P}{Q} \implies T=0$$
We say that $L_{T}$ is regular (for contraction) if for regular $(x,y)$ and $T \in \Omega_{L}$, we
have $dL_{T}/dT \neq 0$.

\vspace*{.5cm}

Assume $\alpha(T)=\psi(\zeta + T)/\psi(\zeta)$, for $\zeta \in \mathcal{B}$ and $\psi=e^{\varphi}$
reduced. Assume $h$ algebraic and $\gamma(T)=h(\psi)(\zeta + T)/h(\psi)(\zeta)$ with $d h(\varphi)
\sim_{0} d \varphi$. Over $\Omega_{L}$ we have that (\ref{f_contr}) is obviously independent of $T$.
As a consequence of Kohn (\cite{Kohn}) if $B(\zeta + T)=\gamma(T)B(\zeta)$ for $B$ real and
reduced, we have that $\gamma$ is a polynomial in $T$ $\Leftrightarrow$ $L_{T}(\zeta)$ reduced in
$\zeta$ with respect to lineality. Further, the set of $T$, such that $\gamma$ not constant outside
the origin is symmetric with respect to the origin and bounded. We will denote this set
$\Omega_{\gamma}$.
For $0 \in \Omega_{\gamma}$, any regular approximation of the origin must go through
$\Omega_{\gamma}$. For $T \notin \Omega_{L}$ but close and $T \in \Omega_{\gamma}$
$$ \frac{-d L_{T}/d y}{d L_{T}/dx}=\frac{-(\gamma(T)-1)dB/dy}{
\eta_{1}'(x)-\eta'(x)}$$ where $\eta'(x)=dh(x)(\zeta)/x$ and $\eta_{1}'(x)=dh(x)(\zeta + T)/x$.
For $h$ algebraic and such that $h^{2}=1$, we have that $$\frac{dh(x)(\zeta + T)}{x} \sim_{0}
\frac{\alpha(T)dh(x)(\zeta)}{x}$$
Thus, over $\{ y=B(\zeta) \}$ we have that $ \frac{d L_{T}/dy}{d L_{T}/d x} \sim_{0}
\frac{\gamma-1}{\alpha-1}{\eta'}^{-1}(x)$.
 Outside $\Omega_{\gamma}$, the quotient is
bounded by a constant. We have
$\frac{\gamma}{\alpha}(T)=\frac{\eta_{1}}{\eta}(x)\frac{x(\zeta)}{x(\zeta + T)}$ $\sim_{m}
\eta \big( \frac{\tau_{T}x}{x} \big)$. Over the half-space $\{x \quad dh(x) \geq \mu x \}$, we know
that ${\eta'}^{-1}$ does not change sign, if $x \neq 0$.
We have earlier seen that if $\eta$ is bounded in the infinity, we have that $\eta_{1} \sim_{m} 0$
and the associated ideal ( with respect to $h$ ) has a global pseudo-base.

\vspace*{.5cm}

Assume now right-hand sides $P,Q$ to the dynamical system such that $P,Q$ becomes reduced by
iteration. We can prove that $\Delta_{P/Q} \downarrow \{ 0 \}$ for iteration. Define
$L_{(2),T}(\zeta)=B^{2}(\zeta + T)-B^{2}(\zeta)$ for $\zeta \in \mathcal{B}$. We obviously have that
$\Omega_{L_{(j),T}} \downarrow \{ 0 \}$ as $j \uparrow \infty$. We can say that $L_{(j),T}$ is
reduced for contraction. We can assume that on $\Omega_{\gamma}$ we have $dL_{T}/dT \neq 0$, if $L_{T}$
is reduced for contraction. According to the condition $\Omega_{L} \downarrow \{ 0 \}$ for
iteration, we can prove $\{ (L_{T}-c_{T})=\frac{d}{dT}L_{T}=0 \}$ is nowhere dense outside
$\Omega_{\gamma}$. Thus over $L_{T}=c_{T}$, we can assume $L_{T}$ regular (for contraction).

\newtheorem{f_contr_red}[f_J]{ Proposition }
\begin{f_contr_red}
Given a finitely generated system with polynomial right hand sides $P,Q$, assume that
$\Delta_{P/Q}$, the lineality corresponding to $P/Q$, is $=\{ 0 \}$. Then $L_{T}$ is reduced for contraction.
If $P,Q$ are in the radical to reduced polynomials, we have $L_{T,(j)}=B^{j}(\zeta +
T)-B^{j}(\zeta)$, becomes reduced for contraction, for $j$ sufficiently large.
\end{f_contr_red}

\section{ Puiseux-index }

Assume $(I)_S$ an ideal of holomorphy, with induced topology of Schwartz-type and with a reduced
pseudo-base. Let $(I)_S^c$ denote elements developed over the pseudo-base in $(I)_S$ with constant
coefficients and $(I)_S^h$ correspondingly with holomorphic coefficients. If $f \in (I)_S^c$ and $h(f)
\sim_0 f$ (geometrically equivalent), then $h(f) \in (I)_S^h$. If also $h(f) \sim_{\infty} f$
(equivalent in strength), then $h(f) \in (I)_S^c$.

Consider now the situation where $f$ is reduced for lineality and $h(f),f$ are entire functions, such that locally
$h(f)=\alpha f$, for $\alpha \in H(\Omega)$, where $\Omega$ is an open set in $\mathbf{C}^{\nu}$. If $h$ is
locally injective, we have again the situation $h(f) \in (I)_S^c$, otherwise we assume $h^N$ injective, for
large $N$. Assume $N=2$, we then have $h(\alpha) \beta$ constant, that is $h(\alpha)=C\beta^{-1}$.
Define $ord_{0,L} F_1(\zeta)=k$, where $\zeta \in L$, $L$ a line-segment through $0$ and where the pseudo-base
is assumed indexed after increasing order of zero's. If $k'=ord_{0,L} F_2(\zeta)$, we have the Puiseux-indexes
for deviation $dist(A,C(A,0)) \sim Cr^{k'/k}$ as $ r \rightarrow 0$ and $\mid \zeta \mid \leq r$
\cite{Chirka}.

\vspace*{.5cm}

Assume $h(f)=\alpha f$ with $\alpha \in H$ and $h^N$ locally injective, for some $N$.
Then $$\frac{1}{(\mid 1 - h \mid)}(f)=\mid \sum_j \alpha_jf + Cf \mid \leq \sum_j \mid \alpha_j \mid
\mid f \mid + C_N \mid f \mid.$$ Thus $\mid f \mid^{-2} \leq 1 + C_N$, for a real constant $C_N$.
If $h^2$ is injective and $f$ reduced, for $\zeta \in Z_{\alpha}$, $1/f=f+c_2f+\ldots$, for constants
$c_2, \ldots$. Thus $f^2=1/(1+C)$, for a $f$. If $h$ is defined over $J^{\epsilon}_{h}$, we can apply
Proposition \ref{c_monotropy}.

\vspace*{.5cm}

Assume $(I)$ an ideal, not necessarily geometric and $(J)=\mbox{ ker }h$, for a homomorphism $h$,
then $N(I \cap J)=\{ \zeta \quad h(f)(\zeta)=0 \quad f \in (I) \}=N_I(J)$, $N_{I_N}(J)=\{ \zeta \quad
h(f^N)(\zeta)=0 \quad f \in (I) \}$. Assume $h$ injective, for $N \geq N_0$, then $N_{I_N}(J)=N(I)$
for $N \geq N_0$. Thus, we have a division in components, $$N_I(J)=N(I) \cup N_{I_2}(J) \cup N_{I_3} \cup \ldots \cup N_{I_{N_0-1}}.$$
Further, $J_2=I(N_{I_2}(J))=\{ g \in (I_2) \quad h(g^2)=0 \}$. We have seen that $V_{c}=\{ \zeta
\quad f(\zeta)=c \}$ and this means that $g \in (J_2) \implies V_{c}(g^2) \neq \emptyset$, for some
$c$ and for all $j < N_0$, the sets $V_{c}(g^j)$ are nontrivial for suitable constants $c$.

\vspace*{.5cm}

Assume now $f \in (I)_S^c$ and reduced and that $h^2$ is injective. Then,
$$\mid \zeta - z \mid^{\sigma} \leq \mid f(\zeta-z) - f(z) \mid + \mid f(z) - f(\zeta) \mid \leq
\mid c \mid + \sum^{\infty}_k \mid f_{\alpha} \mid \mid \zeta \mid^{\mid \alpha \mid}$$ and if $\mid
\zeta \mid < 1$, $$\mid \zeta - z \mid^{\kappa+\sigma} \leq C + C_1 \sum_k^{\infty} \mid f_{\alpha}
\mid,$$ where $\kappa$ is dependent on the modulus of $z$. More precisely, assume $\mid z \mid < 1/R$, for
$R>1$. Then, since $f$ is reduced, there is at least one $\alpha$ such that $1/R^{\kappa} <
C_{\alpha}\mid f_{\alpha}(z) \mid$ and it is clear that $\kappa$ depends on the value of $R$.
 Assume $\mid f_{0} \mid = \max_{\mid \alpha \mid \geq k} \mid f_{\alpha} \mid$ and
$f_0=\sum_j A_j F_j$ for $A_j \in H$, where the pseudobase is reduced and indexed by increasing
order of $0$ in $z$. Then, $$f_0=\sum_jA_j'F_1^{m_j}=(\sum_j A_j'F_1^{m_j-m_0})F_1^{m_0}$$
according to the quantitative version of R\"uckert's Nullstellensatz (cf. \cite{Cart}),
we can find a $m_0$ such that $\mid \sum_j A_j' F_1^{m_j - m_0} \mid \leq C_j$ for all $j$. Thus,
$\mid f_0 \mid^{\frac{1}{\kappa+\sigma}} \leq C \mid F_1 \mid^{\frac{m_0}{\kappa+\sigma}}$ and we can
conclude, for $V$ the zero's to $dh(f)$ over a set $U_{c}$, there are positive constants
$C_1,C_2,\rho,c$ and where $\rho < 1$
$$ \mbox{dist}(\zeta,V) \leq
\mid \zeta \mid^{c \rho}( C_1 \mid F_1 \mid^{\frac{m_0}{\kappa+\sigma}} + C_2)^{\rho} $$

\newtheorem{d_alg}{ Lemma }[section]
\begin{d_alg}
If $h$ is an algebraic homomorphism over an ideal $(I)$ with $(I_N)=\{f^N \quad f \in (I) \}$
such that $h$ is locally injective over $(I_N)$
for $N \geq N_{0}$ and $(J)=\mbox{ ker }h$, then $N(J)$ is locally an algebraic variety.
\end{d_alg}

Note that if $h$ is the extension-mapping defining the Schwartz-spaces, we have seen for our ideals
$(I)$, that $r(I_{weak}) \sim I_{strong}$ and if $h$ is
algebraic and $(I')=h(I)$, $h(r(I_{weak})) \sim r(I_{weak}') \sim (I_{strong}')$.

\newtheorem{d_algh}[d_alg]{ Lemma }
\begin{d_algh}
Assume $X$ a reduced analytic variety and $I(X)=\mbox{ ker }h$ for a homomorphism $h$.
Then $h$ can be selected as algebraic.
\end{d_algh}

Proof:\\
Through the conditions that $X$ is reduced we know that $h$ can be selected as locally injective,
thus within multiplication with a scalar we have that $h$ is locally algebraic. $\Box$

\vsp

We have earlier studied symbols on the form $P(z)=e^{\frac{1}{z}}p_{N}^{\frac{1}{N}}$ with $p_{N}$ reduced.
We have then, for $N$ suitable that $P(\xi) \rightarrow 0$ as $\xi \rightarrow - \infty$. In the same way
for all derivatives to $P$, why $P(-z)$ has an infinite zero in the real infinity. In the same way, $Df(-z)$ has an infinite zero in the real infinity. These remarks can immediately extended to the case $n > 1$.

\newtheorem{d_qube-rep}[d_alg]{ Lemma }
\begin{d_qube-rep} \label{d_qube-rep}
Assume $p_{N}$ reduced in a neighborhood of $a$, $V_{a}$ and $P \equiv 0$ in $V_{a}$, we then have
$\mid P - e^{\frac{1}{a-z}}p_{N}(z-a) \mid < \epsilon$
on a $V \supset V_{a}$
\end{d_qube-rep}

If a constant value is preserved for $P$ in $V$ we also have that $\mid P(z) \mid \leq Q(\frac{1}{\mid z - a \mid})$
in $V \backslash \{ a \}$. We thus have negative indicator for $P$ in $a$. Further using the relation in Lemma
\ref{d_qube-rep} $$e^{-C\frac{1}{\mid z \mid}} \mid P \mid \leq \epsilon e^{-C\frac{1}{\mid z \mid}} + \mid p_{N}(z) \mid$$
in a neighborhood of the origin. The corresponding result as $\mid z \mid \rightarrow \infty$ is that $P$ is of
type $C$.

\newtheorem{d_qube-rep2}[d_alg]{ Lemma }
\begin{d_qube-rep2} \label{d_qube-rep2}
Assume $p_{N}$ as in Lemma \ref{d_qube-rep} and $P^{N} \equiv 0$ on $V_{a}$, then we have $$\mid P-e^{\frac{1}{a-z}}p_{N}^{\frac{1}{N}}(z-a) \mid \leq \epsilon$$
\end{d_qube-rep2}

Consider $I=\int_{\overline{\Gamma}_{C}} e^{\overline{z}} d z$ where the integral is taken
along $\overline{\Gamma}_{C}$ clockwise. This can also be seen as $\overline{\int_{\Gamma_{C}} e^{z} d z}$
where now the integral is taken in counter clockwise direction why $0$ is considered as outer to
$\Gamma_{C}$ in the latter integral. Thus $P$ can on the boundary to a cube with a finite distance to an
infinite zero, be seen as a measure as in Lemma \ref{d_qube-rep}. If the infinite zero corresponds to $P^{N}$ and not $P$, we have the representation as in Lemma \ref{d_qube-rep2}.

\vspace*{.5cm}

In the case with the weighted lineality we define $\Gamma$ as the boundary of a cube such that
$\overline{\mathcal{F}}$ $f$ has point support on $\Gamma$ (finitely many) and $f=const.$ on $\Delta^{Q}_{x}$.
If $N(J)=\Delta^{Q}_{x}$ with $(J)$ as before, we must have $I(\Delta^{Q}_{x})=rad(J)$. Further
$I(\Delta^{Q}_{x})$ is a closed ideal. This means that if $f$ is constant on the weighted lineality,
then $\overline{\mathcal{F}}$ $f$ has point support on $\Gamma$ and conversely.

\newtheorem{d_J-meas}[d_alg]{ Proposition }
\begin{d_J-meas}
Assume $\mu_{\beta}$ are measures with point support on finitely many cubes and $(J)$ the ideal as in section \ref{sec:d_trans}.
Then we have $f \in rad(J)$ $\Leftrightarrow$\\ $\overline{\mathcal{F}}$ $f=\sum_{\beta} D^{\beta}\mu_{\beta}$
\end{d_J-meas}

\subsection{ Estimates }
Assume $J^{\infty}_{V}$ denotes the ideal of
symbols with an infinite zero on a set of positive complex dimension $V$.
Consider phases on the form $h(\varphi)=P/Q$ with $P \in J_{\infty}$ (cf. section \ref{section:d_org}), where $J_{\infty}=\{ g \quad D(f-g)=0 \quad \exists f \in J^{\infty}_V \}$,
that is, $g=f+c$ for a constant $c$ in the infinity, $Q$ is assumed hypoelliptic and
self-adjoint. We know that for $zP(z)=F(\frac{1}{z})$, $F$ has an infinite zero in the infinity
($z \rightarrow 0$). We also consider $P$ such that $z^{\rho}P(z)=F_{\rho}(\frac{1}{z})$, for a real number
$0 < \rho <1$, where $F_{\rho}^{\frac{1}{\rho}}$ is assumed to have an infinite zero in the infinity. This
would be the case for instance if $P^{N} \in J_{\infty}$. We define $I_{\rho}=\{ P \quad z^{\rho}P(z)=F_{\rho}
(\frac{1}{z}) \quad F_{\rho}^{\frac{1}{\rho}} \in J^{\infty}_V \}$ where as before $V$ is a set of positive
complex dimension containing the origin. If we let $V_{r}=\{ z \quad \mid z \mid < r\}$ for $r < 1$ and $F_{\rho}'(\frac{1}{z})=F_{\rho}(\frac{1}{z^{\rho}})$ then if $I_{\rho}'$ corresponds to $F_{\rho}' \in J^{\infty}_{V_{r^{1/\rho}}}$, we have $I_{\rho}'=I_{\rho}$.

\newtheorem{d_r-est}[d_alg]{ Lemma }
\begin{d_r-est}
If $P \in I_{\rho}$ on $V_{a}$, a neighborhood of $a$, we have for $p_{N}$ reduced on $V_{a}$ that
$$ \mid P - e^{{(a-z)}^{-\rho}}p_{N}(z-a) \mid < \epsilon$$
and if $P^{N} \in I_{\rho}$ on $V_{a}$, that for $N \geq N_{0}$
$$ \mid P - e^{{(a-z)}^{-\rho}}p_{N}^{\frac{1}{N}}(z-a) \mid < \epsilon \text{  on  } V_{a} $$
\end{d_r-est}
Define a contour $\Gamma$ according to $\Gamma=\{ \zeta \quad h(e^{\varphi})(\zeta)=0 \}$ and the neighborhood
$\Gamma_{\epsilon}=\{ \zeta \quad \mid h(e^{\varphi})(\zeta) \mid < \epsilon \}$. The R\"uckert Nullstellensatz
now gives that $$ \mid \lambda \mid = \mid e^{h(\varphi)}/Q \mid < \epsilon M$$ on a poly-disc neighborhood of
$\Gamma$. Define the ideal $(I_{\Gamma})$ so that $N(I_{\Gamma})=\Gamma$. Let $I_{\Gamma}(u)=\int_{\Gamma}
e^{\varphi} u d x$. Then there exists a homomorphism $h_{1}$ such that $I_{\Gamma}(h_{1}(u))=0$. Note that
for $\overline{(I_{\Gamma})}$ with $e^{\varphi} \in (I_{\Gamma})$ the closure taken in Whitney-norm we have that $dh(e^{\varphi})=0$ on $\Gamma$
and $e^{\varphi}dh_{1}=0$ on $\Gamma$, for $e^{\varphi} \in \overline{(I_{\Gamma})}$.

\vspace*{.5cm}

If $h$ is not reduced, we denote by $I_{\Gamma,h}$ the integral $\int_{\Gamma} e^{h(\varphi)} u dx$ and
by R\"uckert's Nullstellensatz we can realize $I_{\Gamma,h}$ as a measure in $\Gamma_{\epsilon}$
$\mid I_{\Gamma,h}(u) \mid \leq C_{\Gamma} M \epsilon \sup_{\Gamma} \mid u \mid$. We consider particularly
neighborhoods of $\Gamma$ such as $U_{c}=\{ \zeta \quad \mid \zeta \mid^{c}d(\zeta,\Gamma) < C \}$
for positive constants $c,C$ and the associated sets $U_{c}'=U_{c} \backslash \Gamma$. If $h(\varphi) \in J_{\infty}$
there are $g \in J^{\infty}_{V}$ such that $D(h(\varphi)-g)=0$ on $U_{c}$ and $h(\varphi)=g+a$ for a constant
$a$. We also consider $h(\varphi)=f+g+a$ where $g$ and $a$ are as before and $f$ has negative order $- \rho$,
for $0 < \rho < 1$. Immediately $d(\zeta,\Gamma)$ has negative order on $U_{c}$ and by studying
$1/(\zeta - \zeta_{0})$ where $\zeta_{0} \in \Gamma$ we can compare with indicator theory and we see on
$U_{c}'$ that $\limsup  \frac{1}{r} \log  \mid d \mid = - \widetilde{\rho}$ where $\widetilde{\rho}$ is close
to $1$. We know that if $V$ is an algebraic set, then $U_{c}$ is a domain of holomorphy. Assume $h$ such that
$h^{N}$ is locally injective for $N \geq N_{0}$, we have a decomposition of $U_{c}$
according to $1=P_{k}+R_{k}$, for $k < N_{0}$ with $P_{k}$ polynomial and $R_{k}$ reduced. We have the same decomposition
if the set $\{ h^{N}(e^{\varphi})=c \widetilde{h}^{N}(e^{\varphi}) \} \uparrow U_{c}$ as $N \uparrow \infty$.
Assume $\zeta \in \Gamma \cap V=\Gamma \cap E^{c}$ (the complementary set to the polar set to $h(\varphi)$),
then $\widetilde{h}(e^{\varphi})$ is a polynomial in a neighborhood of $\zeta$. The decomposition of $U_{c}$ gives a decomposition of the distance-function $d(\zeta,\Gamma)=d_{R}(\zeta,E)+d_{V}(\zeta,\Gamma)$ where the terms are extended by zero on the respective complementary set. Note that $d_{V}$ defines a norm. Using an a-priori estimate by Palamodov (cf. \cite{Palamadov}) we have
$$ \mbox{ log }\mid d \mid - c \leq  h(\varphi)  \leq \mbox{ log }\mid d \mid + c$$ and by taking
limes sup. as $r \rightarrow \infty$ we get
$$ - \widetilde{\rho} - c \leq \limsup  \frac{1}{r}  h(\varphi)  \leq - \widetilde{\rho} + c$$
Thus for a sufficiently small neighborhood $U_{c}'$, $e^{h(\varphi)}$ has negative indicator. If $h$ is such
that $h^{N}$ is locally injective for $N \geq N_{0}$ we see that also $U_{c}$ is decreasing by iteration as is $V_{N}$. More precisely we can write $U_{c/N} \downarrow$ as $N \uparrow \infty$ and
$$ - \widetilde{\rho} - c/N \leq \limsup \frac{1}{r}  h^{N}(\varphi)  \leq - \widetilde{\rho} + c/N$$ and for a sufficiently large iteration index, we have that $h^{N}(\varphi)$ has negative order.

\section{ Transversality }
\label{sec:d_trans}

Assume $T$ is a compact operator and of type $0$ over a (reduced) ideal $(I)$.
Then, $S(f)=\frac{T(f)}{f}$ is of negative type over $(I)$ (of type $-\infty$).
Let $\mathcal{K}$ $=\{ g \quad T(g)-R(g) \text{ of type 0 } \quad R \text{ of type } -\infty \}$.
For $g \in \mathcal{K}$, we have $\frac{T(g)}{g}$ is of type $0$ and also of type $-\infty$.
Assume, $T^{N} - I \sim_{m} 0$. This means that there is an entire function $R$, such that modulo monotropy,
$T^{N}-I=R$ in $Exp_{0}$.

\newtheorem{d_appr}[d_alg]{ Lemma }
\begin{d_appr}
Assume $R$ entire and $R \sim_{m} 0$, then there exists a $R_{1}$ of type $-\infty$ such that $R-R_1$ is of type $0$.
\end{d_appr}

This means that for $K$ a compact operator, such that $K \sim_{m} 0$, we have \\$\mathcal{K}$ $=(I)$ and for a compact operator such that $K^{N} \sim_{m} 0$, from some
$N \geq N_{0}$, we have that $\mathcal{K}$ $=rad(I)$.

\subsection{ Transversality modulo monotropy }

Assume $h(g)=w\mathcal{F}(\textit{g})-\mathcal{F}\textit{w(g)}$ for $g \in J$ and given a homomorphism $w$, such that $w$ is locally factorized by $\tau_{\epsilon} v$, for an algebraic homomorphism $v$.

\newtheorem{d_rep}[d_alg]{ Lemma }
\begin{d_rep}
For any $g \in (I)$ an ideal in $H$ and a homomorphism $w$ defined on $(I)$,
there is a $g_1 \in H$ such that $w(g)=\widetilde{w}(g_1)=I_{\big[ W,I \big]}(g_1)$
\end{d_rep}

Thus $g \in (I)$ can be realized through $Id \in H'$ and $Id(g_1)=g$. The composition of analytic functionals,
here realized through measures, is given as convolution of measures. $Id(\mu)(g_{1})=\int_{V} g_{1} d \mu$ and
$Id(\mu)^2(g_{1})=Id(\mu * \mu)(g_{1})$. We write $g^{2}=g_{1}*g_{1}$ for short. Let
$h'(g^{2})=$\\ $v \mathcal{F}$ $g^{2}-\mathcal{F}$ $v(g^{2})=v^2 \mathcal{F}$ $(g_{1})-\mathcal{F}$ $v^2(g)$.
Further, we note that $h^{'2}(g)=$ $(v \mathcal{F}-\mathcal{F}$ $v)^{2}(g)=$ $(\mathcal{F}$ $v)^{2}(g)-v \mathcal{F}\mathcal{F}$ $(g)$.

\newtheorem{d_rep2}[d_alg]{ Lemma }
\begin{d_rep2} \label{d_rep2}
For any $g \in (I)$ and a homomorphism $v$ acting on $(I)$, there exists a
$g_{0}$ such that $\mathcal{F}$ $v^2(g_{0})=-(\mathcal{F}\textit{v})^{\textsl{2}}$ $(g)$.
\end{d_rep2}
The proof is analogous to the proof of the previous Lemma. Finally, assuming as above that
$v \mathcal{F}$ $(g)=\mathcal{F}$ $v(g) + c$, for a constant $c$. If $(\mathcal{F}$ $g_{1})^{2}=$ \\ $\mathcal{F}$ $(\mu * \mu)(g)=$ $\mu*\mu(\mathcal{F}$ $(g_{0})) + c$, where $\mu$ corresponds to the functional $Id$, thus there is a $g_{2}$ such that $-\mathcal{F}\mathcal{F}$ $v(g)=v\mathcal{F}$ $g_{2}+ c$, for a constant $c$. Thus, we could say that $d (h'(g^{2})-h^{'2}(g))=0$ and that $h'$ modulo monotropy and within representation in $H'$ is an algebraic homomorphism in the tangent space.

\subsection{ A lifting principle }

We can define a vector-field along a morphism $h$ according to
\begin{displaymath}
   \left. \begin{array}{ll}
   X \underrightarrow{\phi}   &  T(Y) \\
   h   \searrow & \downarrow{\pi} \\
                & Y
   \end{array} \right.
\end{displaymath}
We then know that $\phi^{*}=\phi$ is equivalent to the proposition that $dh=hd$. If this
is the case we have for our application that $d\psi$ and $\psi$ are not in the respective "eigen-
spaces" $(I_{\mu}),(I_{\frac{1}{\mu}})$ on the domain $\Omega$. 

\vspace*{.5cm}

Consider the intersection between the transversals and the leaves to the foliation ( modulo monotropy
). We know that if
\begin{equation} \label{d_Poincare}
P/Q -y/x
\end{equation}
 changes sign in a domain $\Omega$, then this intersection is not
discrete (that is a characteristic $\gamma$ to the system runs partially in this intersection).
The proposition is thus that there is a $\psi$ non-constant such that
$$ \frac{h(\psi)}{\psi}=\mu \qquad \frac{d h(\psi)}{d \psi}=\frac{1}{\mu}$$
If $h$ is such that $h^{2}$ is locally injective ( modulo monotropy algebraic ), we can write (for
all $j$) $ \frac{d}{d \zeta_{j}}\big( h^{2}(\psi)-\psi^{2} \big)=0$ or locally
$\frac{d}{d \zeta_{j}}\big( c \psi - \psi^{2} \big)=0$, but then $\psi$ would be constant in
contradiction with the conditions.
We conclude that

\newtheorem{d_fol}[d_alg]{ Proposition }
\begin{d_fol}
Given a homomorphism $h$ algebraic modulo monotropy and such that $h^{2}$ is locally injective
then the intersection between the transversals and the leaves to the foliation is discrete.
\end{d_fol}

\subsection*{ Fredholm operators }
Differential operators have a natural interpretation as Fredholm-operators.
We consider the standard projections $P: L^2 \rightarrow R(L_{\lambda})$,
$Q:H^{s,t} \rightarrow N(L_{\lambda})$, where $H^{s,t}$ is the Sobolev-space with tensorized weights and
$L_{\lambda}$ the differential operator $L(D)-\lambda$, for a real or complex parameter $\lambda$.
$L_{\lambda} \in \Phi(H^{s,t},L^2)$ ($\Phi$ denotes the class of Fredholm operators) gives a
decomposition $$H^{s,t}=X_0 \oplus N(L_{\lambda}) \text{ and } L^2=Y_0 \oplus
R(L_{\lambda}),$$ where $N(L_{\lambda})$ denotes the solutions to the homogeneous equation
and $R(L_{\lambda})$ the range of $L_{\lambda}$. It is not difficult to prove that, for
$\sigma > 0$, $L_{\lambda}$ is homogeneously $L^2$-hypoelliptic, which means
that $N(L_{\lambda})=Q(H^{s,t}) \subset C^{\infty}$, that is $Q$ is
regularizing on $H^{s,t}$ and $E_{\lambda}$ is a left parametrix to
$L_{\lambda}$. If $L_{\lambda}$ is homogeneously $L^2$-hypoelliptic, then also its Hilbert-space adjoint is
homogeneously $L^2$-hypoelliptic, that is $P^{\bot}(L^2)=N(L_{\lambda}^{*}) \subset C^{\infty}$,
and $P^{\bot}$ is regularizing on $L^2$. We conclude that $E_{\lambda}$ is a left and right
parametrix to the operator $L_{\lambda}$.
Also, if
$\sigma$ arbitrary, that is $L_{\lambda}$ not necessarily
homogeneously $L^2$-
hypoelliptic, we can extend the definition of $E_{\lambda}$ with a
regularizing term and we get a hypoelliptic action on, at least part of
$L^2$. If $E_{\lambda} \in \Phi$, is on the form $I-K$, then there is a
positive integer $N_0$ such that $N(E_{\lambda}^N)=N(E_{\lambda}^{N_0})$ for all $N \geq N_0$.  Since adding a compact
operator to the Fredhom-inverse, does not change the index or the form $I-K$, we let $E_{\lambda}$ be
  defined as
\begin{displaymath}
      \left \{
\begin{array}{lr}
  L_{\lambda}^{-1} + X \text{ on } R(L_{\lambda}) \text{ with } X \in C^ {\infty} \\
  X' \in C^{\infty} \text{ on }  R(L_{\lambda})^{\bot} \cap N(E_{\lambda}^M)^{\bot}  \text{ for a suitable } M \\
  f \in L^2  \text{ otherwise }.
  \end{array} \right.
\end{displaymath}
Thus, this operator $L_{\lambda}$ is hypoelliptic on
   $N(E_{\lambda}^M)^{\bot}$.

\section{ Analytic hypoellipticity in $L^2$ for the iterated operator }
\label{sec:b_wf-a}
We shall in this section only consider constant coefficients, differential operators, with $D_{x}$ replaced by $(i D_{x})$.
Consider the translation operator on $H(E_{ \mathbf{R}})$, $a_{\Gamma}u(\xi_x)=u(\xi_x + i t \eta_0)$, for $t \eta_0 \subset
\Gamma$. Let $$A_{\Gamma}=  \mathcal{F}^{-\textsl{1}} \textit{a}_{\Gamma}  \mathcal{F} \text{:} \textit{H(E}_{ \mathbf{R}}{)} \rightarrow
\textit{H}'(\textit{E}_{ \mathbf{R}}{)}.$$ Here, $ \mathcal{F}$ denotes the Fourier-Borel transform and $\Gamma$ some simple cone in $\Delta_{ \mathbf{C}}$
and $E_{ \mathbf{R}}$ denotes in this context some subset of $ \mathbf{R}^{\nu}$ or the  entire $ \mathbf{R}^{\nu}$.
If $\eta_0 \in \Delta_{ \mathbf{C}}(P)$, we have $P(i D)A_{\Gamma}u=A_{\Gamma}P(i D)u$. Consider the real
quotient homomorphism, $\varphi:$ $E_{ \mathbf{R}} \rightarrow$ $E_{ \mathbf{R}} / i \Delta_{ \mathbf{C}}$.
Let $a_{\varphi(\Gamma)}u(\xi_x)=u(\varphi(\xi_x))$, we then have $b_{\varphi(\Gamma)}=\lim_{t \rightarrow
+0}a_{\varphi(\Gamma)}$. We define $B_{\Gamma}=\lim_{t \rightarrow
+0}A_{\Gamma}$.
 The definition of $b_{\Gamma}$ (cf. \cite{Sjo}), is not dependent on the neighborhood of 0 when we argue
that this means that $b_{\Gamma}$ is quasi-portable by $E_{ \mathbf{R}} \times i 0$. The same arguments give that $b_{\varphi(\Gamma)}$
is quasi-portable by $E_{ \mathbf{R}} \times i \Delta_{ \mathbf{C}}$.

\vsp

For $u \in L^2$, we can make a composition, $u=u_+ + u_-$, where $\widehat{u}_-$, has support contained in cones
with indicator $\delta_-(\eta_0) < 0$. Note that the indicator for a cone $\Gamma$ is here defined as
$\delta(\eta_0)= \sup_{ y \in \Gamma} < y , \eta_0 >$. According to ( \cite{Sjo}), $u_-$ can be considered as
regular elements and it is sufficient to consider $u_+=u-u_-$.
We argue that, to conclude that the operator $P(i D)$ is analytically hypoelliptic on $L^2$, it is sufficient
to study the set of complex lineality and the action on $u_+$- terms.

\vsp

Assume $h_F(\eta)$ the indicator of growth, corresponding to $B_{\Gamma}$. Using  \cite{Mart}(Ch.2,Theorem 4.3), it
can be shown, that for $\eta \in \Delta_0=\Delta_{ \mathbf{C}} \backslash 0$, $h_F=g$, where $g$ is the indicator of growth corresponding to the frequency component
$\Sigma$ in $\mbox{WF}_A(u)$. Thus any cone in $\Delta_0$ has indicator $\geq 0$.
More precisely, let $W$ denote the convex hull of the real support
for $B_{\Gamma}$, that is $W=\{ y; <y,\eta> \leq h_{F}(\eta) \
\parallel \eta \parallel=1 \}$. Let $$W_+=\{ y \in W; <y,\eta> \geq
0 \parallel \eta \parallel=1 \}$$ and let $W_{-}$ be the complementary set. Let $$V_+=\{ \eta; <y,\eta>
\geq 0 \quad y \in W_+\}.$$ Further, $\Delta_0 \cap V_+= \{ \eta; <y,\eta>=0
\ y \in W_+ \}$. Since $g(\eta)=h_F(\eta)=0$ over $\eta \in
\Delta_0$, we must have $\Sigma \cap \Delta_0 \cap V_-= \emptyset$. Next,
using  \cite{Mart} (Ch.2, Prop. 5.2) and conjugated indicators,
we can prove a support theorem in $H'(E_{ \mathbf{R}})$, that is

\newtheorem{b_support}{Proposition }[section]
\begin{b_support}
If $P_{\lambda}(i D)$ is a constant coefficients differential operator and
$u \in H'(E_{ \mathbf{R}})$, $P(iD)u \in L^2(E_{ \mathbf{R}})$, we have in $H'$, $\mbox{ supp }P(i D)u=\mbox{ supp }u$.
\end{b_support}

Let $\Sigma$ be such that $b_{\Gamma} \neq 0$, $\Sigma'$ the frequency component
corresponding to $\mbox{ WF}_A(P_{\lambda}u)$. Knowing that $\Sigma' \subset \Sigma$ and $\Delta_0 \subset \Sigma \backslash \Sigma'$ (that is
a necessary condition for hypoellipticity is that the operator is reduced for complex lineality), we have $$\Sigma=\Delta_0 \cup \Sigma' \backslash \Delta_0 \text{ and } \Delta_{ \mathbf{C}} \cap \Sigma'=\{ 0 \}$$ and the
lineality can be regarded as the "boundary" for $\Sigma$.
We claim that cones $\Gamma_j$ exist, such that $\Delta_0 \subset \cup_j \Gamma_j \subset \Sigma$. If we assume the operator $P_{\lambda}$, such that $\mbox{ Re }P_{\lambda} \sim
P_{\lambda}$, it will be sufficient to study the real part of the polynomial, which is considered as a real analytic
function. This means that a proposition (cf.  \cite{Mart} ch2, sec.4,Cor.2) can be applied on $B_{\Gamma}$
which gives that if $B_{\Gamma}$ is quasi-portable by $E_{ \mathbf{R}} \times i \Delta_{ \mathbf{C}}$, then it is portable by the same set.
Thus, we  get a corresponding
representation $B_{\Gamma}=\sum_j B_{\Gamma_j}$. We also note that $b_{\Gamma}$
is in $L^2$ strictly portable by $E_{ \mathbf{R}} \times i \Delta_{ \mathbf{C}}$, because of the properties of the restriction-mapping
in this space.

\newtheorem{b_WF-reg}[b_support]{Lemma }
\begin{b_WF-reg}
 Assume $P_{\lambda}$ with constant coefficients and self-adjoint and that
 $P_{\lambda}$ is reduced with respect to the set of lineality, $\Delta_{ \mathbf{C}}$, $P_{\lambda}u \in L^2$ for some distribution $u$ and such that $\mbox{ WF}_A(P_{\lambda}u)=\emptyset$.
 Then $b_{\Gamma}(u)$ is strictly portable
 by open connected cones with indicator $< 0$.
\end{b_WF-reg}

Proof:
Assume contrarily, that $L_+$ is a line in $\Gamma_+$, an open connected cone with indicator $\geq 0$ and that $b_{L_+}(u) \neq
0$. We can write
\begin{equation}
\label{b_b-gamma}
P_{\lambda}b_{L_+}(u)=b_{L_+}(P_{\lambda}u) + C_{P_{\lambda},L_+}(u)
\end{equation}
where $C_{P_{\lambda},L_+}(u)= \big[ P_{\lambda}b_{L_+} - b_{L_+}P_{\lambda} \big](u)$. We can assume that $u$ is such that $\widehat{u} \neq 0$
on $L_+$. We have either $C_{P_{\lambda},L_+}(u)=0$, in which case the reducedness condition implies $L_+=0$ or $C_{P_{\lambda},L_+}(u) \neq
0$, which contradicts that $\Delta_{ \mathbf{C}} \cap \Sigma'=\{ 0 \}$.$\Box$

\vsp

Consider the canonical homomorphism $ \mathcal{O}_{\textsl{X}} \rightarrow  \mathcal{C}_{\textsl{X}}$, which assigns to sections $s \in
 \mathcal{O}_{\textsl{X}}(\textsl{U})$, their continuous functions $[s] \in  \mathcal{C}_{\textsl{X}}(\textsl{U})$. According to  \cite{Gra}, as a consequence of the R\"uckert Nullstellensatz, the kernel to this homomorphism, is the nilradical $ \mathcal{N}_{\textsl{X}}$.
Using  \cite{Mart} and  \cite{Gra}, the homomorphism can be regarded as an evaluation functional on $H_U$. That is
$E(\varphi)(x)=\varphi(x)$, for $\varphi \in H_U$. Using the nuclearity of $H'_U$, we could write
$I_E(\varphi)(x)=\varphi(x)$. The cited result means that $I_E^N(\varphi)(x)=0$, for $N$ large and $\varphi \in
H_U$ and iteration in this context corresponds to composition of convolution kernels. More
precisely, let $\Delta \in H'( \mathbf{C}^{\nu} \times  \mathbf{C}^{\nu})$, be the kernel
representing $\mbox{ Id }_x$ ($x \neq 0$). It follows from the argument below, that in
a neighborhood of $x$, $V$, we have for this kernel,
$$<I^N_{\Delta}(\varphi),\psi >_V=<I_{\big[ \Delta,\Delta
\big]_N} (\varphi),\psi >_V=0, \text{ with } \varphi,\psi \in
H( \mathbf{C}^{\nu}), \text{ for some } \textit{N}.$$

\vsp

Define an ideal  $ \mathcal{J}$, as holomorphic functions, invariant for complex
translation with $t\eta$, $\eta \in \Delta_{ \mathbf{C}}$ and $t\eta \in \Gamma$, that is $ \mathcal{J}$ $=\{ \varphi \in
 \mathcal{O}$ $(U)$ ; $a_{\Gamma}\varphi=\varphi \}$. We then have $$V_1=U_1 \backslash \Delta_{ \mathbf{C}}
\subset N( \mathcal{J})=\textit{U}_{\textsl{1}}=\textit{N}( \mathcal{N}_{\textsl{U}_{\textsl{1}}}).$$ The construction can be
iterated, $ \mathcal{J}^{\textsl{2}}=\{ \varphi \in  \mathcal{O}$; $a_{\Gamma}\varphi^2=\varphi^2 \}$, with $V_2=U_2 \backslash \Delta_{ \mathbf{C},\textsl{2}} \subset
U_2$, where $\Delta_{ \mathbf{C},\textsl{2}}$ is the lineality corresponding to the operator $P^2$. Note that since $U_j=V_j \cup
\Delta_{ \mathbf{C},\textsl{j}}$, we have $ \mathcal{J}(\textsl{U})= \mathcal{\textsl{J}_{\textsl{1}}}  \mathcal{\textsl{J}_{\textsl{2}}}$ and $ \mathcal{\textsl{J}_{\textsl{2}}}$
$=\{ \varphi \quad a_{\Gamma} \varphi= \varphi \}$ $N( \mathcal{\textsl{J}_{\textsl{2}}^{\textit{j}}})$ $=V_j$. Thus, as $a_{\Gamma}$ is a homomorphism,
 we have that if $\Delta_{ \mathbf{C},\textsl{j}} \downarrow \{ 0 \}$, then also $V_j$ are decreasing, as the index
grows. The sequence $U_N$, constructed in this way, is thus for a partially hypoelliptic operator,
decreasing (with respect to inclusion of sets) and $0 \notin V_N$, for all $N$, but $P \in
 \mathcal{J}^{\textsl{N}}_{\textsl{0}}$, for all $N$. From the construction, $$N( \mathcal{J}^N)=\textsl{N}( \mathcal{N}_{\textsl{U}_{\textsl{N}}})
\text{ which implies } \mbox{ rad } \mathcal{J}^{\textsl{N}}= \mathcal{N}_{\textsl{U}_{\textsl{N}}}, \text{ for all } \text{N}.$$ Using that $\Delta_{ \mathbf{C},\textsl{N}}$ is an analytic set,
 for $U_N$ sufficiently small,the corresponding nilradical will be the zero ideal. The conclusion must be that $P \in  \mathcal{J}^{\textsl{N}}_{\textit{t}}$
implies that $t=0$ and the iterated operator $P^N$, is in this sense reduced.

\newtheorem{b_transl}[b_support]{ Proposition }
\begin{b_transl}
In the case of several dimensional translations, (product topology) we have that if $P$ is a
partially hypoelliptic operator with constant coefficients, then there is a $N_0$ such that for
$ N \geq N_0$, $P^N$ is reduced for complex lineality.
\end{b_transl}

\section{ $b_{\Gamma}$ considered on $\mbox{ Exp }_{\parallel \cdot \parallel,A}$}
\label{sec:b_b-gamma}
Assume, that for a holomorphic function $g$

\begin{equation} \label{b_Hartogs}
\mid g(\zeta) \mid \leq C e^{B \mid \mbox{ Im } \zeta \mid + A \mid \mbox{ Re } \zeta
\mid} \qquad \zeta \in  \mathbf{C}^{\nu}
\end{equation}

for a positive constant $C$.  If the real and positive numbers $A,B$ can be
chosen arbitrarily small, the corresponding operator allows real support. According to  \cite{Mart},
if $g$ is entire in $ \mathbf{C}^{\nu}$ and $\mid g(i y) \mid \leq C(\epsilon) e^{{\epsilon}
\parallel y \parallel}$, $\mid g(x) \mid \leq C'(\epsilon) e^{(A+\epsilon) \parallel x \parallel}$,
for $x,y \in   \mathbf{R}^{\nu}$, then $\mid g(x+iy) \mid \leq C''(\epsilon)e^{(A+\epsilon) \parallel x+iy \parallel}$, where $A$ is
positive number, determined by the indicator, corresponding to $g$ and $\parallel \cdot \parallel$ is a norm in $ \mathbf{C}^{\nu}$.
This means that it is possible to consider $b_{\Gamma}$ in $\mbox{ Exp }_{\parallel \cdot
\parallel,\textit{A}}$.
The following limits can be determined.

$$
\lim_{y \rightarrow 0} \big[ e^{-(A+\epsilon)\parallel x+iy \parallel} g(x+iy) - e^{-(A+\epsilon) \parallel
x+ iy
\parallel}g(x) \big]$$
$$\lim_{y \rightarrow 0} \big[ e^{-(A+\epsilon) \parallel x+iy \parallel} g(x)- e^{-(A+\epsilon) \parallel x
\parallel} g(x) \big]$$

which means, (cf. the Lemma below), that
$$b_{\Gamma} \big[ e^{-(A+\epsilon) \parallel x \parallel} g(x) \big] =
e^{-(A+\epsilon)\parallel x \parallel} g(x)
\text{ in } \mbox{Exp}_{\parallel \cdot \parallel,\textit{0}}$$
The corresponding operator can now be considered as portable by a ball of radius $A$ with respect to
some complex norm or by a tube \\ $\{ x \in  \mathbf{R}^{\nu}$ $, \mid x_j \mid < A+\epsilon' \} + i
 \mathbf{R}^{\nu}$

\newtheorem{b_e-type}{ Lemma }[section]
\begin{b_e-type} \label{b_e-type}
Given an entire function $g$ of finite type $A$, the entire function  $g(\zeta)-g(x)$, $\zeta=x +
iy$, is of (complex) type $0$,that is, for all $\zeta \in B_{\epsilon} + i  \mathbf{R}^{\nu}$, there is a positive constant $C$,
$$ \mid e^{-\epsilon \parallel \zeta \parallel} \big[ g(\zeta) - g(x) \big] \mid \leq C$$
\end{b_e-type}

Proof: Let $M_{\alpha}(t)=\sup_{\parallel \eta \parallel \leq t} \mid g_{\alpha}(i \eta ) \mid$. In general,
if $g_{\alpha}$ is entire for all $\alpha$ (analogous to Weierstrass M-test) we have $M_{\alpha}^{1/\mid
\alpha \mid}(t) \rightarrow 0$ as $\mid \alpha \mid \rightarrow \infty$, for every $t$. This means that
 $\sum_{\alpha} g_{\alpha}(i\eta)x^{\alpha}$ is an entire function of
$\zeta=x+i \eta$. Let
$$\Gamma_{g}=\{ \eta, e^{-\epsilon\parallel x \parallel} \big[ g(x+it\eta) - g(x)
\big]=0 \ \forall x \ \forall t \}=\Delta_{ \mathbf{C}}(g)$$
Then, for $\eta \in \Gamma_{g}$, $g_{\alpha}(it\eta)=C_{\alpha}$,
a constant independent of $t$. Thus $M_{\alpha}^{1/\mid \alpha \mid}(t) \rightarrow 0$ over
$\Gamma_{g}$ as $\mid \alpha \mid \rightarrow \infty$, for all $t$. Assume $\eta$ close to $\eta_0 \in \Gamma_{g}$. Then,
$g_{\alpha}(it\eta)=\sum_{\beta}g_{\alpha \beta}(it\eta_0)(it(\eta - \eta_0))^{\beta}$.
Now, to see that $g_{\alpha}$ is of exponential type $0$, let
$M_{\alpha\beta}(t)=\sup_{\parallel \eta \parallel \leq 1} \mid g_{\alpha \beta}(it\eta) \mid$ and
consider
$$\frac{\log{M_{\alpha\beta}(t)}}{t} \leq \epsilon
\mid \alpha \mid \ \forall \alpha$$
 The left side has the indicator for $g_{\alpha\beta}$ as limes superior, as $t \rightarrow
\infty$. As $\mid \beta \mid \rightarrow \infty$, the limit on the left side, times $1/\mid \beta \mid$, is $0$, for all $t >
1$. Thus, $g_{\alpha}$ is holomorphic and of type $0$, in any conical neighborhood of $\Gamma_{g}$. $\Box$.

\vsp

This may be regarded as a generalization of the proposition that some iterate of the polynomial defining the lineality to a symmetric operator,
is a hypoelliptic polynomial. Given $g$ of type $A$, we then have that $e^{-A\parallel x
\parallel}b_{\Gamma}(g) \in \mbox{ Exp }_{\parallel \cdot \parallel,0}$.

\vsp

Concerning self-adjointedness, the following propositions are immediate.

\newtheorem{b_self-a}[b_e-type]{ Proposition }
\begin{b_self-a}
Assume $b_{\Gamma_{+}}f=g_A^{+} \in \mbox{ Exp }_{\parallel \cdot \parallel,\textit{A}}$.
A sufficient condition for self-adjointedness of $I_E$ is that $2E=g^{+}_A + \overline{g^{+}_A}$. A
necessary condition, is that $\mbox{ Im }g_A^{+}=b_{\Gamma_{-}} \mbox{ Im }f$.
\end{b_self-a}

It is a consequence of the previous argument, that if $f$ is of type $0$, then the corresponding
operator will be self-adjoint in $Exp$. It is also a consequence, that hypoelliptic operators can be
considered as of (complex) type $0$ (consider for instance them as acting on parametrices to hypoelliptic polynomial operators).

\section{ Some remarks on Weyl-calculus }

In any theory involving the strong Fourier-transform, particular care has to be taken with the sets
$V_{c}$. Consider for instance Weyl-calculus. The fundamental relations are
$$ a(x,\xi)= \int K(x+t/2,x-t/2)e^{-i <t,\xi>}d t$$
$$ K(x+t/2,x-t/2)=\frac{1}{(2 \pi)^n} \int a(x,\xi)e^{i <t,\xi>} d \xi$$
where $a \in S^m$ and the operator $a(x,D) \textsl{:} \mathcal{S} \rightarrow \mathcal{S}$.
Assume $a(x,t)=a(x,0)$ for $t$ on a line. Then, $K(x+t/2,x-t/2)=K(x,x)$ or
$\tau_{t/2}I_K\tau_{t/2}=I_K$. This situation may occur if $h(a)(\zeta)=0$, for a radically
injective homomorphism $h$ and $N(I_h)$ corresponds to a generalization of the lineality. On these
lines $K(x,x)=c \delta_{\textsl{0}}$ (the Dirac-measure) and $D^{\alpha}_ta(x,t)$ is a polynomial,
for $x$ fixed.

\newtheorem{c_RED}{ Lemma }[section]
\begin{c_RED}
Assume $F$ a reduced kernel in $L^2 \cap H$, then\\
 $\mathcal{F}^{-\textsl{1}}$ $F \in L^2 \cap H$
\end{c_RED}
Otherwise, if $X=X' \backslash \Gamma$ reduced, for a line $\Gamma$ and $X' \subset \mathbf{R}^{\textit{n}}$,
we have $\mathcal{F}^{-\textsl{1}} \textsl{F}$=
$c \delta_{\Gamma} + \mathcal{F}^{-\textsl{1}}\big[ \textsl{F} \mid_{\textsl{X}} \big]$, that is $\mathcal{F}^{-{\textsl{1}}}
\big[ \textsl{F} -$ $I \big] \in L^2 \cap H$. Assume now, $a(x,\xi)=\sum_j a_jf_j=\sum_j I_{\big[ A_j,F_j \big]}$
over $L^2 \cap H$, with $F_j$ reduced in both $\xi$-variables separately. Then $\mathcal{F}^{- \textsl{1}}$ $a(x,\xi)$
$=\sum_j I_{\big[ A_j, \mathcal{F}^{-\textsl{1}} \textsl{F}_{\textsl{j}} \big]}$ and $ \big[ A_j, \mathcal{F}^{-\textsl{1}}
\textsl{F}_{\textsl{j}} \big]$ $\in L^2 \cap H$. Further, $I_K=\sum_j I_{\big[ \tau_{-t/2}A_j\tau_{-t/2},\mathcal{F}^{-\textsl{1}} \textsl{F}_{\textsl{j}}\big]}$, that is
$$ \frac{1}{(2 \pi)^n} \int \int A(x - t/2,y + t/2)F(y,\xi)e^{i <t,\xi>} d \xi d y$$
If $A_j$ is constant on $t/2$-lines, then $I_K=\sum_j c_jI_{\big[ \delta,F_j \big]}$. This
phenomena does not occur for reduced symbols (the coefficients reduced), since the sets $V_{c}$
then reduce to points. If thus $a \in \overline{J}_{h}$ and $a \sim_{0} \sum_j a_j f_j$, with $f_j$
reduced and $a_j$ constants, we have,
$$ \tau_{t/2}I_{K}\tau_{t/2}=\sum_j a_j I_{\mathcal{F}^{-\textsl{1}} \textsl{F}_{\textsl{j}}}=\mathcal{F}^{-\textsl{1}} {\textit{a}}$$
with $\mathcal{F}^{-\textsl{1}}$ $F_j \in L^2 \cap H$. Note again that the sum does not have to be
reduced. Also, note that in analogy
with earlier results, we have "locally" a polynomial representation $a^w_t-p_t=H_t$, with $H_t$
regularizing and $p_t$ polynomial.

\section{ Link between the dynamical system and the wavefront set }

Given existence of $\psi$ not constant such that $h(\psi)=\mu \psi$ and if also $Q=\mu P$ for a constant $\mu$ and where $P,Q$ are polynomial right hand sides to the system under consideration, we have that
$H_{h}(\psi)=0$ why $H_{h}$ has a real integral curve. Assume $W_{c}=\{ (\eta_{1}(\zeta),\zeta) \quad \zeta \in V_{c} \}$
where $V_{c}$ is the foliation corresponding to the symbol to a differential operator with analytic coefficients
and $\eta_{1}(\zeta)$ is a real-analytic function (assume order one). Assume $\Gamma$ a tangent-line to the foliation
(restriction to one variable) and $\Gamma \rightarrow L$ mapping dual line-segments in $Exp_{\rho,A}$ with respect
to the Fourier-Borel transform. Let $W_{\Gamma}$ be the
corresponding set $\{ (\eta_{1}(\zeta),\zeta) \quad \zeta \in \Gamma \}$. For the situation with $\zeta \in V_{c}$
there are corresponding functions $\psi \in (I_{c})$
and a corresponding tangent equation $Q-cP=0$ where $P,Q$ are assumed polynomials. Let $t$ be the mapping $(I_{c}) \rightarrow \Gamma$,
where we can obviously assume that singular points to the system are mapped on foliation points and that regular
points are mapped on the outside of the foliation. Further, for solutions
$u$ to the homogeneous equation to the differential operator above, it is not
difficult to find a mapping $ v \text{:} W_{\Gamma} \rightarrow WF_{a}(u)$.
Denote with $\widetilde{\gamma}$ the image of $\gamma$ under the mapping $\gamma \rightarrow W_{\Gamma} \rightarrow WF_{a}(u)$, then according to Hanges' theorem
(cf. \cite{Sjo}), there is a neighborhood of the origin $U$ such that
$ \widetilde{\gamma}(U) \subset WF_{a}(u)$.

\vspace*{.5cm}

Starting with a complex line $\Gamma$, we consider the transversal $L$ as linearly independent.
$\Gamma$ corresponds as before to an ideal $(I_{\mu_{1}})$, for a real number $\mu_{1}$. Let $(I_{\mu_{1}})'$
be the ideal corresponding to the transversal $L$. In a neighborhood of a singular point to the
associated dynamical system, we can in this way construct a symplectic base. If
$(I)=\bigoplus_{\mu}(I_{\mu})$ (finite sum) and if $\eta$ is cyclic we have thus $\eta(\psi)=\phi$ for
a $\phi \in (I_{\mu_{1}})'$, that is $\eta$ maps tangents on transversals and we could claim that
the transversal is locally unique. Note that in the case where $\psi=e^{\varphi}$ if also
$\frac{d}{d \varphi} h = h \frac{d}{d \varphi}$, then $(I_{\mu})=(I_{\mu})'$, but we are still
assuming polynomial right hand sides for the associated dynamical system, why there are only finitely
many tangents and we will denote half of them transversals and also in this case we can construct a
local symplectic base. More precisely assume $h$ algebraic such that $h(g^{2})=h(g)^{2}$ for $g$ real
and $h^{*}=h$. The corresponding
eigen-values must be real and positive, that is $h(g^{2}) \geq 0$ and $h(g^{2})=cg^{2}$ with $c$
real and non-negative. Assume $c=\mu^{2}$ for $\mu$ real, then the "eigen-spaces" can be seen as
linearly independent $(I_{\mu})_{+} \bigoplus (I_{\mu})_{-}$. We can set
$(I_{\mu})'=(I_{\mu})_{+}$. This means that $\eta$ maps tangents $(I_{\mu})_{-}$ on transversals
$(I_{\mu})_{+}$ or $(I_{\mu}) \rightarrow (I_{-\mu})$. Further we have argued that the decomposition
of tangents and transversals is linearly independent ($\eta_{1}$ is assumed of order $1$ and
locally injective). Thus $t$ and $v$ maps
$$ (I)_{+} \bigoplus (I)_{-} \rightarrow W_{\Gamma} \bigoplus W_{L} \rightarrow WF_{a}$$

\section{ Multi-valuedness}
We note that if $\gamma$ is defined by an algebraic homomorphism $h$, then multi-valuedness for
$h(\psi)$ corresponds to multi-valuedness for $\psi$, since the corresponding Puiseux-series is
mapped on a Puiseux-series. Consider $\Sigma=\{ P=0 \}$ for a polynomial $P$ such that $d_{\zeta}P \neq 0$
on $\Sigma$ and $F(\psi)(\zeta)d_{\zeta}P(\psi)(\zeta)=\sum_{j}b_{j}(\zeta)\psi^{j/\mu}$ where $F$ is assumed
holomorphic with respect to $\Sigma$ and $\psi$ is assumed reduced ( assume a pseudo-base element ). For $\phi=\psi^{1/\mu}$, we get
$h(\phi)^{\mu}=h(\psi)$ and if $h$ maps reduced (pseudo-base) elements on reduced (pseudo-base) elements, we get an expansion
$\sum_{j}h(b_{j}\phi^{j})=\sum_{j}h(b_{j})h(\phi)^{j}$ with multi-valuedness of the same order as the
given set. Consider
$$ P_{c}(\psi)=\prod_{\zeta \in Eq}(c \psi - h(\psi))$$
where $h$ is not assumed reduced. If $h^{2}$ is locally injective we consider the locally discrete set
$$Eq=\{\zeta \quad h(\psi)(\zeta)=h(\psi)(\zeta_{0}) \}$$
and the corresponding Puiseux-expansion $\sum b_{j}(\zeta)\phi^{j}$. This is mapped by $h$
on another Puiseux-expansion and a corresponding discrete set
$$Eq'=\{ \zeta \quad \psi(\zeta)=\psi(\zeta_{0}) \}$$
Let $\eta(\psi)=h(\psi)/\psi$.
\newtheorem{e_sgn}{ Lemma }[section]
\begin{e_sgn} \label{e_sgn}
Assume $h$ such that $h^{2}=1$, then $\eta$ does not change sign over $Eq$.
\end{e_sgn}

The proof is trivial. Assume $Eq$ is a discrete set and that the number of elements is $2$. We then have
$hP_{c}(\psi)=(c h(\psi)-\psi)^{2} \geq 0$ (with $h^{*}=h$ and $\psi$ real). Thus $P_{c}(\psi) \geq0$ under
these conditions. In the same manner, if the number of elements is odd ($=\mu-1$) and $\zeta_{0}$ such that
$\eta(\psi) \geq c$, then we have $P_{c}(\psi) \leq 0$. The sign for $P_{c}(\psi)$ is dependent on $\zeta_{0}$
and not on $\psi \in \mbox{ nbhd }(I_{c})$. Note that if $\eta(\psi)(\zeta)=c + P(\zeta)$ (in a domain of holomorphy)
where $P$ is a first degree polynomial. In this case $P_{c}(\psi)=P(\zeta_{0})^{\mu - 1}=(a_{0}+a_{1}\zeta_{0})^{\mu - 1}$.

\vspace*{.5cm}

Consider now $\eta(\psi)(\zeta)=\mu + P_{1}(\zeta)$, with $P_{1}$ polynomial, then
$h(\phi)=b_{0}+b_{1}\phi + b_{2}\phi^{2}$. We can assume $b_{0}=0$ ($h$ algebraic) and apply $h$ again to get
$\eta(\psi)(\zeta)=(1-c_{2}(\zeta))/c_{1}(\zeta)$. We could consider $\{ \phi,\phi^{2} \}$ as
a pseudo-base for an ideal $(I)$ and as before consider
the distance to $V=N(h(I))$, $d(\zeta,V) \leq C \mid \zeta \mid^{c \rho} \big( c_{1}\mid \phi
\mid^{\frac{1}{\kappa + \sigma}} + c_{2})^{\rho}$ with $\rho < 1$, where
$\kappa + \sigma$ is a positive constant dependent on choice of neighborhood of the origin. Let $\Omega_{\mu}=\{ \zeta \quad
h(\psi)(\zeta)=\mu \psi(\zeta) \}$ for a constant $\mu$. Assume $\psi=e^{\varphi}$ such that
$\varphi \rightarrow - \infty$ as $\zeta \rightarrow 0$, over $H_{h}(\psi)=0$ we further have
$d_{\zeta}h(\psi)=\mu d_{\zeta} \psi$. For $\zeta \notin \Omega_{\mu}$ we have $h(\psi)=\sum_{j} c_{j} \psi^{j/\mu}$ and the image under
$t$ lies outside the tangent-line $\Gamma$. The corresponding trajectory determines a spiral
approximation of $\Omega_{\mu}$. Starting with the Puiseux-expansion $h(\psi)=c_{0}+c_{1}\psi^{\frac{1}{\mu}}
+ c_{2}\psi^{\frac{2}{\mu}}$ and assuming $h \sim_{m} g$ where $g$ is an analytic homomorphism, we
have $\tilde{t} h(\psi) \sim \sum_{j} c_{j}z_{1}^{\frac{j}{\mu}}$. Monotropy for the trajectory $\gamma$
corresponds to translation of $\zeta$ that is $\gamma_{g}=\tau_{\epsilon} \gamma_{h}$, why
$\tilde{t}g(\psi)=\sum \big[ \tau_{\epsilon}c_{j}(\zeta) \big] z_{1}^{\frac{j}{\mu}}$ where
$\tilde{t}$ is the projection on the second variable of $\gamma$.

\vspace*{.5cm}

Remark:\\
The condition $dw$ locally reduced for a homomorphism $w$ deserves explanation. If $\psi$ is
singular for the dynamical system, it is through our conditions isolated. If $\{ \psi_{j} \}$ is a
sequence such that $d \psi_{j} \neq 0$ and $dw(\psi - \psi_{j}) \rightarrow 0$, then $\{ \psi_{j}
\}$ must be a regular approximation of $\psi$, that is $\psi - \psi_{j} \rightarrow 0$. In this
sense $d w$ is reduced over regular approximations.

\section{ Hanges' theorem}
We can now generalize Hanges' theorem to pseudo-differential operators $A_{\lambda}$ with foliation
that is included in the zero-set to the polynomial corresponding (using $t$) to $X=\{\gamma \quad P(\gamma)(\zeta)=0 \}$.
Assume $\gamma$ an integral curve with a real branch in $X$ and that the associated dynamical system has
polynomial right hand sides, we then have $\tilde{\gamma}(U) \subset WF_{a}(u)$, for solutions to the
homogeneous equation to $A_{\lambda}$ and for $U$ a neighborhood of the origin. Conversely, if $\gamma$
is not a tangent to the origin and does not have a real branch in $X$, we know that $\tilde{\gamma}(U)
\cap WF_{a}(u)=\emptyset$

\vspace*{.5cm}

Assume for a symbol $a$, holomorphic and of finite type and corresponding a self-adjoint operator,
the intersection of the foliation and the transversal set is included in the zero-set to a
holomorphic function $F$. Assume further that $F$ is in a finitely generated ideal $(I)=\mbox{ ker
}h$, where $h$ is assumed ( modulo monotropy ) algebraic and such that $h^{2}$ is locally injective. We will here  consider the case where $h(F)$ real. Further
that $A \sim_{\infty} \tilde{F}+\Psi$ where $\Psi$ is assumed regularizing and $\tilde{F}$ is the operator
corresponding to $F$. Hanges (\cite{Han}) studies the characteristic set corresponding to the
principal part of the symbol $a$. We propose that if $F$ is defined by $h(F)=0$, the
homogeneous condition can be replaced by the following.

\vspace*{.5cm}

Assume $t$ a real variable and consider $\zeta \rightarrow t \zeta$ such that $h(t^{2} F) \sim_{m}
t h(F)$. Define a continuous function $s$ according to $s(t \zeta)=h(t F)(\zeta)$, then $s$ is
homogeneous of order $\frac{1}{2}$. Let further $r^{2}(t \zeta)=t r^{2}(\zeta)$, where $r$ is analogous to $s$.
A conical neighborhood
of $z_{0}$ is mapped by a canonical homogeneous transformation on an
open conical set in $(x,t,\xi,\tau) \in T^{*}(\mathbf{R}^{\textit{n}+\textsl{1}}) \backslash \textit{0}$. We can
now argue as
Hanges (\cite{Han}) that is the operator corresponding to $r^{2}(\zeta)$ is reduced to canonical form $t \frac{\delta}{\delta t} -
\tilde{B}(x,D_{x})$ close to $(x_{0},0,\xi_{0},0)$ and $\tilde{B} \in L^{0}_{c}$ with proper
support. We study a mapping
$T \text{:} (I) \rightarrow \mbox{ nbhd }\Gamma$, where $\Gamma$ is the tangent line to the
foliation and where the associated dynamical system is assumed to satisfy Bendixson's regularity
conditions (cf. \cite{Ben}). We consider $\Sigma=\{ \zeta \quad Tr(\zeta)=0 \}$ and $\tilde{\Sigma}=\{ \psi \quad
h(\psi)=0 \}$. According to Bendixson (cf. \cite{Ben} Ch.2, Theorem 9 ) if the dynamical system
has polynomial part right hand sides of order $1$, we have exactly four trajectories that run
through the origin. Hanges considers correspondingly four half-bicharacteristics defined through
$\gamma_{j} \text{:}  J \rightarrow \Sigma_{j}$,  $j=1,2$ $\gamma_{j}(0)=z_{0}$, where $J$ is an
open interval in $\mathbf{R}$ containing $0$. For $J_{k}=\{ t \in J \quad (-1)^{k}t > 0 \}$ $k=1,2$
we write $\gamma_{j,k}=\gamma_{j} \mid J_{k}$ for these. Note that in our case $\Sigma_{j}=\Sigma$.
The following proposition is a direct consequence of (\cite{Han} Prop.1)

\newtheorem{e_Hanges}[e_sgn]{ Proposition}
\begin{e_Hanges}

Assume $\gamma_{j,k}$ a characteristic such that $\gamma_{j,k}(0)=z_{0}$. If for a pseudo-differential operator
$A$ described above, $u \in \mathcal{D}'$,
solves the homogeneous equation $Au=0$ and $\gamma_{j,k} \cap WF_{a}(u) = \emptyset$, then $z_{0}
\notin WF_{a}(u)$.
\end{e_Hanges}

\subsection{ The non-homogeneous equation }
Assume $P$ a polynomial operator and $P(D)u=f$ in $\mathcal{D}'$ for $u \in \mathcal{D}'$ and $\gamma$ is a
trajectory to the dynamical system such that $\tilde{\gamma} \cap WF_{a}(u)=\emptyset$

\newtheorem{e_reg_appr}[e_sgn]{ Lemma }
\begin{e_reg_appr}
If $\gamma' \rightarrow P_{0}$ singular through regular points, we have that $\tilde{\gamma'} \cap
WF_{a}(u)=\emptyset$
\end{e_reg_appr}

Proof:\\
Since $P_{0}$ is mapped by $t$ on the foliation to an operator of finite type, we have
$\tilde{\gamma'} \rightarrow$ a finite value. The indicator of growth to an operator $F$ is then
the same as the indicator to $\tau_{\tilde{\gamma'}}F$ as $\gamma' \rightarrow P_{0}$ ($\varphi \rightarrow \pm
\infty$). Note that we have not assumed that $\gamma'$ is a trajectory to the dynamical system.$\Box$

\vspace*{.5cm}

Note that if $P^{2}$ is hypoelliptic and self-adjoint, we have $WF(Pu) \subset WF(u)$ and if $u=Pv$ in
$\mathcal{D}'$ we have $WF(P^{2}v)=WF(v) \subset WF(Pv)$. The geometrical objects are considered as the
same, it is only the names that are swapped. We assume that the latter case is used. We are going to prove
that $WF(Pu)=\{ \tilde{\gamma} \} \cup WF(u)$, where $\tilde{\gamma}$ is as in the lemma and $P$ is
a pseudo-differential operator with foliation in the zero-set to a polynomial $F$.

\newtheorem{e_exist}[e_sgn]{ Lemma }
\begin{e_exist}

Assume $\gamma$ a trajectory to the dynamical system and $\gamma \rightarrow P_{0}$ singular, with
a tangent determined in $P_{0}$, then there is a regular approximation of $P_{0}$, for $P$ as above.
\end{e_exist}

Proof:\\
The conditions give that all approximations $\gamma$ of $P_{0}$ have a tangent determined
$(I_{c})$, $c$ constant, in $P_{0}$. As the foliation to the operator is finitely generated
(in an algebraic variety), we have the minimally defined case and there always is a regular
approximation.$\Box$

Remark:\\
Note that the conditions in the second lemma are satisfied if the dynamical system has polynomial
right-hand sides.

\newtheorem{e_transappr}[e_sgn]{ Lemma }
\begin{e_transappr}
Assume $P_{0}$ singular and mapped by $t$ onto $N(I)$ for a finitely generated ideal $(I)$. If
$\gamma$ is a regular approximation of $P_{0}$ with tangent determined $(I_{c})$ $c$ constant,
then $\gamma$ is in the set of transversals.
\end{e_transappr}

\vspace*{.5cm}

Assume now $A$ a pseudo-differential operator of finite exponential type and with foliation in
$N(F)$  for $F \in (I)$ polynomial and $(I)$ a finitely generated ideal (the minimally defined
situation). Assume $\tilde{\gamma} \cap WF_{a}(Au)=\emptyset$ and $\tilde{\gamma}$ transversal, that is
$\gamma$ is a regular approximation of $P_{0}$, a singular point. If we have the representation
$Au=\tilde{F}u+Ru$ where $R$ is regularizing and $\tilde{F}$ has zero's that include the foliation
to $A$, then the generalized Hanges' result gives that $WF_{a}((A-\tilde{F})u) \cap
\tilde{\gamma}=\emptyset$ and $WF_{a}(u) \cap \tilde{\gamma}=\emptyset$.

\newtheorem{e_HangesII}[e_sgn]{ Proposition }
\label{e_HangesII}
\begin{e_HangesII}
Assume $A$ a pseudo-differential operator, self-adjoint and of finite exponential type. Further
that the foliation is included in the zero-set to a polynomial or a difference of two polynomials
in a finitely generated ideal. Assume every foliation point mapped on a singular point for the
associated dynamical system, such that every trajectory that reaches such a singular point, has
a tangent determined in this point. Then every regular approximation of the singular point is
mapped in the set of transversals $\tilde{S}$, so that $WF_{a}(Au)=\tilde{S}\cup WF_{a}(u)$.
\end{e_HangesII}

\section{ Foliation in a semi-algebraic set}
Assume an ideal of holomorphy $(I)$ defined by a polynomial $P_{\mu}$ and consider
$(\overline{I})$ the Whitney-closure. Define the subspace $(I_{\lambda})=\{ \gamma \quad P_{\mu}(\gamma) \leq \lambda
\}$ and the corresponding ideal $(\tilde{I}_{\lambda})$ under the mapping $t$. We will now discuss some properties of the class of operators that has foliation in intersection with transversal set included
in the semi-algebraic set $\{ \zeta \quad \tilde{P_{\mu}} \leq \lambda \}$.

\vspace*{.5cm}

If the symbol to our pseudo-differential operator has foliation in $\{ \zeta \quad \tilde{P}_{\mu}
\leq \lambda \}$ we can prove that the corresponding ideal $(I)$ containing the polynomial
$P_{\mu}$, has a global pseudo-base. Assume $(I)$ finitely generated such that the minimally
defined zero-set case holds, that is assume over the semi-algebraic set that $P_{\mu}=0$ $\implies
d P_{\mu} \neq 0$. Singular points are now such that $\{ P_{\mu}=\lambda \}$, $\lambda$ constant and also the semi-algebraic
set is minimally defined, since the singular points are the foliation to a polynomial $P_{\mu}$,
which is assumed to be with lineality (not reduced) and we have regular approximations over
transversals. Note that for $\psi=e^{\varphi}$, as $(I_{c})=(I_{c})'=(I_{c})''$ we have
$P(\delta_{j} \gamma)=P(d \gamma)$ for all $j$, why in the minimally defined case, we must have for all $j$,
$P_{\mu} \delta_{j}-\delta_{j} P_{\mu} \neq 0$.

\vspace*{.5cm}

Given exactness for $\eta$ ($=h(\psi)/\psi$) we can show existence of a global pseudo-base also for the ideal that
defines the foliation $(I)$. Assume $P_{\mu} \in (I)$ and that the semi-algebraic set $\{ P_{\mu}
\leq \lambda \}$ is minimally defined such that we always have regular approximation of singular
points $(P_{\mu}(0)=\lambda)$ and under these conditions every regular approximation has a tangent
determined.

\newtheorem{e_HangesIII}[e_sgn]{ Proposition }
\begin{e_HangesIII}
Consider for a pseudo-differential operator $A$, $Au=f$ in $\mathcal{D}'$, for $u \in \mathcal{D}'$, where we
assume $A$ has symbol with a foliation in a semi-algebraic set $\{ P_{\mu} \leq \lambda \}$ that is
minimally defined and where $P_{\mu} \in (I)$ is a polynomial and $(I)$ is assumed to have a global
pseudo-base. Then for $\tilde{S}$ the transversal set, $WF_{a}(Au)=\tilde{S} \cup WF_{a}(u)$.
\end{e_HangesIII}

\vspace*{.5cm}

Assume now $\gamma \in L_{\mu}$ and the foliation in a semi-algebraic set as before, this means
that $\{ \mu(d \gamma) \leq \lambda \}$ that is a bounded positively definite measure defined close to the
foliation. Particularly if the ideal $(I)$ is finitely generated and the index $\rho \leq 1/N$ for
an integer $N$, we have for $M_{\mu}=\{ \gamma \quad \mu(d \gamma) \leq \lambda \}$ that $\mid
\mu(d \gamma) \mid \leq C \mid d \psi \mid^{\rho}$ for a positive constant $C$ and conversely if $\mid d \psi
\mid \leq \lambda^{N}$ we must have $\gamma \in M_{\mu}$.

\vspace*{.5cm}

Assume existence of a global pseudo-base for $(I)$ and $\Phi_{\mu}$ the generalized
Cousin-integral (cf. \cite{jag_0}), then there are $\widetilde{\Phi}_{\mu}$ of type $0$. Assume $\widetilde{\mu}$ such
that $\widetilde{\Phi}_{\mu}=\Phi_{\widetilde{\mu}}$ ( addition of measures with point-support ),
then the corresponding polynomial can be chosen reduced, that is with no foliation why the only
singular point possible is the origin. The conclusion is that for every trajectory to the
associated dynamical system that is tangent to $M_{\widetilde{\mu}}$ holds that this trajectory
will stay in a bounded set in phase space as $\varphi \rightarrow \pm \infty$. For $M_{\mu}$ the
same conclusion holds $\mu-$ a.e.

\section{ A case of infinitely generated foliation }
We have earlier established that in the case where $h^{2}$ is locally injective and $\eta$ ($=h(\psi)/\psi$) exact,
that we have a global pseudo-base for the ideal that generates the integral curves. Note that since
$h$ is reduced over the tangent space, if we form $(I)=\{ \phi \quad \delta_{j} \phi \in (J) \quad j=1,\ldots,n \}$ we
have a local pseudo-base for $(I)$ and a global pseudobase for $(J)$. Assume now
$\eta_{j}(\psi)=h^{j+1}(\psi)/h^{j}(\psi)$ without improvement (or degeneration) of behavior by iteration,
that is we are assuming a not-finitely generated ideal. We give the following proposition, where $\tau$ denotes translation

\newtheorem{e_bounded_psi}[e_sgn]{ Proposition }
\begin{e_bounded_psi} \label{e_bounded_psi}
Assume with $S_{\tau} \psi=\tau \psi/\psi$, $h(S_{\tau} \psi) \sim_{m} 0$ in the $\mid \zeta \mid -$infinity and $\eta(\psi)$ not constant, then $\psi$ is in a bounded
set, symmetric with respect to the origin.
\end{e_bounded_psi}

Proof:\\
If $\eta_{j}$ is constant then $\psi$ must be on a transversal. Consider
$\psi_{k}$ in a regular approximation of $\psi$ such that $\psi_{k}$ are non-constant as well as
$\eta_{j}(\psi_{k})$. Assume further
\begin{equation} \label{e_slow_osc}
h(S_{\tau_{z}} \psi) \sim_{m} 0 \quad \mid \zeta \mid \rightarrow \infty
\end{equation}
particularly we have seen before that if $h$ is algebraic, $S_{\tau_{z}} \psi \sim_{m} 0$ as $\mid
\zeta \mid \rightarrow \infty$.

\vspace*{.5cm}

Obviously we have $\eta_{j}(\psi) \sim_{m} h(\psi)$, why if $\mid \eta_{j}(\psi)-h(\psi) \mid <
1/\mid \zeta \mid$ in the $\mid \zeta \mid$-infinity, we have $\mid \eta_{j}(\psi) \mid < 1/\mid
\zeta \mid + \mid h(\psi) \mid$ and using the condition (\ref{e_slow_osc}) we have that $\eta(\tau_{z} \psi /
\psi) \sim_{m} 0$ as $\mid \zeta \mid \rightarrow \infty$ (modulo monotropy).
Consider now $\psi(z)=\tau_{z}\phi(x)/\phi(x)$ non-constant, why if
$\phi$ is reduced, then $\psi$ is bounded $\forall x,z$.
According to
(\ref{e_slow_osc}) we have that $\psi \eta(\psi) \sim_{m} 0$ why we can claim that $\eta$ is bounded in
the $\mid \zeta \mid$-infinity.  Let $(I)$ be the $\phi$ for which $\psi$
is bounded when $\tau$ is compact. For reduced $\phi$ we know that the corresponding $\phi$ have a
global pseudo-base. Obviously, the two ideals are locally equivalent in the $\mid \zeta
\mid$-infinity.

\vspace*{.5cm}

Let $h(\phi^{2})-h^{2}(\phi) \sim_{m} 0$ such that $h(\phi^{2}-h(\phi)) \sim_{m} 0$. Assume further
that $\phi \neq 0$ in the $\mid \zeta \mid$-infinity such that $\phi - \eta(\phi) \sim_{m} 0$. We
then have that modulo monotropy $\mid \phi \mid < 1/\mid \zeta \mid + \mid \eta(\phi) \mid$. In
this sense $\eta$ is downward bounded in the infinity. Further symmetry with respect to the origin
in $\zeta$ can be shown for regular approximations of $(I_{c})$, $c$ a constant.$\Box$

\newtheorem{e_global_p_b}[e_sgn]{ Proposition}
\begin{e_global_p_b}
Let $(I)$ be the ideal of $\phi$ such that $\tau_{z} \phi/\phi \sim_{m} 0$ in the infinity and $\tau_{z}$ compact. Then
$(I)$ has a global pseudo-base.
\end{e_global_p_b}

Proof:\\
Assume Schwartz-type topology and a compact translation. We can then assume
the ideal $(I)$ finitely generated. Further it is not difficult to prove, for instance using the quantitative
version of R\"uckert's Nullstellensatz, that $\eta \sim_{m} 0$ over $(I)$. Thus the indicator for $\eta$ over
$(I)$ is $0$ and we can find a transform $\gamma \sim_{m} \eta$ over $(I)$, that is exact so that the ideal
defined by $\gamma$ has a global pseudo-base.$\Box$

\vspace*{.5cm}

 Assume the conditions in Proposition \ref{e_bounded_psi} further $d \eta \neq 0$ and continuous. Then
 $\psi \in B$ where $B$ is a bounded set, symmetric with respect to the origin. Considered as a
 cone, $tB=\Gamma$ always has negative indicator, that is for $\xi_{0}$ finite and $s$ a real number,
 we have $< s y,\xi_{0} > < 1$ $\Leftrightarrow <y, \xi_{0} > <1/s\rightarrow 0$ as $s \rightarrow \infty$.
 We now know (cf. \cite{Sjo} Theorem 6.6 ) that in this case we have no contribution to the
 wave-front-set.

 \newtheorem{e_main}[e_sgn]{ Proposition }
 \begin{e_main}

 Assume $\psi$ such that $h(S_{\tau} \psi) \sim_{m} 0$ in the infinity and $\eta$ such that $d\eta(\psi) \neq 0$
 and continuous, we then have that $\{ t \psi \}$ is not mapped on to the wave-front-set.
 \end{e_main}

\section{ Implications on microlocal analysis }
Assume $\zeta_0$ a point in the plane, such that $h(f)(\zeta_0)=0$ of order $N$ and consider a neighborhood $W$
of $\zeta_0$, such that $\zeta_0$ is an isolated zero. Then $(\zeta - \zeta_0)^N v(\zeta)=h(f)(\zeta)$,
where $v$ is real analytic and such that $v(\zeta_0) \neq 0$. Assume
$\gamma$ a real-analytic function such that $\gamma(x_0)=\zeta_0$ (implicit function theorem).
Then, in $x_0$, $f(\gamma(x_0))=$
constant. Let $$X_c=\{ x \quad f(\gamma(x))=c \} \text{ for a constant c }.$$ If $f=\widehat{u}$, for $u \in
\mathcal{S}$ and if $u$ has support on $X_c$, then $X_c \subset \mbox{ sing supp }u$. Further, if we
consider full lines in $V_c$, we can chose $\gamma$ of order $1$ and such that
$h(f)(\gamma(x))=c_N(x-x_0)^N v(\gamma(x))$. More precisely, in the Cousin integral approach, if
$f=\widetilde{\Phi}-\Phi$ and $e^{\varphi}=f$, we can give $V_{c}$ as $\{ \zeta \ \varphi(\zeta)=const. \}$.
Consider $V= \mbox{ nbhd }\zeta_0$ and $f(\zeta_0)=c$.
Let $\eta(V_{c} \cap V ) \subset X_{c} \cap U$, for
$U=\mbox{ nbhd }x_0$, where $\gamma,\eta$ are real-analytical functions ( of order $1$ ) and
injective on these sets. Let $\zeta \in V_{c}$ and $u(\eta(\zeta))=\widehat{v}(\zeta)$, for
$\widehat{v} \in \mathcal{E}'$ $(V)$. Assume $\widehat{v}$ is constant on $\Gamma \cap V$, for a line $\Gamma$.
Then, $\mathcal{F}^{-\textsl{1}}$ $\big( u(\eta( \zeta)) \big)=$ $c \delta_{\Gamma}(x)$, for $\Gamma=\{ \gamma(x)
\}$ with $\delta_{\Gamma}(x)=\delta_{0}(x-\Gamma)$. Note that if $u \in \mathcal{S}'$, then
$u(\eta(\zeta)) \in \mathcal{S}'(\textsl{V})$.

\vspace*{.5cm}

The restriction (inclusion)-homomorphism $r_{\Gamma}$ considered on $\mathbf{R}^{\textsl{n}}$
$C^{\infty}_r \rightarrow C^{\infty}_r \rightarrow C^{\infty}$, can be extended to $$\mathcal{E'}
\rightarrow (C^{\infty}_{\textit{r}})'(\textit{U}) \rightarrow (C^{\infty}_{\textit{r}})'(\Gamma),$$ where $C^{\infty}_r$ denotes
real-analytic functions. Assume $\beta$ another homomorphism on $\mathcal{E'}$, such that $\beta(T)(x)=0$ for
$T \in \mathcal{E'}$ and $x \in \Gamma$, a line-segment, implies $T$ is constant on this line-segment, that is
$r'T=c \delta_{\Gamma}$. Let now $\beta$ be defined with its support outside $\mbox{ sing supp }T$.
Consider for instance the class of distributions $\mathcal{E'}_{\Omega}(\textsl{V})$ for open sets $\Omega
\subset V$ as the class $\{ T \in \mathcal{E'}(\textsl{V}) \quad \mbox{ sing supp }\textsl{T} \subset \Omega \}$ and
assume $\beta$ with support on $V \backslash \Omega$. Let $$W_h=\{ (x, \gamma(x)) \in U \times V \quad \beta(T)(x)=0 \}$$
and $$W_c=\{ (\eta(\zeta),\zeta) \in U \times V \quad \zeta \in V_c \}.$$ Thus, over $\mathcal{E'}_{\Omega}(\textsl{V})$ we have that
$\beta$ defines the foliation and $W_{c}=W_{h}$

\vspace*{.5cm}

Assume $r'_L$ the restriction-homomorphism to a convex, compact set $L$ and acting on $\mathcal{E'}$.
If $\Gamma$ is the polar to $L$, a line-segment through $x$, we have that also $\Gamma$ is a line-segment and the
functional $\mathcal{F}^{-\textsl{1}}$ $r_{\Gamma}'\mathcal{F}$ $T$ is portable by $L$. It also follows that
$r'_L$ is portable by $L$. We let $\mathcal{F}$$r'_L=$ $r'_{\Gamma}\mathcal{F}$ in $\mbox{ Exp
}_{\rho,A}$. Assume $\gamma,\eta$ real-analytic functions mapping $L$ to $\Gamma$ and back (inverse
function theorem), we have the following proposition,

\newtheorem{c_Delta}[c_RED]{ Proposition }
\begin{c_Delta}
If $u \in \mathcal{E'}_{\Omega}(\textsl{V})$, $x \in \mbox{sing supp }u$ and $x$ on a line-segment $L$
in $\Omega$. Then there is a constant $C$, such that $r'_L u=C\delta_{L}$. If $f=\widehat{u}$ and
$\zeta \in V_{c}$ and on a line-segment $\Gamma$ in $V_{c}$, then $r'_L u=C\delta_{L}$, for a constant $C$.
\end{c_Delta}

Assume $\beta\text{:} C^{\infty} \rightarrow C^{\infty}$ and
$\alpha= \beta \mathcal{F} - \mathcal{F}$ $\beta \text{:} \mathcal{E'} \rightarrow C^{\infty}$ and for the restriction to
a line-segment, $\mathcal{E'} \rightarrow C^{\infty}_{\textsl{0}}$. As $\beta u(x)-\beta u(-x)=\alpha \mathcal{F} \textit{u} + \mathcal{F}
\alpha \textit{u}$, we have that $\beta u(x)-\beta u(-x) \in C^{\infty}(\text{nbhd}L)$. Assume
$\beta u(x)-\beta u(-x) \neq 0$ on $L$ and $\zeta \in U_c \backslash V_c$, $x=\eta(\zeta) \in L$.
Then $$\frac{1}{\mid x \mid^c} \leq C \mid \beta u(x)-\beta u(-x) \mid \leq C'( \mid \beta u(x) \mid +
\mid \beta u(-x) \mid ) \text{ on } L,$$ thus $x \notin \mbox{ sing supp }u$ or $-x \notin \mbox{ sing supp
}u$. Assume $M$ a set, symmetric around $0$, $M=\mbox{ supp }F$ for a real and real-analytic
function $F$. Define $$Eq^{\pm}=\{ \zeta \in M \ F(\zeta)=F(-\zeta) \}.$$ Then if $F_1=\Pi_{\xi \in Eq}(\zeta -
\xi(\zeta))F(\zeta)$ and $\mbox{ supp }F_1 \subset M \backslash Eq^{\pm}$. In the same manner, if
$Eq'=\{ \zeta \in M \ F_2(\zeta)=0 \}$, for $F_2$ real and real-analytic, we can define $F_1'$ with
$\mbox{ supp }F_1' \subset M \backslash Eq'$. The conditions above are satisfied for $L \subset
\mbox{ supp }F_1'$.

\vspace*{.5cm}

 Let $A=\{ \zeta \quad h(\widehat{u})(\zeta)=dh(\widehat{u})(\zeta)=0 \}$.
Note that $\overline{J}_{h}$ as a closed ideal, can be considered as radical. Assume $h(f)=\alpha
f$, for $\alpha \in H$ and $h^2$ injective. Then, $\mid f \mid^2 \in \overline{J}_{h}$ is
equivalent with the proposition $h(\mid f \mid^2)=constant$. Further, if $h(\mid f \mid^2)=\beta
\mid f \mid$, for a constant $\beta$, we conclude that $\mid f \mid=constant$.

\newtheorem{c_mindef}[c_RED]{ Proposition }
\begin{c_mindef}
 Assume $N=\{ \zeta \quad h(f)(\zeta)=0\}$, $C_{\Gamma}=\{ \zeta \quad dh(f)(\zeta)=0 \}$, $C=\cap_{\Gamma} C_{\Gamma}$ and $A=C \cap N$,
then for $\Gamma_{1}$ a line through $\zeta$, we have
\begin{displaymath}
      \left \{
\begin{array}{lr}
\mid r_{\Gamma_{1}} f \mid=constant \text{ for } \zeta \in C \\
\quad r_{\Gamma_{1}}f=constant \text{ for } \zeta \in A
\end{array} \right.
\end{displaymath}

\end{c_mindef}

\section{ Symmetry and the singular support }
Let $\Omega_{A}=\{ \zeta \quad \mid \zeta \mid < A \}$ and for a small set $U$, $\Gamma_{A}^{U}=\{ \zeta \in U \quad \mid \zeta \mid=A \}$.
Let $\Omega_{\epsilon}=\{ x \quad \frac{1}{2}\epsilon < \mid x \mid < \epsilon \}$ and
$x=\eta(\zeta)=\sum_j a_jt^{j}$ the Puiseux-expansion corresponding to the singular points on $\Omega_{A}$. Then for $\zeta \in \Omega_A$ and $\eta(\zeta)-\xi \in
\Omega_{\epsilon}$, $\mid \eta(\zeta)\mid < \mid \xi \mid +\epsilon$. Assume as before,
$f(\zeta)=\widehat{c}(\zeta)=u(\eta(\zeta))$ and further $\zeta \in \Omega_A$ implies $\eta(\zeta)
\notin \mbox{ sing supp }u$. In the construction of the Green-function, we use
$T_{\epsilon}=\alpha_{\epsilon}g_{\lambda}$ for $\alpha_{\epsilon} \in \mathcal{D}$ with support in
$\mid x \mid \leq \epsilon$ and with $\alpha_{\epsilon}(x)=1$ for $\mid x \mid \leq \frac{1}{2}
\epsilon$, it is a function only of $\mid x \mid$. $g_{\lambda}$ is a tempered fundamental solution
to $D_{\lambda}=\Delta - \lambda$. We then have $T_{\epsilon}(\eta(\zeta)-\xi) \in
C^{\infty}(\Omega_{\epsilon})$. Let $N_{\lambda}=N(D_{\lambda})$. Then on $U_{c} \backslash
N_{\lambda}$, $1 \leq \mid \widehat{D_{\lambda}} \mid \mid \widehat{g}_{\lambda} \mid \leq C \mid
\xi \mid^{q-c} \mid \widehat{g}_{\lambda} \mid$, why if $c > q$, we can regard
$\widehat{g}_{\lambda}$ as reduced. If further $f$ is reduced, for $\zeta \in U_{c} \backslash
N_{\lambda}$, $\mid \eta(\zeta) - \xi \mid^{\sigma} \leq \mid \check{T}_{\epsilon}*c(\eta(\zeta))
\mid$. If $\xi=\eta(v)$ for $v \in \mbox{nbhd} \zeta \subset \Omega_A$, $v=\zeta +
\epsilon e^{i \phi}$ we have as $\eta$ is real-analytical, for $\mid \zeta \mid < r <A$, that
$\mid \eta(\zeta)-\xi \mid \leq \sum_{\alpha} C_{\alpha}'r^{-\mid \alpha \mid}\epsilon^{\mid \alpha
\mid}$

\vspace*{.5cm}

We can regard $T_{\epsilon} \in \mathcal{E}'$ as very regular in $\mathcal{D'}^{F}$. Thus, $T_{\epsilon}$
is on the form $\delta_{0} - \gamma_{\epsilon}$, with $\gamma_{\epsilon}$ regularizing and
$\gamma_{\epsilon}(0)=0$. Assume $\Gamma \mbox{:} \{ \xi=\eta(\zeta) \}$. We have
$$ \int_{\Gamma} T_{\epsilon}(\eta(\zeta)-\xi)u(\xi)d \xi = \int_{\Gamma} \delta_{0}(\eta(\zeta) -
\xi)u(\xi)d \xi= \sum_{\zeta_j} u(\eta(\zeta_j)) = \sum_{j}^{\text{ finite sum}} f(\zeta_j)$$
This is also the expression we expect, when we let $\epsilon \rightarrow 0$ in
$T_{\epsilon}$.

\vspace*{.5cm}

Assume $V=\mbox{nbhd} \zeta$ and $u(\eta(\zeta)) \in C^{\infty}(V)$. Let $U=\eta(V)$ such that
$u(x) \in C^{\infty}(U)$ (also $\gamma(U)=V$). Let
$$ W'_{L^2}(u)=\{ \zeta \quad u(\eta(\zeta)) \in C^{\infty}(V') \quad \exists V'=\mbox{nbhd} \zeta
\quad u \in L^2 \}$$
Assume $u$ has the property that

\begin{equation} \eta(\zeta) \notin \mbox{ sing supp }_{L^2} \mbox{ Re }u \Rightarrow \eta(\zeta) \notin
\mbox{ sing supp }_{L^2}u \label{c_cond}
\end{equation}

Then, obviously for $u$ according to (\ref{c_cond}), we have $\zeta \in W'(\mbox{ Re }u)$ implies
$\zeta \in W'(\overline{u})$. We regard, for a polynomial $P(\zeta)$, $P\chi_{\Omega_A}$ as a
bounded operator in $L^2$. So, if for a constant $\beta$, $\beta \overline{P}=u(\eta)$ in $L^2$, $\mbox{ Re
}P=\frac{1}{2}(P+P^{*})$ and $\beta(P+P^{*})=2\mbox{ Re }u=u(\eta(\overline{\zeta})) +
\overline{u}(\eta(\zeta))$

\newtheorem{c_Sym}{ Lemma }[section]
\begin{c_Sym}
For $u$ according to (\ref{c_cond}), we have $\zeta \in W'_{L^2}(\mbox{ Re }u)$ implies
$\overline{\zeta} \in W'_{L^2}(u)$
\end{c_Sym}

In the general case, we have that the size of $V \subset W_{L^2}'$ is directly proportional to the
distance from $\zeta$ to $0$. Conversely to (\ref{c_cond}),

\newtheorem{c_Symb}[c_Sym]{ Proposition }
\begin{c_Symb}
For $V \subset W'_{L^2}(u)$ and $V \subset W'_{L^2}(\overline{u})$, where $V$ is symmetric with
respect to the real axis, then $V \subset W'_{L^2}(\mbox{ Re }u)$
\end{c_Symb}

\section{ On preservation of flatness }

Let for a homomorphism $\beta$, $$I_C=\{ g \in \overline{J}_h \quad \beta(f-g)=0 \quad \exists f \in J^{\infty}_C \}$$
such that any $g \in I_C$ can be represented as $g=f+c$, for a constant $c$ on $\Gamma^{U}_{A}$, where
$f_{\alpha}(\eta)=0$ for all $\alpha$ and $\eta \in C$. Let $C$ be a cone with vertex in the origin
and $f \in J^{\infty}_{h}(T^{C})$, where $T^{C}=\mathbf{R}^{\textit{n}} + \textit{i} \textsl{C}$ and
$J^{\infty}_{h}(f)=\{ D^{\alpha}I_{H}(f) \}_{\alpha}$, then $I_{H}(f)(x+ i t \eta)-I_{H}(f)(x)=0$,
for $t \eta \in C$, why $C \subset \Delta_{\mathbf{C}}(I_{H}(f))$ and $I_{H}(f)$ is flat on $T^{C}$.

\vspace*{.5cm}

Assume $W=\{ \zeta \quad D^{\alpha}I_{H}(f)(\zeta)=0 \quad \forall \alpha\}$. If
$D^{\alpha}I_{H}(f)=I_{H}(D^{\alpha}f)$ for all $\alpha$, we have that $h$ preserves flatness.
Now consider
$$ I_{B_{\alpha}}(f)=D^{\alpha}_{x}I_{H}(f) - I_{H}(D^{\alpha}_{y}f)$$
(also consider $I_{B_{\alpha}'(f)=I_{(D^{\alpha}_{x}H-D^{\alpha}_{y}H)}(f) }$, if
${}^tI_{H}=I_{H}^{*}$, we have $I_{B_{\alpha}'}=I_{B_{\alpha}}$). Immediately,
$I_{B_{\alpha}}=\big[ D^{\alpha}_{x}I_{H} - \overline{D^{\alpha}_{x}}I_{H}^{*} \big] +
\big[ \overline{D^{\alpha}_{x}}I_{H}^{*} - \big( \overline{D^{\alpha}_{x}}I_{H}^{*} \big)^{*} \big]$.
Let $J(I_{H})=I_{H} - I_{H}^{*}$. Then
$$ I_{B_{\alpha}}=D^{\alpha}_{x}J(I_{H}) + J(D^{\alpha}_{x})I_{H}^{*} +
J(\overline{D^{\alpha}_{x}}I_{H}^{*})$$
If $I_{H}$ is a bounded operator in $L^2$, we have $J(I_{H})=2 i \mbox{ Im }I_{H}$. If
${}^tI_{H}=I_{H}^{*}$ we have $J(I_{H})=0$ and if $I_{H}$ is an integral operator with distribution
kernel, then $B_{\alpha}$ can be neglected. With the condition that $h^2$ is locally injective, it
is clear that $h^2$ preserves flatness and if $D^{\alpha}I_{H}(f) \in J_{h}$, $\forall \alpha$, we have that
$cD^{\alpha}f=I_{\big[ B_{\alpha},H \big]}(f)$, for a constant $c$ and since $f \in \overline{J_{h}}$,
we see that $h$ preserves flatness.

\newtheorem{c_h-flat}[c_Sym]{ Proposition }
\begin{c_h-flat} \label{c_h-flat}
If for all $\alpha$, $B_{\alpha} \prec \prec H$, we have in $\mathcal{D'}^{F}$, that $h$ preserves flatness. Further,
$\mbox{ Re }h^N \sim_{\infty} h^N$, for some $N$ over $V_{c}$, so $h^N$ preserves flatness. If also
$h^2$ is locally injective, we have that $h$ preserves flatness over $\overline{J}_{h}$.
\end{c_h-flat}

Proof: \\
Assume $r^{t}$ the restriction homomorphism to $\Omega_{t} \supset V_{c}$, so that $\Omega_{t}
\downarrow V_{c}$ as $t \uparrow \infty$. The proposition that $(rh)^{t}$ is reduced, is then a
proposition that $(rh)^{t}=const.$ on $\Omega_{t}$, that is if we let $2 i \mbox{ Im
}(rh)=$ $2(rh)-2\mbox{ Re }(rh)=$ $(rh)-(rh)^{*}$, using geometric equivalence, $\mbox{ Re }(rh)^{t} \sim_{\infty}
(rh)^{t}$. $\Box$

\vspace*{.5cm}

Assume $p$ and $\beta$ homomorphisms, with $p= \frac{\delta}{\delta x_1} + \beta$, then $I_{p}=\mbox{ ker }p$
is an ideal. If $g \in I_{p}$ we see that
\begin{equation} \label{c_p-def}
(\frac{\delta}{\delta x_1})f \big[ g-\beta g \big] + (\frac{\delta}{\delta x_1})g \big[ f - \beta f \big]=
(\frac{\delta}{\delta x_1})f(\frac{\delta}{\delta x_1})g + \beta f \beta g
\end{equation}
so $(\frac{\delta}{\delta x_1})(fg)=(\frac{\delta}{\delta x_1})f(\frac{\delta}{\delta x_1})g - \beta(fg)$.
If $\zeta$ is a zero of order $\geq 2$ for $g$, (\ref{c_p-def}) gives that
$$ \beta g \big[ \beta f + (\frac{\delta}{\delta x_1})f \big]=0$$
and if $\beta g \neq 0$ we have that $f \in I_{p}$. Let $V_2=\{ \zeta \quad g(\zeta)=\beta g(\zeta)=0 \}$,
then $(\frac{\delta}{\delta x_1})g \big[ f- pf \big]=0$, so if $\mbox{ ord }_{\zeta} g=1$, we have
$f=pf$ and thus $(\frac{\delta}{\delta x_1})f=(1-\beta)f$. Note that every geometric ideal $I(V_{c})$
can be divided into $I^{(1)}(V_{c}) \cup I^{(2)}(V_{c})$, where $I^{(2)}$ denotes $f \in I(V_{c})$,
such that $\mbox{ ord }_{\zeta} f \geq 2$, for $\zeta \in V_{c}$ and $I^{(1)}(V_{c})$ denotes $f
\in I(V_{c})$ with simple zero's. Over $I^{(2)}$, $\frac{\delta}{\delta x}$ can be considered as a
homomorphism, why it is sufficient to consider simple zero's.

\section{ On the tangents to the Weyl group, $W$ }

An analytic mapping between real vector spaces $f \text{:} E \rightarrow X$ can be factorized as
\begin{displaymath}
   \left. \begin{array}{ll}
    E \underrightarrow{\eta} &
    TS^{n}(E)   \\
    \downarrow & \downarrow{\widetilde{f}}   \\
    E^{*}  \underrightarrow{\widehat{u}} & X
   \end{array} \right.
\end{displaymath}

so that $\widetilde{f}(\eta(x))=\widehat{u}(\xi)$. Further, if for $U \subset X$, $\varphi \text{:} U \rightarrow
E$, we have an induced linear mapping $\varphi_{*}$ between the spaces of symmetric tensors and
equivalence classes of symmetric tensors, with respect to change of coordinate maps, that are called point
distributions. More precisely, given a chart $(c,t)$ in a neighborhood of $x$, $\theta_{c}
\text{:} TS^{(r)}(X) \rightarrow T^{(r)}_{x}(X)$ maps symmetric tensors on point distributions
bijectively. There is uniquely, in $T^{(0)}_{x}(X)$ an element $\epsilon_{x}=\theta_{c}(1)$,
independent on $c$, with center $x$. The field of point distributions $D_{t}$ defined by $t$,
with action by a Weyl group $W$, is given by $x \rightarrow t*\epsilon_{x}$. We have that
$$\widetilde{f}(\theta_{c}(t))=D_{t}f(x)=<t*\epsilon_{x},f>=<\epsilon_{x},\check{t}*f>=<t,f \circ \rho(x)>$$
and if $\rho$ traces $\Gamma^{U}_{A}$ and if $D_{t}f=const$, we have $t=1$, which implies
$\theta_{c}(t)=\epsilon_{x}$ and if $X$ is locally compact, we have $\epsilon_{x}=\delta_{x}$, the
Dirac-measure. According to \cite{BourL} (ch.III, Prop. 23), we can identify $V_{W}=\{ f \quad
\tau_{g}f=f \quad g \in W \}$ for $f=D_{t}f'$ with $f' \in C^{\infty}_{r}$ in a field of
point distributions, with $U(W)=$ $\{ t \in \mathcal{T}^{\infty}$ $(W)$ $\quad \text{supp} t \subset \{ 0 \} \}$
where $\mathcal{T}^{\infty}$ denotes distributions with finite support, thus $V_{W} \cong U(W)$.
Further, $T_{e}(W)$ are the primitive elements in $U(W)$, that is $T_{e}(W)=\{ f \otimes 1 + 1
\otimes f \quad f \in V_{W} \}$ and $T(W)=W \times T_{e}(W)$.

\vspace*{.5cm}

The field of point distributions defined by $t$, with action given by the Weyl-group $W$, ( $W$ can
be generated by the lineality and for instance symmetry in a canonical way ), is defined as $x
\rightarrow t*\epsilon_{x}$, where with our conditions we can assume $\epsilon_{x}$ the
Dirac-measure. If $D_{t'}h(f)=h_{1}(D_{t}f)$, then $<t',dh(f)>=h_{1}(<t,df>)$.
Concerning the Weyl-group, if $g \in W_{h_{1}}$, we have $\tau_{g}h_{1}=h_{1}$. Assume $\tau_{g}
\pi= \pi \tau_{g}$. We then have $h \pi= \pi \tau_{g} h_{1}$. If $g' \in W_{h}$, we have
$\tau_{g'}h \pi=h \pi=\pi \tau_{g} h_{1}$ implies $g=g'$. If $\tau_{g'} \pi=\pi \tau_{g}$ and
$\tau_{g} h_{1}=h_{1}$. Then $\pi \tau_{g} h_{1}=\tau_{g'} \pi h_{1}=\tau_{g'} h \pi=h \pi$,
that is $\tau_{g'} h=h$ and $g' \sim g$. Essentially, existence of lineality and equivalent
Weyl-groups, mean that the corresponding operators have the same micro-local properties.

\vsp

Assume $D_{t'}h(f)=h_{1}D_{t}(f)$, then $h_{1}(<t,d f>)=< t,h_{1}(d f)>$ and if $h_{1} d=d h_{1}$,
we have $0=< t'-t,d h(f)>+<t, d(h-h_{1})(f) >$. If $f \in \overline{J}_{h}$, we have $t'=t$. Assume
$D_{\varphi_{*}(t)}=D_{t'}h=h_{1}D_{t}$. Then $0=<\varphi_{*}(t)-t',d f>+<t',d( f-h(f))>$,
thus if $f \sim h(f)$, we have $\varphi_{*}(t) \sim t'$. Let $D_{\varphi_{*}(t)}g=D_{t'}h$, for an
algebraic homomorphism $g$, as in section \ref{sec:c_mono}. Then $f \in J^{\epsilon}_{h}$ means
 $\varphi_{*}(t) \sim t'$.

\vspace*{.5cm}

Consider the sequences
\begin{displaymath}
   \left. \begin{array}{llll}
   (I)  & \underrightarrow{D_{t}} \  X  & \underrightarrow{p} T(W) & \underrightarrow{\pi}
\ W  \rightarrow 0 \\
   \downarrow_{h} & \qquad \downarrow_{h_{1}} & \qquad \downarrow_{h_{2}} & \qquad \downarrow_{\phi}  \\
   (I') & \underrightarrow{D_{t'}} X' & \underrightarrow{p} T(W') & \underrightarrow{\pi'} \ W'  \rightarrow 0
   \end{array} \right.
\end{displaymath}
 We note that $h_2=T(\phi)$. If $\pi \text{:} W \times T_{e}(W)
\rightarrow W$ is the canonical projection, $\pi^{-1}(e)=T_{e}(W)$. We then have
$\pi'h_{2}\pi^{-1}=\phi$ and if $\phi(e)=e(=0)$, we have $h_{2}(T_{e}(W)) \subset T_{e}(W')$.
If $\phi=id$, we have $h_{2}(T_{e}(W))=T_{e}(W)$.

\section{ Some remarks }

Assume now the factorization
\begin{displaymath}
   \left. \begin{array}{ll}
   X  & \underrightarrow{\eta} S(X) \\
   \downarrow{p} & \quad \downarrow{T_{c}} \\
   T(X) & \underrightarrow{\theta_{c}} TS(X)
   \end{array} \right.
\end{displaymath}
and $p=\pi^{-1} \text{:} X \rightarrow T(X)$, using the symmetric group $S$ (cf. \cite{Bour}).
Thus $p(x)=\theta_{c}(T_{c}\eta(x))=T_{x}(X)$. The proposition that $\eta$ is constant, is then a
proposition that $p(x)=0$ or equivalently $T_{c}\eta(x)=0$. If $\eta$ is a submersion ( for
instance the orbital-mapping in the symmetric group), then using foliation, we have for all $x \in
X$, existence of a subspace $V_{x} \subset X$ such that $T_{x}(V_{x}) \oplus T_{x}(X)=\mbox{ ker
}_{x} T_{c}\eta$.

\vspace*{.5cm}

Assume $f \in J^{\epsilon}_{h}$ and if $D_{\phi_{*}(t)}g=D_{t'}h$ and if $h=g_1g$, we have over
$J^{\epsilon}_{h}$, that $g_1$ is a translation. Monotropy as we have defined it, means
particularly a factorization of an analytic homomorphism, by an algebraic homomorphism and a
translation. Thus if $h \text{:} (I) \underrightarrow{g} (I') \underrightarrow{\tau_{\epsilon}}
(I'')$ and $h_{1} \text{:} X \underrightarrow{g_1} X' \underrightarrow{\tau_{\epsilon}'} X''$.
We have, $X'=\tau_{\epsilon}X''$ and $<t', d h(f) >=\tau_{\epsilon}(<t',d g(f) >)$. Note that the same
$t'$ can be used in the representation of monotropic functions.

\vsp

Consider $h(\psi)=\tau_{\epsilon}w(\psi)$ such that
\begin{displaymath}
   \left. \begin{array}{ll}
   w(\psi) \longrightarrow   &  \alpha(\epsilon)w(\psi) \\
   \downarrow & \nearrow \\
   h(\psi) &
   \end{array} \right.
\end{displaymath}

Let $\Omega_{\lambda}=\{\epsilon \quad \frac{h(\psi)}{w(\psi)}=\alpha(\epsilon) \quad \alpha \in H
\text{ not constant} \}$ and $\Omega_{\lambda}'=\{ \epsilon \quad \frac{\tau_{\epsilon} \psi}{\psi}
\text{ not constant }\}$. When $\psi$ reduced, $\Omega_{\lambda}'$ is bounded and symmetric with
respect to the axes. If $w \tau_{\epsilon}=\tau_{\epsilon}'w$ where also
$\tau_{\epsilon'}=\tau_{\epsilon}'$, we have
$w(\frac{\tau_{\epsilon}\psi}{\psi})=\frac{\tau_{\epsilon'}w(\psi)}{w(\psi)}$. As $w$ preserves
geometric convexity we see that $\Omega_{\lambda}$ is symmetric with respect to the axes. The same
type of argument gives that it preserves boundedness. Note that $\Omega_{\lambda}'$ can be
described as the set where $\alpha$ is not a polynomial outside the origin. As $w$ maps polynomials
on polynomials as does $w^{-1}$, the same description holds for $\Omega_{\lambda}$.
The same conclusion is for $\psi \in (J_{A})$ that $h$ and $w$ have the same
micro-local properties.

\section{ Algebraicity of the ideal }
\subsection*{ The mapping $t$ }
Given an ideal $(J)$ we have if $f,g \in (J)$, then $\alpha f + \beta g \in (J)$ for constants
$\alpha,\beta$. But we do not necessarily have that $t(\alpha f + \beta
g)$ is univalent. A necessary condition for geometric convexity in the plane must be local
injectivity for $t$. Assuming this, if the ideal is defined by an algebraic homomorphism, we know
that if $\{ t \psi \}$ is a line, a closed curve or a spiral, then $\{ t h(\psi) \}$ is mapped on the same type of
geometrical object. Since also $h^{-1}$ is algebraic, the proposition can be reversed. We need however
to assume  $\{ \psi \}$ reduced locally for univalentness, but to get the same geometric object,
$h$ need not be reduced. A necessary condition for an algebraic homomorphism to preserve geometric
type of object in the plane is that the preimage is locally reduced and that $t$ is locally injective, which is
the case if we consider for instance regular approximations of singular $\psi$.

\vspace*{.5cm}

Assume $\psi \in (I)$ with a reduced pseudo-base $F_{1}, \ldots, F_{\rho}$ ordered after increasing
order of zero and $h: F_{j} \rightarrow c_{j} F_{j}$ for constants $c_{j}$. Consider $Eq_{1}=\{
F_{1}=cF_{2}  \quad \exists c \}$. As $F_{1}$ has a representation with locally isolated zero's, if $F_{1}-cF_{2} \neq 0$ on
the boundary of a small disc $\Omega$ containing a zero to $F_{1}$ then according to Rouch$\acute{e}'$s
theorem, $Eq_{1}$ must be a locally discrete set. In the same manner for holomorphic coefficients,
even if we have non-trivial foliation in the coefficients.

\vsp 

It is not difficult to prove that if $(I')$ is finitely generated through an algebraic homomorphism, then
$N(I')$ is algebraic. Let $\Omega'=N((\overline{I'}))$ (closure in Whitney-norm) and $h: (I)(\Omega)
\rightarrow (I')(\Omega')$. This must mean that $\Omega'$ is domain of holomorphy, if $\Omega$ is a
domain of holomorphy. Since $h$ maps constants on constants and polynomials on polynomials, we have
that $\Omega'$ is a Stein-domain if $\Omega$ is a Stein-domain. Note that it is sufficient for the last result, to consider
symbols on the form $1/R$, why it follows since $h$ is locally bounded. In the case where $(I')$ is infinitely
generated, we can settle for studying regular approximations.

\subsection*{ The measure $\mu$ }

The disjoint decomposition of $\Gamma=\{ \gamma(x) \}$ that we have found,
can be used to define a submersion and corresponding transversal set.
Define $$J=\{ g \quad dw(\mathcal{F} \textit{g)=d}\mathcal{F} \textit{w(g)} \}$$ and assume that $dw$ is reduced,
for a homomorphism $w$. We have $$dw \mathcal{F}-\mathcal{F} \textit{dw=d}\mathcal{F} \textit{w + B}_{\textsl{1}} +
\textit{B}_{\textsl{2}}$$
and $B_{1}=(dw\mathcal{F}$ $-d\mathcal{F}$ $w)=0$ on $J$, also $B_{2}=(d\mathcal{F}$ $w-\mathcal{F}$ $dw) \in C^{\infty}$
for $g \in \mathcal{E'}$. We thus have $dw(\mathcal{F}\textit{g})=$ $d\mathcal{F}$ $w(g) + C^{\infty}$.
For $g \in J$, we have $w(\mathcal{F}$ $g)=\mathcal{F}$ $w(g) + c$, for a constant $c$. Assume $\mathcal{F}$ $w(g)=0$
and $$w(\mathcal{F} \textit{g})(\gamma(\textit{x}))=\widehat{\textit{u}}(\xi) \text{ and }  \Gamma=\{ \xi=\gamma(\textit{x}) \}$$ then for $g \in J$,
we have $\widehat{u}=const.$ and $\mathcal{F}^{-\textsl{1}}$ $(w\mathcal{F}$ $g)=\delta_{0}$. Define
$$\Sigma'=\{ x \quad \mathcal{F}^{\textsl{-1}}\textit{w}(\mathcal{F}\textit{g})= \textit{w(g)} \quad \textit{g} \in \textit{J} \}$$
Since $\widehat{H*'\varphi} \rightarrow \widehat{H}$ as $\widehat{\varphi} \rightarrow 1$ and analogously for $\widehat{H*''\varphi}$, we see that $\Sigma=\Sigma'$. Thus, for $x \in \Sigma$, we have $w(g)=\delta_{0}$ and for $x \in \Delta^{Q}_{x}$, we have $w(g)=0$. Further, if $h(g)=w(\mathcal{F}$ $g)-\mathcal{F}$ $w(g)$, we have $h(g)=0$ on $\Sigma$ and if $g \in J$, we have $h(g)=const.$ on $\Delta^{Q}_{x}$.

\vspace*{.5cm}

Assume $\mu=\mathcal{F}$ $w-w\mathcal{F}$ a positively definite measure for an algebraic homomorphism $w$.
We have the following proposition.

\newtheorem{e_Rycket}[e_sgn]{ Proposition }
\begin{e_Rycket}
The measure $\mu$ is finitely generated for algebraic homomorphisms $w$ and infinitely generated
for analytic homomorphisms on the form $h=\tau_{\epsilon}w$.
\end{e_Rycket}

Proof:\\
According to R\"uckert's Nullstellensatz, if
$<I_{\Delta}(\varphi),\psi>=<\delta_{x}(\varphi),\psi>$ for $\varphi,\psi \in H(V)$ where $V$ is a bounded
neighborhood of $0$, then $<I^{N}_{\Delta}(\varphi),\psi>=0$ for large $N$. It is not difficult to
prove that there is a $g_{1} \in H$ such that \\ $\mathcal{F}$ $w(g^{N})=\mathcal{F}$ $\big[ w,I_{\Delta}
\big](g_{1}^{N})=\mathcal{F}$ $w(I_{\Delta}^{N}(g_{1}))=constant$. In the same manner there is a $g_{0}
\in H$ such that $w \mathcal{F}$ $g^{N}=\big[ w, I_{\Delta} \big](\mathcal{F}\textsl{g}_{\textsl{0}})^{\textsl{N}}=$ $0$. Note
that if $w$ is algebraic $\big[ \widehat{w},I_{\Delta} \big]=\big[ w,\widehat{I}_{\Delta}
\big]=\big[ \widehat{I}_{\Delta},w \big]=\big[ I_{\Delta},\widehat{w} \big]$.$\Box$

\vspace*{.5cm}

In the case of $h=\tau_{\epsilon}w$, to prove that the measure $\mu$ is infinitely generated,
it is sufficient to consider the one-variable case. For instance, $\mu=\big[ \widehat{I}_{\Delta},\tau_{\epsilon}w
\big]-\big[ \widehat{I}_{\Delta},e^{ix \cdot \epsilon}w \big]$. Let $T_{1}=\frac{d}{d x}-ix$, $T_{2}=\frac{d^{2}}{dx^{2}}+x^{2}$ and so on.
Thus $\mu \sim \sum \frac{\epsilon^{n}}{{n!}}T_{n}w$
and we see that $\mu$ is infinitely generated for $\epsilon >0$. It is clear that $\mu$ modulo
monotropy is finitely generated.$\Box$

\vspace*{.5cm}

Assume now that $\mu$ is algebraic in the tangent space, that is
$d \mu^{2}(g)=d \mu(g^{2})$ and consider $\mu_{1}=\mathcal{F}$ $d \mu-d \mu$ $\mathcal{F}$. Thus $\mu_{1}$
is finitely generated. The proof is completely analogous with the one given above and further it
follows that if $w$ is algebraic in the tangent space, then $d \mu$ is finitely generated. Let
$\mu_{(2)}$ be the measure corresponding to $w^{2}$, then $d \mu_{(2)}(g) \sim_{m} d \mu(g^{2})$ for $w$
algebraic. That is assume $g_{1}*g_{1}$ a regular approximation of $g^{2}$, such that
$\widehat{g_{1}} \sim_{m} \widehat{g}$. Then
$$d \mu_{(2)}(g) \sim_{m} d \widehat{w^{2}}(g)-d w^{2}(\widehat{
g}_{1}) \rightarrow d \widehat{w^{2}}(g)- d w(\widehat{g^{2}})=d(\widehat{w}(g^{2})-w(\widehat{g^{2}}))$$
Given existence of $g_{1}$ as above we have that $\mu_{(\cdot)}$ is algebraic in the tangent
space. Consider $h=\tau_{\epsilon}w$ with $w$ algebraic. We have earlier considered
$\Phi_{\mu}(\gamma)=\int_{\gamma} d \mu(\gamma)$ with $\gamma=(\psi,h(\psi))$ and $\psi$ reduced.
This is finitely generated and we have an estimate $\mid \Phi_{\mu}(\gamma) \mid \leq C \mid \psi
\mid^{\rho}$ for the Puiseux-index $\rho$ and a constant $C$.

\vspace*{.5cm}

Assume $\gamma=(\psi,h(\psi))$ with $\psi$ reduced (or $h$ reduced). If $\gamma \rightarrow P$,
then $P$ is on the curve and also $\overline{P}$ is on the curve $\overline{\gamma}$. Assume
that $\overline{P}$ singular for $\gamma$. Then $d \gamma(\overline P)+d \overline{\gamma}(P)=0$.
This means that $\mbox{ Re }\gamma$ is singular in $P$. Conversely, if $P$ is singular for $\mbox{ Re
}\gamma$ and $d \gamma \rightarrow p$, then $p$ must be purely imaginary. In the same manner if $P$
is singular for $\mbox{ Im }\gamma$ and $d \gamma \rightarrow p$, then $p$ must be purely real. The
conclusion is that if either $\mbox{ Re }\gamma$ or $\mbox{ Im }\gamma$ has $P$ as singular point,
then it can not be singular for the other unless it is singular for $\gamma$.

\vspace*{.5cm}

Denote $dh-hd=\big[ dh \big]-\big[ dh \big]^{*}$. If $h=h^{*}$ we have locally $dh-hd \sim 2 i
\mbox{ Im }\big[ dh \big]$. Sufficient conditions for a regular approximation are $h(d g) \neq 0$
($\Rightarrow d g \neq 0$) and $dh(g) \neq 0$. Assume $g$ satisfies these conditions and

\begin{equation} \label{e_imag}
dh(g)-h(d g) \sim \mbox{ Im }\big[ dh \big](g)
\end{equation}

If $h \sim \overline{h}$, then locally $d \big[ \mbox{ Im }h \big] \sim \mbox{ Im }\big[ dh \big]$.
Through the conditions for the dynamical system, if $P$ is singular for $\mbox{ Im }\big[ dh \big]$
that is $((\ref{e_imag})=0)$, then $P$ can not be singular for $\gamma$ if it is not singular for
$\mbox{ Re }\gamma$. Thus if we assume real dominance, then singular points are given by the ones
for $\mbox{ Re }\gamma$. If $h^{2}(g)=h(g^{2})+F$, where $d F \sim_{m} 0$, then there are $d F_{1}$
of type $-\infty$ such that $dF - dF_{1}$ is of type $0$. We could say using the generalized Cousin integral
representation, that $h^{2}(g)=h(g^{2})$ in the tangent space modulo regularizing action.

\subsection{ Algebraicity in the tangent space }
If $h$ is not algebraic but algebraic in the tangent space, we have
$$ d h^{-1}(g)=d h(1/g)=dh(\frac{1}{1-g'})=d \sum_{j} h^{j}(g')=d \sum_{j} h({g'}^{j})$$
If $h$ is algebraic in the tangent space, the corresponding $\eta'$ is algebraic in the sense that
${\eta'}^{-1}(g)=g d h^{-1}(g)=(1-g')d \sum h^{j}(g')$. We can form the ideals $$(J_{\lambda}')=\{ g \quad dh(g) \geq \lambda g\}$$
for $\lambda > 0$. We have that $d h^{j}(g')/g'$ does not change sign over
$(J_{\lambda}')$. If $g'$ is real then ${g'}^{2} \geq 0$ and also $g' d \sum h^{j}(g')$ does not change
sign over $(J_{\lambda}')$. For ${\eta'}^{-1}(g)$ we have constant sign on $(J_{\mu}') \cap V_{\pm}$
where $V_{\pm}=\{ g \quad \pm g >0 \}$. Thus ${\eta'}^{-1}(g)$ only changes sign on $\{ g=0 \}$.

\section{ Analytic sets and symmetry }

\subsection{ A weighted lineality }
Assume $h(f)=P/Q$ for polynomials $P,Q$, where $Q$ is assumed hypoelliptic and self-adjoint.  Let
$I_{G_{\tau_{\eta}}}=\tau_{\eta}I_{H}\tau_{\eta}$. Thus
$$ \frac{Q(x+\eta)}{Q(x)} I_{G_{{\tau}_{\eta}}} = \Big( 1 + \sum_{\alpha \neq 0}
\frac{Q_{\alpha}(x)}{Q(x)} \eta^{\alpha} \Big) I_{G_{{\tau}_{\eta}}}$$
for $\mid \eta \mid < R$, $\eta$ real and $R$ finite and where $Q_{\alpha}/Q \rightarrow 0$ as $\mid
x \mid \rightarrow \infty$ . Note that the condition $Q$ hypoelliptic, implies that there are $\mid x \mid$ large
such that $Q(x+\eta) \neq Q(x)$ for $ \eta \neq 0$. We now have
$$ I_{G_{{\tau}_{\eta}}}(f)-P/Q=\Big( \sum_{\alpha \neq 0} \frac{Q_{\alpha}(x)}{Q(x)} \eta^{\alpha}
\Big) I_{G_{{\tau}_{\eta}}}(f)$$ We now introduce a "weighted" lineality
$\Delta^{Q}_{x}=\{ \eta \quad \frac{Q(x+\eta)}{Q(x)} I_{G_{{\tau}_{\eta}}}=I_{H} \quad \forall x \}$
where the term $\frac{R}{Q}h(f)$ is assumed insignificant. Let
$\tau_{\eta}I_{H}\tau_{\eta}=\widehat{g}_{\eta}$, so that
$ \Delta^{Q}_{x}= \{ \eta \quad \sum'_{\alpha} \Big( q_{\alpha}(x)\eta^{\alpha} \Big)
\widehat{g}_{\eta}(f)=P/Q \}=\{ \eta \quad \sum_{\alpha}' q_{\alpha}(x)\widehat{(D^{\alpha}g_{\eta})}(f)=P/Q \}$,
where $\Sigma'$ denotes a finite sum, since we have assumed $Q$ a polynomial. Further, we assume
$q_{\alpha}(x) \rightarrow 0$ as $\mid x \mid \rightarrow \infty$ and $\mid \eta \mid < R$.
Let $\tilde{\Delta}^{Q}_{x}=\{ \eta \quad \sum_{\alpha}' q_{\alpha}(x) \big( D^{\alpha}_{\eta}g_{\eta}
\big)(f)=(P/Q)\delta_{0}(\eta) \}$. Thus, $I_{G_{{\tau}_{\eta}}} - P/Q=\sum_{\alpha}' q_{\alpha}(x)
\widehat{ (D^{\alpha}_{\eta}g_{\eta} }(f))$. Note that $\{ \eta \quad \sum'=0 \quad \forall x \}$ corresponds to
the set $\Delta_{x,y}(h)$ in the Weyl-calculus, that is invariance for translation with $\eta$ in both
variables separately for the kernel $H$. Further if $\eta$ is considered on a line-segment $L$,
then $\tilde{\Delta^{Q}_{x}}$ is $\Delta^{Q}_{x}$ considered on the dual line-segment $\Gamma$,
$\mathcal{F}$ $r_{L}= r_{\Gamma} \mathcal{F}$ in $Exp_{\rho,A}$.

\subsection{ R\"uckert for monotropic homomorphisms }

Assume $h(f)=\tau_{\epsilon}g(f)$ as earlier and $h(f)$ with a zero in a point $a \in
\mathbf{C}^{\textsl{n}}$, which is not of infinite order. We assume $g(f) \neq 0$ in $a$ and in a
neighborhood of $a$. Consider now translation along a transversal emanating from the point $a$.
We have seen that $\tau_{\epsilon}^{N}g(f)=0$ is modulo monotropy equivalent with
$\tau_{\epsilon}g^{N}(f)=0$ and further $h^{N}(f)=0$. (Note that $\epsilon$ depends on $N$,$f$ and
$a$). Assume $N(J_{g})$ the zero-set to the ideal corresponding to $g$, is an analytic set. We assume
that it has positive dimension locally (that is near the translated $a$). Then there is an irreducible,
given a stratification of $N(J_{g})$, of dimension $>0$. This gives a neighborhood of the point $a$,
where $h^{N}(f)=0$

\newtheorem{d_inf-zero}{ Lemma }[section]
\begin{d_inf-zero} \label{d_inf-zero}
If $h$ is an analytic homomorphism with a zero in a point $a$, with factorization
$h(f)=\tau_{\epsilon}g(f)$, where $g$ is an algebraic homomorphism with $N(J_{g})$ locally of positive
dimension, then $h^{N}$ has an infinite zero in $a$, for $N \geq N_0$, with $N_0$ positive.
\end{d_inf-zero}

\subsection{ The set of symmetry, $\mathbf{\Sigma_{x,y}}$ }
Let $\Sigma_{x,y}=\{ \eta \quad H(x+\eta,y-\eta)=H(x,y) \quad \forall x,y \in \mathbf{C}^{\textit{n}} \}$ for the kernel to a
continuous linear operator $h$ on $L^2 \cap H$. Assume to begin with that the kernel $H(x,y)$ is in $L^2
\cap H$, such that $h$ is compact. Then $H(x+\eta,y)=H(x,y+\eta) \Leftrightarrow$ $\eta \in \Sigma_{x,y}$
 for all $x,y$. Particularly, $\Delta_{x,y} \subset \Sigma_{x,y}$ and with
the conditions on $H(x,y)$ above, we have that $h$ has slow oscillation over $\Sigma_{x,y}$. If
$\Sigma^{*}_{x,y}$ is the set corresponding the adjoint operator $h^{*}$, we have for $h^{*}=h$
that $\Sigma^{*}_{x,y}=\Sigma_{x,y}$. Assume
$$\Sigma^{As}_{x,y}=\{ \eta \quad H(x+\eta,y)=H(x,y+\eta) \quad \mid x \mid \rightarrow \infty, \quad \mid y \mid \rightarrow \infty
\}$$ and correspondingly for $\Delta^{Q}_{x,As}$, we then have for $\eta$ such that $\mid \eta \mid
<R$, $R$ finite, that $\Sigma^{As}_{x,y}=\Delta^{Q}_{x,As}$, if $I_{G_{\tau_{\eta}}}=\tau_{\eta} I_{H} \tau_{\eta}$ is bounded.
Assume $$\widetilde{\Sigma}^{As}_{x,y}=\{ \eta \quad H(x+\eta,y)=H(x,y+\eta) \quad x,y \in V \}$$ where $V$ is
a set of positive dimension containing the infinity. If $\eta \in \Sigma_{x,y}$ and $h^{*}=h$, we
must have $\eta \notin \Delta^{Q}_{x}$ and conversely if $\eta \in \Delta^{Q}_{x,As} \backslash
\Delta^{Q}_{x}$ then $\eta \in \Sigma^{As}_{x,y}$.

\newtheorem{d_As}[d_inf-zero]{ Lemma }
\begin{d_As}
If $\eta \in \widetilde{\Sigma}^{As}_{x,y}$ and $h^{*}=h$ over $L^2 \cap H$, then $\eta \in \Sigma_{x,y}$
\end{d_As}
Proof: \\
Assume $h$ is such that $I_{G_{\tau_{\eta}}}=I_{H}$ on a set of positive
dimension in the infinity, for $\eta$ finite, then the proposition follows from Lemma
\ref{d_inf-zero} $\Box$.

\vspace*{.5cm}

Assume $h_1(x,y,\eta)=H(x+\eta,y)-H(x,y+\eta)$. If $h_1$ is entire in $(x,y)$, we must have for
$\eta \in \Sigma^{As}_{x,y}$, that $h_1$ is identically $0$. Otherwise, we assume $h_1$ holomorphic
in a domain $\Omega$ containing the infinity. If $h_1$ is $0$ on a set $V \subset \Omega$ of
positive complex dimension (we assume $V$ contains the infinity), then $h_1=0$ on $\Omega$. The
definition of $H$ can be extended outside $\Omega$ to a constant, or for
$H$ developed in a pseudo-base, with constant coefficients, to an entire function.
Again, if $\eta \in \Sigma^{As}_{x,y}$, where $h$ is according to the conditions in Lemma \ref{d_inf-zero},
then there is a $N$ sufficiently large, such that
$(I_{G_{\tau_{\eta}}}-I_{H})^N=0$ on a set of positive dimension close to the infinity. Taking
monotropy in to account, this means that $I_{G_{\tau_{\eta}}}^N \sim_{m} I_{H}^N$.

\vspace*{.5cm}

The conclusion is thus, that for $h=h^{*}$ bounded over $L^2 \cap H$, if we consider the sets corresponding to
the iterated operators $h^{N}$,
\begin{equation} \label{d_sym} \Delta^{Q}_{x,As} \sim_{m} \Delta^{Q}_{x} \cup \Sigma_{x,y} \text{ a disjoint
union } \end{equation}
Note that given the division $\eta \in \Delta^{Q}_{As}$, $\eta=\eta_1+\eta_2$ with $\eta_2 \in
\Sigma_{x,y}^{As}$ and if we assume $g_{\eta}=\tau_{\eta}F=(\tau_{\eta} I_{G'} \tau_{\eta}=)\tau_{\eta_2}\tau_{\eta_1}F$, where $F$ is a
self-adjoint integral operator with kernel representation such that $F\sim h^{N}$, we have the representation $g_{\eta}=\mbox{ Re
}(\tau_{\eta_1}F)$ for $\eta \in \Delta^{Q}_{x,As}$.

\vspace*{.5cm}

In the case where $h$ does not have a kernel in $L^2 \cap H$, we still assume $\eta \in
\Delta^{Q}_{x,As}$. This means particularly, that given that the translation is bounded, for $\mid
\eta \mid$ finite, $h$ and $h^{*}$ have slow oscillation relative $Q$. Assume instead that $h^{-1}$ is bounded,
then obviously $\Sigma^{As}_{x,y} \subset \Delta_{x,As}^{Q}$. If we assume $h^N$ locally
injective, for large $N$, we also have the opposite inclusion, $\Delta^{Q}_{x,As} \subset
\Sigma_{x,y}^{As}$. Further, we know that under these conditions, if $\eta \in \Delta^{Q}_{x,As}$,
then for $t$ real, $t \eta \in \Delta^{Q}_{x,As}$ as $t \rightarrow \infty$ since the translation is a
compact operator over reduced elements. Naturally, the decomposition (\ref{d_sym}) holds for the
iterated operators also without the condition that the operator is bounded on $L^2 \cap H$.

\vspace*{.5cm}

Assume now $h$ is not self-adjoint, but such that $h \sim_{\infty} \mbox{ Re }h$. Consider
$I_{G_{\tau_{\eta}}}$ as a function of a complex $\eta$ and let $G_{R}$ be the kernel corresponding
to $\tau_{\eta} \mbox{ Re }I_{H} \tau_{\eta}=I_{G_{R}}$
and in the same way, let $H_{R}$ be the kernel corresponding to $\mbox{ Re }I_{H}$. Both operators are
considered as locally defined in $z,w$. Thus,
$G_{R}(z,w;\eta)=\frac{1}{2}( H(z+\eta,w-\overline{\eta})+\overline{H}(z+\eta,w-\eta))$. Over $t \eta \in
\Sigma( \mbox{ Re }h )=\Sigma^{*}( \mbox{ Re }h )$, for all real $t$, we have
$$G_{R}(z-\eta,w+\eta;\eta)=H_{R}(z,w) \Leftrightarrow$$  $$\frac{1}{2}( H(z,w+2 i \mbox{ Im }\eta) + \overline{H}(z,w) )=\frac{1}{2}(
H(z,w) + \overline{H}(z,w) )$$ The necessary condition is thus that $2 i \mbox{ Im } \eta \in
\Delta({}^tI_{H})$. Conversely, if $2 i \mbox{ Im }\eta \in \Delta_{z,w}(I_{H})$, where $I_{H}$ is
locally defined in $z,w$, then $t \eta \in \Sigma_{z,w}(\mbox{ Re }h)$, for all real $t$.

\subsection{ Behavior at the infinity }
Note that for $\eta \in
\Sigma^{As}_{x,y}$, we have that a constant value is preserved for $h_1$ in $(x,y)$, in the sense of
Cousin \cite{Cous}, why it is natural to consider the quotient-topology $I_{h_{1}}(f)=P/Q$, and
$f=P'/R$.

\newtheorem{d_pol-infty}[d_inf-zero]{ Lemma }
\begin{d_pol-infty}
If $P$ is a polynomial with a real zero in the infinity, then the set of complex zero's to $P$ close to the
infinity does not have positive complex dimension.
\end{d_pol-infty}

This means that $P$ does not have an infinite zero in the complex sense. If $h$ is such that $h^2$
is injective (assume $f$ reduced), the proposition that the zero in the infinity has complex
dimension is that it is a regular point in the zero-set to $h(f)=P/Q$. 
Assume $f=P'/R$ an entire function, then there exists $h$, a differential operator or $h \sim_m$
differential operator, such that $f \in \overline{J}_{h}$. Assume further that $R^{*}=R$ and
hypoelliptic. The local condition
\begin{equation} Q(D)(P'/R)=constant \label{d_eq} \end{equation}
means that, for an appropriate test-function $\varphi$, $\sum_{\alpha} Q_{\alpha}\varphi D^{\alpha}(P'/R)=
\sum_{\alpha,\beta} Q_{\alpha}\varphi
\frac{{P'}^{(\beta_1)}R^{(\beta_2)} \ldots R^{(\beta_s)}}{R^{(s)}}=constant$. Assume ${P'}^{(\beta_1)}/R
\rightarrow 0$ as $\mid \xi \mid \rightarrow \infty$ for all $\beta_1 \neq 0$ (slow oscillation
relative $R$). A necessary condition for the condition (\ref{d_eq}) is that $P'$ and $Q$ are in some sense equivalent.
We will settle for $P' \sim_{m} R$, that is $f \in J^{\epsilon}_{h}$ and where $h(f)=constant$ we have $\mid 1 - g(P'/R)
\mid < \epsilon$, for an algebraic homomorphism $g$. Particularly, if $g=1$, $P' \sim_{m} R$.
We can now study the ideals of polynomials $$(J_{Q})=\{ P \qquad \exists R \mbox{
hypoelliptic and self-adjoint } P \prec \prec R \quad P \sim_{m} R \}$$ and
$$(J_{Q}')=\{ P \qquad \exists R \mbox{
hypoelliptic and self-adjoint }\delta P \prec \prec R \quad P \sim_{m} R \}$$

\subsection{ Some results for $(J_{Q}')$ }

Assume $f$ symmetric and $Df$ corresponding to a regularizing operator, for instance
$Df \prec \prec I$, with $f=P/Q$, where $D=D_{x_{1}} \ldots D_{x_{n}}$.

\newtheorem{d_rel_HE}[d_inf-zero]{ Lemma }
\begin{d_rel_HE}
Assume $f$ as above and with $Q$ hypoelliptic and self-adjoint. Then $D^{\beta_j}P \prec \prec Q$ for all
$\beta_j \neq 0$ $\Leftrightarrow$ $\mbox{ Im }P \prec \prec Q$
\end{d_rel_HE}

Proof:\\
It is not difficult to prove that for
$C_{\varphi,f}=f \varphi - \varphi f$, where $\varphi$ is a real test-function, arbitrary in
$C^{\infty}_{0}$,
\begin{equation} \label{d_comm}
C_{\varphi,f} \prec \prec I \quad \forall \varphi \in C^{\infty}_{0} \Leftrightarrow D^{\beta_j}f \prec \prec I \quad \forall \beta_j \neq 0
\end{equation}

Particularly, we have that,

\newtheorem{d_rel_smooth}[d_inf-zero]{ Lemma}
\begin{d_rel_smooth}
If $f$ has finite type and $f=R/S$, where $S=S^{*}$ and hypoelliptic, if $Df \prec \prec I$,
then $\delta R \prec \prec S$.
\end{d_rel_smooth}

If $Q$ is real, we obviously have $\mbox{ Im }f \sim_{\infty} (\mbox{ Im }P)/Q$. If $Q$ is only self-adjoint,
we have if we assume $P \prec Q$ and $Q$ hypoelliptic, we get the same conclusion.
With these conditions on $Q$ further, $D(\frac{P}{Q}) \sim_{\infty} \frac{D P}{Q}$ and through the
conditions on $f$, $f \sim_{\infty} \overline{f}$, so if we assume $(D^{\beta_j}P)/Q \sim_{\infty}
D^{\beta_{j}}(\frac{P}{Q}) \prec \prec I$, for all $\beta_{j} \neq 0$, then according to
(\ref{d_comm})
$C_{\varphi,f} \prec \prec I$ with $\varphi$ arbitrary in $C^{\infty}_{0}$, why $\mbox { Im }f
\sim_{\infty} ( \mbox{ Im }P)/Q \prec \prec I$. The converse follows in the same manner. Thus,
$\mbox{ Im }f \prec \prec I$ implies $C_{\varphi,f} \prec \prec I$ and according to (\ref{d_comm}),
$ D^{\beta_{j}}f \prec \prec I$ for all $\beta_{j} \neq 0$ and we are done. $\Box$

\vspace*{.5cm}

Note that if for $P,Q$ polynomials, $P-Q$ has a real zero in the infinity and if $(P-Q)^N$ has a
zero of complex dimension in the infinity, then $(P-Q)^N$ can not be a polynomial. However,
$(P-Q)^N$ has a ps.d.o realization with a regularizing action and with a polynomial part.
Note that if $P_j$ is reduced for all $j$, we only have to deal with isolated zero's, but $\sum_j
P_j$ may well have a zero of infinite order. The condition, $\sum_j P_j \in rad(I_{RED})$ implies
that we do not have a zero in the infinity, but $(\sum_j P_j)^N \sim_{m} \sum_j P_j^{\mu_j}$, for
iteration indexes $\mu_j$, where the right side may have a zero of infinite order. Assume $V$ a set
of positive complex dimension containing the infinity and $h(f)=\sum_j P_j$, such that $h$
satisfies the conditions in Lemma 5-4, then $\sum_j h^{N}(P_j)$ has an infinite zero in the
infinity and we can assume $h^N(f) \in J^{\infty}_V$. Note that $J^{\infty}_{V}$ denotes the ideal of
symbols with an infinite zero on a set of positive complex dimension $V$. This corresponds to an operator with
regularizing action in a ps.d.o realization.

\subsection{ Behavior in the origin }
\label{section:d_org}
Assume $g=R/S$, with $S$ hypoelliptic and self-adjoint and $Dg \rightarrow 0$ as $\mid x \mid
\rightarrow \infty$. Let $J_{\infty}=\{ g \quad D(f-g)=0 \quad \exists f \in J^{\infty}_V \}$,
that is, $g=f+c$ for a constant $c$ in the infinity. The conditions on $g$ mean particularly that
$g \in (I)_{\mathcal{Q}}'$. Note that $J_{\infty}$ is locally defined and that $V$ is assumed to have
positive complex dimension and that $V$ contains the infinity. Particularly, if $g \rightarrow 0$
as $\mid x \mid \rightarrow \infty$, we have that $g \in J_{\infty}$ can be approximated by
$J^{\infty}_V$.

\vspace*{.5cm}

Assuming $zP(z)=F(\frac{1}{z})$, for a polynomial $P$ and $F$ holomorphic for $0 < \mid z \mid <A <
\infty$. If $F$ has an infinite zero in the infinity, we can not draw a conclusion concerning the
order of zero for $P$ in the origin, since this corresponds to an essential singularity in the
general case, as $z \rightarrow 0$. However, since $zP(z) \rightarrow 0$ as $z \rightarrow 0$, for all
polynomials $P$, we could at least formally say that if $P \in J_{\infty}$, $zP(z)$ has a
zero of infinite order in the origin.

\vspace*{.5cm}

Assume $V$ a complex domain that contains the infinity and $u \in C^{0}(V)$. We know that $\mbox{
log }\mid u \mid$ is a sub harmonic function on $V \backslash Z_u$ and that $\mbox{ log } \mid u \mid \equiv - \infty$ on
$Z_{u}$. Let $u=e^{g}$, then $Z_{u}$ becomes the polar set to $-g$ ($u$ is also holomorphic on
$V$). Let $(I)=\{ g \quad e^{g} \rightarrow 0 \quad \mid x \mid \rightarrow \infty \}$.
For a complex-valued $g$, the corresponding ideal is not closed, but we can assume $g \sim_{\infty}
\mbox{ Re }g$ or consider $rad(I)$, to get a closed ideal.

\vspace*{.5cm}

Assume $\widetilde{h}(e^{\varphi})=e^{h(\varphi)}$ with $\widetilde{h}(f)=\tau_{\epsilon}\widetilde{g}(f)$ and
$f=e^{\varphi}$. Then we have that $N(\widetilde{h}(f))$ is the polar set to $h(\varphi)$. A
comparison between the dimensions for $N(\widetilde{h}(f))$ and $N(\widetilde{g}(f))$ corresponds to a comparison
between the dimensions for the polar sets to $h(\varphi)$ and $g(\varphi)$ respectively. Note that we may have
$\mbox{ dim }V=1$ for $V \sim_{m} V'$ and $\mbox{ dim }V'=0$, this because $\epsilon$ depends on
the choice of point $x$. Further, note that
$D_{t}\widetilde{h}(f)=\tau_{\epsilon}D_{t}\widetilde{g}(f)=0$, means $\widetilde{h}_{*}(t)=\widetilde{g}_{*}(t)=0$ or
$d h(\varphi)=0$. Further, $\widetilde{g}(f)=\widetilde{g}\widetilde{h}^{-1}(e^{h(\varphi}))$, so
the latter condition is $d h(\varphi)=\tau_{\epsilon} d( h^{-1}gh(\varphi))=0$. If
$\widetilde{h}$ is self-adjoint, we have $d\widetilde{h}(f)=0$ iff $f=0$ without reducedness for
$h$. Further, $\mbox{ dim }N(f)=\mbox{ dim }A$, where $A=\mbox{ sng } N(\widetilde{h}(f))$.

\subsection{ Schwartz-type topology on the phase-ideal}
 We will in this section refer to \cite{Sjo} for notation and results. A contour $\Gamma$, is a
 $C^{\infty}$-mapping $W \rightarrow \mathbf{R}^{\textit{n}}$, where $W \subset \subset \mathbf{R}^{\textit{n}}$
 is open and we have $\Gamma$ 1-1 and also $d \Gamma$ 1-1. A contour is ``good'' for a function $\varphi$,
 if $\tau_{\eta}\varphi - \varphi \leq -C \mid \eta \mid^2$ for a positive constant $C$ and for
 $\eta$ traces $\Gamma$. Here $\varphi$ is assumed continuous $\Omega \rightarrow \mathbf{R}$, with
 $\Omega \subset \mathbf{C}^{\textit{k}}$. For $u \in H^{loc}_{\varphi,x}$, we define
 $I_{\Gamma}(\lambda)=e^{- \lambda \varphi(0)} \int_{\Gamma} e^{\lambda \varphi(y)}u_{\lambda}(y)d y$. Over such $\varphi$, we claim that $\tau_{\eta}$ is a compact operator. Note that $\varphi$ can be regarded as the monotropic function to a corresponding holomorphic function, further that the germs are
locally finite-dimensional. Further, if $f \sim_{m}g$ and $g \rightarrow 0$, along complex lines, then
$\mid f-g \mid \leq Q_{\epsilon}(\frac{1}{t})$, for $Q(0)=\epsilon$, as $t \rightarrow \infty$.
Particularly, we have $f \sim g$ in $H^{loc}_{\varphi}$, that is the difference gives an exponentially
small contribution to the integral $I_{\Gamma}$.

\vspace*{.5cm}

The condition that $\eta \in \Delta^{Q}_{x,As}$ implies that the translation is weakly compact over $h$
and $\varphi \in (I)_{\mathcal{Q}}'$, that is $\frac{Q(x + \eta)}{Q(x)} \tau_{\eta} h (\tau_{\eta} \varphi) \rightarrow h(\varphi)$
as $\mid x \mid \rightarrow \infty$ and symmetrically as $\mid y \mid \rightarrow \infty$. Thus $\tau_{\eta} \varphi$
has to converge as $\mid x \mid \rightarrow \infty$ and $\tau_{\eta}$ is weakly closed over $\varphi$. By
realizing $\tau_{\eta}$ as a measure over $(I)_{\mathcal{Q}}'$, $\tau_{\eta}$ can be considered as weakly compact.
Particularly, the extension used to define the Schwartz-type topology, is weakly compact.

\vspace*{.5cm}

Note that, if $T_{\eta}(h)=\tau_{\eta}h\tau_{\eta}$, for $\eta \in \Delta^{Q}_{x,As}$, we have that
$T_{\eta}$ is a compact operator, even though $\tau_{\eta}$ is not compact unless $h$ is reduced. We define
$\Delta^{Q}_{x}=\{ \eta \quad (T_{\eta}-I)(Qh)(x)=0 \quad \forall x \}$, where
$T_{\eta}$ is assumed compact over $(Qh)$. We have seen that modulo monotropy,
$\Delta^{Q}_{x}=0$, for the set corresponding to $h^{N}$, for all $N \geq N_{0}$, for some $N_{0}$.
Thus, $\Delta^{Q}_{x}$ has, for all $N$, representation with locally isolated points and can
be realized by a locally injective, also $C^{\infty}$ mapping. As $\Delta^{Q}_{x} \cap d \Delta^{Q}_{x} \subset \Delta^{Q}_{x}$, as an analytic set, the same argument can be applied on this set. For an operator on $L^{2}$, on the form $A=T_{\eta}-I$, under the assumption that $\eta \notin \Delta^{Q}_{x}$, we have $N(A^{N}_{\eta})=N(A^{N+1}_{\eta})$, $N>N_{0}$, for some $N_{0}$. Modulo monotropy, we have $\varphi \in N(A^{N}_{\eta})$ implies $\eta \in \Delta^{Q}_{x,(N)}$ and through the conditions, we see that
$\eta=0$.

\vspace*{.5cm}

Assume $h$ localized to $2 \mbox{ Re }h=h(x + \overline{\eta}) + \overline{h}(x+\eta)$ as earlier.
We then know that $\Delta^{Q}_{x}(\mbox{ Re }h)=\{ 0 \}$. Thus, $\eta \in \Delta^{Q}_{x}(h^{*})
\Leftrightarrow \overline{\eta} \in \Delta^{Q}_{x}(h)$. If we have $h^{*}=h$, then $\Delta^{Q}_{x}$
is symmetric with respect to the real axis. Further, $\tau_{\eta}Qh=Qh\tau_{-\eta}$, why $\Delta^{Q}_{x}$
is symmetric with respect to the origin. Assume $\Delta^{Q}_{x} \subset U_{c}$, where $U_{c}$ is a convex
neighborhood of the origin. This neighborhood can be seen as a porteur
for $S_{\tau_{\eta}}$, and $\widehat{S_{\tau_{\eta}}}$ has indicator $-\infty$ over the porteur and $(I_{1})$.

\subsection{ A generalization of the Cousin integral }

The problem is to determine if $I_{\Gamma}$ can be used to generalize the Cousin integral. Assume $\varphi(x)=\mbox{ log }(x-y)$. If $\Gamma$ describes a circle-segment, then $\Phi=\int_{\Gamma} \big( \frac{\delta}{\delta x} e^{\lambda \varphi}\big) u_{\lambda}(x,z) d x$ is on the form of a Cousin integral. $\Phi$ is a representation of $F_{\lambda}=(x-y)^{\lambda}u_{\lambda}$ and if $u_{\lambda} \neq 0$ in a point $y$, then $F_{\lambda}$ has a zero of order $\lambda$ in this point. Further, $\Phi$ has analytic continuation along a transversal to the circle-segment $\Gamma$. Note that if $\Gamma$ has a disjoint decomposition as an analytic set, for instance using properties for the phase-function, then the representative $\Phi$ has a corresponding decomposition. Assume the phase
defined, modulo monotropy, through a homomorphism, such that $\Gamma=\Delta^{Q}_{x} \cup \Sigma_{x,y}$ disjointly and $\varphi^{k} \sim_{m} \tau_{x'} \mbox{ Re }\tau_{x'}\phi$,
for an analytic $\phi$ and $x' \in \Delta^{Q}_{x}$, then
$$ \Phi=\int_{\Gamma} (\frac{\delta}{\delta_{x}}e^{\lambda {\varphi}^{k}(x)})u_{\lambda}(x,z)d x=\int_{\Sigma_{x,y}} \int_{\Delta^{Q}_{x}}(\frac{\delta}{\delta_{x}}e^{\lambda \tau_{x''} \mbox{ Re }\phi(x')} )u_{\lambda}(x,z)d x$$
Note that the integral $\int_{\Gamma} e^{\lambda h^{k}(\psi)}u_{\lambda}d x$, corresponds to $I_{\Gamma}I_{H}^{k}(u_{\lambda}(x,z))$, where $h^{k}=\mbox{ Re }h^{k}$.

\vspace*{.5cm}

Assume $\varphi=h(\psi)$, where $\varphi \in (I)_{\mathcal{Q}}'$ and where $h$ is assumed self-adjoint in $L^2$ and self-transposed in the duality between $H'$ and $H$. We then have $\varphi^{k} \sim_{m} h^{k}(\psi)$. Let $\widetilde{h}$ be defined by $e^{\lambda h^{k}(\psi)}=\widetilde{h}^{k}(e^{\lambda \psi})$. By considering $I_{\Gamma}(u_{\lambda})(y)$ as functional $H_{x} \rightarrow H_{y}'$ and using Stoke's formula,
if $\widetilde{\widetilde{h}}$ is defined by $d \widetilde{\widetilde{h}}=\widetilde{h} d$, we see that
$$ \int_{\Gamma} e^{\lambda h^{k}(\psi)(x)}u_{\lambda}(x,z) d x= \int_{\Gamma} \widetilde{h}^{k}(e^{\lambda \psi(x)})u_{\lambda}(x,z) d x = \int_{\Gamma} e^{\lambda \psi(x)} \widetilde{\widetilde{h^{k}}}(u_{\lambda})(x,z) d x$$
If $k$ is selected so that $\widetilde{\widetilde{h^{k}}}$ is linear, over the "eigen-spaces" to this operator, $\psi$ can be selected as the real part of a holomorphic function. Also note that we can use that $d \widetilde{h}$ is reduced,
that is for $u_{\lambda}$ in the symbol-space, there is a local constant $c$ such that
$$ \int_{\Gamma} e^{\lambda \psi} \widetilde{\widetilde{h^{k}}}(u_{\lambda})d x=
c \int_{\Gamma} e^{\lambda \psi} \widetilde{h^{k}}(u_{\lambda})d x$$
Let $\Delta^{Q}_{x,(2)}$ be the weighted lineality corresponding to $h^{2}$.
We now claim that $\Delta^{Q}_{x,(2)} \subset \Delta^{Q}_{x,(1)}$, but the converse does not necessarily hold. Assume $h^{*}=h$ and let $\tau_{\eta}h=h\tau^{*}_{\eta}$. The condition that $\eta \in \Delta^{Q}_{x,(2)}$, means that $\tau_{\eta} h^{2} \tau_{\eta}=h^{2}$ $\Leftrightarrow$ $\tau_{\eta}^{2}h=h$, thus $\tau_{\eta}^{2}=\tau_{\eta}^{2 *}=1$. Thus, $\tau_{\eta} h \tau_{\eta} = h$. Conversely, $\tau_{\eta} h \tau_{\eta}=h$, means that $\tau^{2}_{\eta} h^{2}=h^{2}$, that should be compared with $\tau_{\eta}\tau_{\eta}^{*}h^{2}=h^{2}$ and if $\mbox{ Im }\eta \neq 0$, we do not necessarily have the opposite inclusion.

\newtheorem{d_w-lin}[d_inf-zero]{ Lemma }
\begin{d_w-lin}
$\Delta^{Q}_{x,(N)}$ is decreasing as $N$ increases.
\end{d_w-lin}

Assume $\Gamma$ is a circle-segment with center $\{ 0 \}$ in a complex variable. Let
$I_{\Gamma}(u_{\lambda})=\int_{\Gamma} e^{\varphi}u_{\lambda} dx$. The Cousin-integral $\Phi$ can now be given
as $\Phi=-\frac{\delta}{\delta y} I_{\Gamma}(u_{\lambda})=\int_{\Gamma} (\frac{\delta}{\delta x}e^{\varphi}) u_{\lambda} dx$
and $\frac{\delta}{\delta y} I_{\Gamma}(u_{\lambda})=-I_{\Gamma}(\frac{\delta}{\delta x}u_{\lambda})$, where
$y$ is assumed to lie in the same plane as $x$, thus $\Phi= \int_{\Gamma} e^{\varphi} \frac{\delta}{\delta x} u_{\lambda} dx$.
We know that $\Phi$ can be analytically continued along transversals $\ni y$, such that $\mbox{ Im }(x-y)=0$
for $x \in \Gamma$ (a real phase). $\widetilde{\Phi}$ is constructed according to Cousin, by adding to
$\Phi$ terms on the form $F_{\lambda}I_{\Gamma}(\frac{\delta}{\delta y}e^{-\varphi})/2 \pi i= F_{\lambda} \mbox{ log }(b-y)$,
where $b$ is the one end point of $\Gamma$, given in local coordinates. Analogously, for the other ($\pm$).
$\widetilde{\Phi}$ is then regular over the transversal, to  $y$.

\vspace*{.5cm}

Assume now instead that $\mbox{ Re }h \sim_{\infty} h$ and $\eta=(x-y)$, with
$\mbox{ Im }\eta \in \Delta({}^{t}I_{H})$. Given a disjoint decomposition of
$\Gamma$, we now have $t \eta \in \Sigma_{x,y}(\mbox{ Re }h)$, $\forall t$ real.
If we represent $\widetilde{\Phi}$ as a $I_{\Gamma}$-integral for $\varphi=h(\psi)$,
where $\psi$ can be chosen as real, we have that for $\varphi^{k}$, it is sufficient to consider
$I_{\Gamma}$ with a real phase. In the case with a complex variable, we thus have
$\Delta^{Q}_{x,(N)} \downarrow \{ 0 \}$ as $N \uparrow \infty$. For the case with several variables,
it follows immediately through the one-variable result, that $\{ 0 \}$ is represented as an isolated point in
$\Delta^{Q}_{x,(N)}$ for $N$ sufficiently large.

\vspace*{.5cm}

We can now repeat the argument concerning the set of lineality $\Delta_{\mathbf{C}}$. Define the ideal $J=I(\Delta^{Q}_{x,(1)})$. Let $V_{1}=U_{1} \backslash
\Delta^{Q}_{x,(1)} \subset N(J)=U_1=N($ $\mathcal{N}_{\textsl{U}_{\textit{1}}} )$. The construction can be iterated,
$J^{2}=I(\Delta^{Q}_{x,(2)})$, with $V_2=U_2 \backslash \Delta^{Q}_{x,(2)} \subset U_2$.
Since $\Delta^{Q}_{x,(N)} \downarrow \{ 0 \}$, as $N \uparrow \infty$, since $\Delta^{Q}_{x,(N)}$ is defined
by a homomorphism $(h \text{:} H \rightarrow H')$, we also have that $V_N$ is decreasing as $N \uparrow \infty$.
Thus, $U_{N}$ is decreasing with respect to inclusion of sets, $0 \notin V_{N}$, for all $N$, but $\varphi \in J^{N}_{0}$
for all $N$. We have $\mbox{ rad }J^{N}=\mathcal{N}$ ${}_{\textsl{U}_{\textit{N}}}$ for all $N$ and $\varphi \in J^{N}_{t}$ means $t=0$ and $\varphi^{N}$ corresponds to a real phase, that is $\varphi^{N}$ is reduced for weighted lineality.

\vspace*{.5cm}

For a generalization to the Cousin-integral in several variables, we can use the same argument. Thus, there is a positive $N_{0}$, such that $\Delta^{Q}_{x,(N)}=0$, $N \geq N_{0}$. Further, $\Delta^{Q}_{x,(j)}$ can, for any $j$, be realized by a locally injective $C^{\infty}$-mapping with a locally injective differential. The corresponding integral will take the form, for a positive constant $C$
$$\int_{\Gamma} \big( \frac{\delta}{\delta_{x_1}} \ldots \frac{\delta}{\delta_{x_k}} e^{\lambda \varphi(x_1, \ldots ,x_k)} \big) u_{\lambda}(x_1,\ldots,x_k,z) d x_1 \ldots d x_k=$$
$$C \int_{\Gamma} e^{\lambda \varphi(x_1, \ldots, x_k)}\frac{u_{\lambda}(x_1, \ldots, x_k,z)}{(x_1 - y_1) \ldots (x_k -y_k)} d x_1 \ldots d x_k$$

\section{ The phase ideal }

\subsection*{ The Hamilton-fields }
Assume $\Omega$ a convex set, in the sense that any straight line between two points in $\Omega$ is
contained in $\Omega$. We assume that the trajectory we consider, $\gamma$, has an interior that is
convex locally, so that in a neighborhood of a point on $\gamma$, in a small circle centered by a
point, the inner of $\gamma$ is situated on one side of the tangent to the point.

\vspace*{.5cm}

Consider the system
\begin{equation}  \frac{d \psi}{d \varphi}=P \qquad \frac{dh(\psi)}{d \varphi}=Q \nonumber \end{equation}
We then have $\frac{\delta \psi}{\delta \zeta_{j}}=P \frac{\delta \varphi}{\delta \zeta_{j}}$ and
$\frac{\delta h(\psi)}{\delta \zeta_{j}}=Q \frac{\delta \varphi}{\delta \zeta_{j}}$ and we consider the
Hamilton-fields $H_{h}=\sum_{j} \frac{d h}{d \varphi} \frac{\delta}{ \delta \zeta^{j}} - \frac{\delta h}{\delta \zeta_{j}}
\frac{d}{d \varphi}$. The Poisson-bracket can now be used to define an ideal over a domain
$\Omega$. Further, if $\Omega$ is a convex domain in the sense above and satisfies the two regularity conditions
and if $h \equiv 0$ on $\Omega$, then $e^{t H_{h}} \in \Omega$ for $\mid t \mid < t_{0}$

\vspace*{.5cm}

Assume $h$ constant over $\psi$, that is $\psi \in (I_{\mu})$, $h(\psi)=\mu \psi$ and $\Omega=N(I_{\mu})$, we then have
$\frac{d h(\psi)}{d \varphi}=0$ and in any singular point to the system, $H_{h}(\psi)=0$.
In this case $\frac{d}{d t}e^{t H_h}=0$ for $\mid t \mid < t_{0}$. If the
order of zero of the point to $P$ is infinite, $H_{h}(\psi)=0$ in a $\zeta$-neighborhood of the singular
point. Consider $h(\psi)=\mu \psi$ and assume $\frac{\delta h (\psi)}{\delta
\zeta_{j}}=\mu_{j} \frac{\delta P}{\delta \zeta_{j}}=\frac{Q}{P} \frac{\delta \psi}{\delta \zeta_{j}}$
and $\frac{d}{d \varphi}h=h\frac{d}{d \varphi}$. Then $\frac{\delta h}{\delta \zeta_{j}} \frac{d \psi}{d
\varphi}=\frac{\delta Q}{\delta \zeta_{j}}$ and $\frac{d h}{d \varphi}\frac{\delta
\psi}{\delta \zeta_{j}}=h(\frac{\delta P}{\delta \zeta_{j}})$. So, $H_{h}(\psi)=\sum_{j} h(\frac{\delta P}{\delta
\zeta_{j}})- \frac{\delta Q}{\delta \zeta_{j}}$. If we assume also $\frac{\delta Q}{\delta
\zeta_{j}}=\mu_{j} \frac{\delta P}{\delta \zeta_{j}}$, we have $H_{h}(\psi)=0$.

\vspace*{.5cm}

Assume $h$ self-adjoint and such that $h^{2}$ is locally injective and consider $\eta(\psi)
\psi=h(\psi)$. With these conditions $\eta$ is a real-valued function and
$$ \psi \frac{d \eta}{d \psi}=\frac{d h(\psi)}{d \psi}- \frac{h(\psi)}{\psi}$$
Assume further $\psi \frac{d \eta}{d \psi}=\mu- \frac{\overline{h(\psi)}}{\psi}$ for a real $\mu$.
Then, $\overline{\int_{\mid \psi \mid=1}h(\psi)d \psi}=\int_{\mid \psi \mid =1}\frac{d \eta}{d \psi}d
\psi$. Further $\overline{\int_{\mid \psi \mid=1} \mbox{ Im }h(\psi)d \psi}=0$ because since $h$ is
locally defined through a polynomial, we have locally $h^{*}=\overline{h}$. Particularly, if we
consider $e^{h(\varphi)}=\widetilde{h}(e^{\varphi})$ and consider the above relations with $h$
replaced by $\widetilde{h}$
$$ \mbox{ Im }\int_{\mid \psi \mid=1}e^{h(\varphi)}d \varphi=0$$
The conclusion is that if the ideal that defines the phase is finitely generated, the imaginary
part of the phase does not contribute to the corresponding operator.

\subsection{ Analytic continuation for $I_{\Gamma}$}
Assume $g=g_{0}+ic$, for a constant $c$ and $g_{0}$ real. If there exist $g_{j} < g$ such that
$\widetilde{h}(u) - v \in H^{loc}_{g_1}$, we say that $\widetilde{h}(u) \sim v$ in $H^{loc}_{g}$.
In the same way $\int_{\Gamma} e^{\lambda h(g)}u dx  \sim \int_{\Gamma} e^{\lambda g}v dx$. Assume
$h_{j}(g)=h(g_{j})$ and using continuity, $h_{j}(g) \rightarrow h(g)$ as $j \rightarrow \infty$,
further $\int_{\Gamma} e^{\lambda g} \widetilde{h_{j}}(u) d x \rightarrow \int_{\Gamma} e^{\lambda g} v dx$
in $(H^{loc}_{g})'$, modulo the equivalence above. Thus if we write $v=\widetilde{h}(u)$ we have
$$\int_{\Gamma} e^{\lambda h(g_{j} - g)} u dx \rightarrow 0 \text{ as } j \rightarrow \infty$$. This corresponds to
analytic continuation for $I_{\Gamma}$ over a real phase. Given lineality, we can chose $g_{j}$ non-constant
and real where $g_{j}$ corresponds to $\widetilde{I}_{\Gamma}$. Now consider
$$\widetilde{I}_{\Gamma}-I_{\Gamma}=\int_{\Gamma} e^{\lambda g}(\widetilde{h}_{j}(u)-\widetilde{h}(u)) dx= \int_{\Gamma} e^{\lambda (h_{j}(g) - h(g))} u
dx$$ and if $h_{j}(g)-h(g) < \varphi$ we get for $\widetilde{I}_{\Gamma} \sim I_{\Gamma}$, a remainder term that
gives an exponentially small contribution to the integral representation.
Trivially we have if $f \sim g$ in $H^{loc}_{\varphi}(\Omega)$ that there exists $g_{0} \in H(\Omega)$
such that $g \sim g_{0}$ and $f \sim_{m} g_{0}$. Further there exist $I^{0}_{\Gamma}$ such that
$\widetilde{I}_{\Gamma} - I^{0}_{\Gamma} \sim_{m} 0$.

\vspace*{.5cm}

Assume as before $\eta \in \Sigma_{x,y}$ implies
$T_{\eta}h=h$ for $T_{\eta}h=\tau_{\eta}h\tau_{\eta}$. Consider $\Sigma_{0}=\{ x \quad <t x,\eta> \leq 0 \quad \forall t \}$
or equivalently $\{x \quad <x,t \eta> \leq 0 \quad \forall t \}$. Let
$\widehat{F_{\eta}}=T_{\eta}-I$, where $F_{\eta}$ is considered in $H'(\Sigma_{0})$.
Further, $\Sigma^{0}=\{ \eta \quad T_{\eta}-I \mbox{ of type } 0 \}$. We see that if $\eta \in \Sigma_{x,y}$, then $\Sigma_{0}$
contains entire lines and we can in this way define a conic neighborhood of $\Sigma_{x,y}$ according to
$\Sigma^{0}$.

\vspace*{.5cm}

 We have noted that $\Delta^{Q}_{x} \downarrow \{ 0 \}$ as $N \uparrow \infty$ and we assume $h$
is algebraic over $I(\Delta^{Q}_{x})=\{ \varphi \quad \varphi=0 \text{ on } \Delta^{Q}_{x} \}$. If $h^{N}$ is
locally injective, we have also that $g$ is real. Let otherwise $L$ be such that $h(L)=I(\Delta^{Q}_{x})$, then
$g^{N} \in L$ and $g \in rad(L)$. Obviously $h^{N}$ is locally injective over $I(\Delta^{Q}_{x})$, for large iteration
indexes, which means that the corresponding $g$ defines a real phase. More precisely, define $(L_1)$ according to
$h(L_1)=I(\Delta_{1})$ and $(L_{2})$ according to $h^{2}(L_{2})=I(\Delta_{2}) $and so on. The proposition is thus
that if $(I)=(I)(0)$ is defined through
$(h,h^{2},h^{3},\ldots,h^{N-1})$ according to a chain of ideals $L_{1} \subset \ldots \subset L_{N-1}=(I)$, then $h^{N+1}$ is locally injective over $(I)$.

\vspace*{.5cm}

If $\varphi \in (I)$ is such that $\varphi=\sum_{j}c_{j}h^{j}(g_{j})$ for $ g_{j} \in L_{j}$ and since we are
dealing with geometric ideals, we can make the decomposition disjoint such that $\varphi=h(\sum_{j} d_{j} g_{j})$
$\exists g_{j} \in J_{j}$ and $(I)=\oplus_{j} J_{j}$. If $\varphi \in I(\Delta_{1})$ then $h(\varphi)=h(g^{2})$
for $g \in (L_{2})$. We now have that $N(h(I(\Delta_{1})))=N(h^{2}(L_{2}))=\Delta_{2} \subset \Delta_{1}$. Thus,
$h(\varphi)=0$ implies $\varphi=0$ and we have also that $h$ is locally injective on $(I)(\Delta_{1})$.

\vspace*{.5cm}

Note finally that if $R=\{ \varphi \quad \varphi \text{ real } \}$, we have $h \text{:} R \rightarrow R$
and $rad(R) \rightarrow rad(R)$ and further $rad(L_{1}) \subset h(L_{2})$. In the more general case,
$h$ is algebraic in the tangent-space and singular points are points for which $g$ is constant.
These constant phases can thus be approximated by real and non-constant ones.

\subsection{ The quantitative version of R\"uckert's theorem }

Assume the ideal $(I)$ defined over a domain of holomorphy $\Omega$ such that for all $\varphi \in (I)$, $h(\varphi)$ is in
the kernel of an elliptic operator. We then have that $e^{h(\varphi)}$ is a polynomial over $\Omega \backslash \mathcal{R}$.
Particularly, if $\widetilde{h}=h$ we have that the polar set to $h(\varphi)$ is an algebraic set.
Assume $h(e^{\varphi})=P/Q$ for $P,Q$ polynomials and $Q^{*}=Q$ hypoelliptic. We are going to prove that
$h(e^{\varphi})=V+\lambda Q$, for a polynomial $V$ in for instance $x_{1}$ and we make the approach
$\lambda Q=\widetilde{h}(e^{\varphi})$ and $$\lambda=\int_{\mid \xi \mid=r'} e^{h(\varphi)}\frac{1}{Q} d\xi$$
If the representation can be proven in the vicinity of the real infinity, then this implies assuming $P \prec Q$ that $e^{h(\varphi)} \prec Q$. Further
$\delta \lambda \prec \prec I$ given that $\delta P \prec \prec Q$. We will use a special case of a Lemma by H. Cartan
(cf. \cite{Cart}). Let $\Delta_{r}=\{ (x_{1}, \ldots, x_{r}) \quad \mid x_{j} \mid \leq r_{j} \quad \forall j \}$ and
$\Delta_{r}'=\{ (x_{2},\ldots, x_{r}) \quad$ $\mid x_{j} \mid \leq r_{j} \quad \forall j \}$ be polycylinders in the complex space. If $Q$ has zero's locally in $\Delta_{r}$ and if $f$ is holomorphic on $\Delta_{r}$ then there are $\lambda,V$ such that $f=\lambda Q + V$ where the polynomial in $x_{1}$
has degree $\mbox{ deg }Q-1$. If $\mid f \mid \leq 1$ on $\Delta_{r}$ we have that the coefficients for $V$
are holomorphic on $\Delta_{r}'$ and of modulus $\leq M$ for a positive constant $M$. Further $\mid \lambda \mid \leq M$
on $\Delta_{r}$. The definition of $\lambda$ can be proven independent of $r'$ close to $r$. We will briefly
sketch the proof. If $v_{p}$ are the coefficients to the polynomial $V$ we have that $\mid v_{p} \mid \leq c'$
on $\Delta_{r'}$ for a positive constant $c'$. Thus there exists a positive constant $c''$ such that
$\mid \lambda Q \mid \leq c''$ on $\Delta_{r}$. For $(x_{2},\ldots,x_{n}) \in \Delta_{r'}'$ and $\mid x_{1} \mid= r_{1}$ we have that $\mid Q \mid > c$ for a positive constant $c$ why using the maximum-principle we have that $\mid \lambda \mid \leq M$ on $\Delta_{r}$ for a positive $M$.

\vspace*{.5cm}

Consider again the polar set $(I)=\{ \varphi \quad e^{\varphi} \rightarrow 0 \quad x \rightarrow a \}$ and its closure with respect to Whitney-norm $\overline{(I)}$. Thus $\varphi \in \overline{(I)}$ $\Leftrightarrow \varphi \rightarrow 0 \quad ( \varphi \rightarrow const. ) \quad x \rightarrow a$ and we assume also that the limit is uniform in a neighborhood of the point $a$. We can now define a module of uniform approximation of $1$ (or of a constant) through non-constant functions $\psi=e^{\varphi}$ with $\varphi$ real. Denote this module $\mathcal{M}_{\textit{a}}$. It is evident in the minimally defined situation that
$\widehat{\mathcal{M}}_{\textit{a}} \subset J_{h}$, for a homomorphism $h$ and that $\overline{\mathcal{M}}_{\textit{a}}=\mathcal{M}_{\textit{a}}$  (Whitney's theorem).

\vspace*{.5cm}

Assume $\Omega$ a domain of holomorphy and $\Gamma=\mbox{ bd } \Omega$ a compact set, further that
$f_{1},\ldots,f_{p}$ are holomorphic functions on $\overline{\Omega}$, constituting a pseudo-base for a
geometric ideal $(I)$. We know that there is a module of coefficients, $c_{j} \in \mathcal{M}$, holomorphic on a compact polycylinder
$\Delta \subset \Omega$ such that $\sum_{j} c_{j}f_{j}=0$ on $\Delta$. We also assume the equation
$\sum_{j}d_{j}f_{j}=1$ solvable on $\Gamma$. If $w$ is an algebraic homomorphism $(I) \rightarrow (I)$
such that $w(f_{j})=constf_{j}$ for all $j$, then we obviously have $\mathcal{N}$ $=w^{-1}\mathcal{M}$ is a module such
that we can solve the equation $\sum_{j}a_{j}f_{j} \in J_{w}$, where $(J_{w})=\mbox{ ker }w$. If $w$ is reduced and $\mathcal{M}_{\textit{a}}$
is the module for uniform approximation of $1$, that is $c_{j} \rightarrow c_{j}^{0}$ such that
$\sum_{j} c_{j}^{0} f_{j} \rightarrow 1$ on $\Gamma$, then the approximating functions can be selected from
$\mathcal{N}$. If we instead consider a homomorphism $w \text{:} (I) \rightarrow (I)$ as factorized through
$\tau_{\epsilon}g$ for an algebraic homomorphism $g$, we can form $\mathcal{N}$ as before and solve the
homogeneous equation modulo monotropy. For the module of uniform approximation, we can consider the equations
$dw(\sum_{j} a_{j}f_{j})=0$ on $\Gamma$ and $w(\sum_{j} a_{j}f_{j})=0$ on $\Omega$. If $w$ is such that
$w^{2}$ is locally injective we have the minimally defined situation and we can always find $g \in (J_{w})(\Omega)$
such that $g \rightarrow g^{0}$ uniformly in a neighborhood of $\Gamma$. In the case where $w$ is an
algebraic homomorphism, it is sufficient to chose $\mathcal{N}=$ $(J_{w})(\Omega)$ and in the case of a factorized homomorphism we have the same proposition modulo monotropy.
\vspace*{.5cm}

We have earlier considered the example $\Omega=\Sigma$ and $\Gamma=\Delta^{Q}_{x}$ and the ideal
$(J)=\{ g \quad dh(g)=0 \}$ for a homomorphism $w$ such that $dw$ is reduced (locally injective) and
$h(g)=w(\widehat{g})-\widehat{w}(g)$. We write $a_{j}=e^{\varphi_{j}}$ and assume that
$w(a_{j})=\widehat{w}(\varphi_{j}) \rightarrow 0$ as $x \rightarrow a$ for $x \in V$, $V$ a neighborhood of
$a$. If $U$ is the dual to $V$ in $\mbox{Exp}$, we see that $w(\varphi_{j}) \in \mathcal{E'}$ ${}^{(0)}(U)$.
The problem is now to determine the support for these measures. We can also assume $\mid \widehat{w}(\varphi_{j})-Q_{j}(\frac{1}{x-a}) \mid < \epsilon$
for $x \neq a$, $x$ close to $a$. This means that the phase $w(\varphi_{j})$ can be seen locally as a
polynomial operator $Q_{j}(D)\delta_{A}$, where $A$ corresponds to $a$ through the duality. Under the condition $Df \prec \prec I$, it is sufficient to consider operators $Q$ of first order.

\vspace*{.5cm}

A comparison with Cartan's version of R\"uckert's Nullstellensatz gives with the same approach as earlier that if $h(e^{\varphi})=0$ we have modulo monotropy that the condition $\mid h(e^{\varphi}) \mid \leq \epsilon$ gives $\mid \lambda \mid \leq M \epsilon$ for a positive $M$ and a small positive $\epsilon$. If $(I)=\{ g \quad e^{g} \rightarrow 0 \quad x \rightarrow a \}$ and $g=h(\varphi)$ we can form $(I)_{m}=\{ h(\varphi) \quad \mid h(e^{\varphi}) \mid < \epsilon \}=\{ h(\varphi) \quad \mid \widehat{h}(\varphi) \mid < \epsilon \}$. We can now use Paley-Wiener's theorem modulo monotropy and we see that $h(\varphi) \in \mathcal{E'}$ ${}^{(0)}$, thus $(I)_{m}$ can be realized through measures.

\cite{Oka_2}
\bibliographystyle{amsplain}
\bibliography{ref}

\providecommand{\bysame}{\leavevmode\hbox to3em{\hrulefill}\thinspace}
\providecommand{\MR}{\relax\ifhmode\unskip\space\fi MR }
\providecommand{\MRhref}[2]{%
  \href{http://www.ams.org/mathscinet-getitem?mr=#1}{#2}
}
\providecommand{\href}[2]{#2}
\begin{thebibliography}{10}

\bibitem{Ben}
I.~Bendixson, \emph{Sur les courbes d{\'e}finies par des {\'e}quations
  diff{\'e}rentielles}, Acta Mathematica (1901).

\bibitem{Bour}
N.~Bourbaki, \emph{Vari{\'e}t{\'e}s differentielles et analytiques}, Hermann
  (1971).

\bibitem{BourL}
\bysame, \emph{Lie groups and {L}ie algebras}, Actualites sci. et industrielles
  Hermann (1975).

\bibitem{Cart}
H.~Cartan, \emph{Id{\'e}aux de fonctions analytiques de $n$ variables
  complexes}, Ann. Sci. de l'{\'E}cole Norm. Sup.,3 (61) (1944).

\bibitem{Chirka}
E.~M. Chirka, \emph{Complex {A}nalytic {S}ets}, Springer Science \& Business
  Media (1989).

\bibitem{Cous}
P.~Cousin, \emph{Sur les fonctions de $n$ variables complexes}, Acta Math., 19
  (1895).

\bibitem{jag_0}
T.~Dahn, \emph{Some remarks on {T}r{\`e}ves' conjecture.}, ArXiv (2013).

\bibitem{jag_I}
\bysame, \emph{On partially hypoelliptic operators. part i: Differential
  operators}, ArXiv (2015).

\bibitem{Gra}
H.~Grauert and R.~Remmert, \emph{Coherent {A}nalytic {S}heaves}, Grundlehren
  der mathematischen Wissenschaften, 265, Springer-Verlag (1984).

\bibitem{Han}
N.~Hanges, \emph{Propagation of singularities for a class of operators with
  double characteristics}, Inst. Adv. Stud., Princeton Univ. Press (1977/78).

\bibitem{Ho_LPDOII}
L.~H{\"o}rmander, \emph{The {A}nalysis of {L}inear {P}artial {D}ifferential
  {O}perators i,ii}, Springer-Verlag. (1983).

\bibitem{Horm}
\bysame, \emph{An {I}ntroduction to {C}omplex {A}nalysis in {S}everal
  {V}ariables}, North Holland Math. Library (1990).

\bibitem{Kohn}
J.J. Kohn, \emph{Hypoellipticity of {S}ome {D}egenerate {S}ubelliptic
  {O}perators}, Journal of functional analysis, 159,(1) (1998).

\bibitem{Mlg}
B.~Malgrange, \emph{Existence et approximation des solutions des {\'e}quations
  aux d{\'e}riv{\'e}s partielles et des {\'e}quations de convolution}, Ann.
  Inst. Fourier, 6 (1955-56).

\bibitem{Mart}
A.~Martineau, \emph{Sur les fonctionelles analytiques et la transformation de
  {F}ourier-{B}orel.}, Journal d'Analyse Math., Vol. XI. (1963).

\bibitem{Ni}
N.~Nilsson, \emph{Lecture notes on linear functional analysis.}, Manuscript,
  Lund (1992).

\bibitem{Ni_92}
T.~Nishino, \emph{Nouvelles recherches sur les fonctions enti{\`e}res de
  plusieurs variables complexes}, J. Math. Kyoto University, 8, (1) (1968).

\bibitem{Oka}
K.~Oka, \emph{Sur les fonctions analytiques de plusieurs variables},  (1960).

\bibitem{Oka_2}
\bysame, \emph{Une mode nouvelle engendrant les domaines pseudoconvexes}, Jap.
  J. Math., 32 (1962).

\bibitem{Palamadov}
V.~Palamodov, \emph{Linear {D}ifferential {O}perators with {C}onstant
  {C}oefficients.}, Grundlehren der mathematischen Wissenschaften, 168,
  Springer-Verlag (1970).

\bibitem{Sjo}
J.~Sj{\"o}strand, \emph{Singularit{\'e}s analytiques microlocales},
  Astr{\'e}risque, 95. (1982).

\end{thebibliography}
\end{flushleft}
\end{document}